\input amstex


\def\s2{\sqrt 2}

\def\eps{\epsilon}
\def\sech{\mathop{\roman {sech}}\nolimits}
\def\grad{\mathop{\bold \nabla\!}\nolimits}
\def\sgrad{\mathop{\bold \nabla_{\!\!s}}\nolimits}
\def\Int#1{\raise 1pt\hbox{$\scriptstyle\int$}\!{#1}\,}
\def\Int#1{\kern 1pt\raise .15pt\hbox{$\char'177$}\!{#1}\kern .3pt}

\def\Der_#1(#2){#2\!_{{}_{(#1)}}}

\def\FT{\mathop{\Cal F}\nolimits}
\def\IFT
   {\mathop{\lower .25pt\hbox{$\Cal I$}\kern -2.1pt%
       \hbox{${\Cal F}$}}\nolimits}

\def\CH{\mathop{\Cal C\!H}\nolimits}
\def\ICH
   {\mathop{{\Cal I}\kern -2.1pt\raise 0.15pt\hbox{${\Cal C\!H}$}}\nolimits}

\def\IST
   {\mathop{{\Cal I}\kern -2.1pt\raise 0.15pt\hbox{${\Cal S\!T}$}}\nolimits}

\def\sph{\hbox{\bf S}}
\def\reals{\hbox{\bf R}}
\def\Cx{\hbox{\bf C}}

\def \norm|#1|{\left \Vert#1\right\Vert}

\def\proof {\noindent{\it Proof.}\enspace}
\def\qedmark{\hbox{\vrule height 4pt width 3pt}}
\def\qedskip{\vrule height 4pt width 0pt depth 1pc}
\def\qed{\nobreak\quad\nobreak{\qedmark\qedskip}}
\def\HN{H_{\kern-1pt{}_N}}
\def\Hinfty{H_{\kern-1pt{}_\infty}}

\def\sphN{\sph_{\kern-1pt{}_N}}
\def\om{\omega}
\def\omN{\om_{\kern-1pt{}_N}}

\def\dsphN{{\hat{\sph}}_{\kern-1pt{}_N}}
\def\Sg#1{\Sigma_{{}_{#1}}}

\def \IP<#1,#2>{\left \langle #1,#2\right\rangle}
\def \PB(#1,#2){\left \lbrace #1,#2\right\rbrace}
\def \LD(#1?#2){{\Cal L}_{\!{}_{#1}}#2}
\def \D{\partial}
\def \C{\hbox{$\Cal C$}}
\def \L{\hbox{$\Cal L$}}
\def \S{\hbox{$\Cal S$}}
\def\halb{ {1\over 2}}

\def\SOn{\mathop{\hbox{\bf {SO}}}}
\def\GLn{\mathop{\hbox{\bf {GL}}}}
\def\SLn{\mathop{\hbox{\bf {SL}}}}
\def\Un{\mathop{\hbox{\bf {U}}}}
\def\SUn{\mathop{\hbox{\bf {SU}}}}

\def\gln{\mathop{\hbox{\bf {gl}}}}
\def\sln{\mathop{\hbox{\bf {sl}}}}

\def\sun{\mathop{\hbox{\bf {su}}}}

\def\es{\epsilon}

\def\a{{\bold a}}
\def\b{{\bold b}}
\def \mat(#1;#2;#3;#4)%
{\pmatrix{{#1}&{#2}\cr{#3}&{#4}}}
\def\G{\hbox{\bf G}}
\def\g{{\hbox{$\Cal G$}}}
\def\Ad{\mathop{\hbox{\rm Ad}}\nolimits}
\def\tr{\mathop{\hbox{\rm tr}}\nolimits}
\def\ad{\mathop{\hbox{\rm ad}}\nolimits}
\def\dip<<#1;#2>>{\hbox{$<\!<$}#1,#2\hbox{$>\!>\!$}}
\def\T{\hbox{\bf T}}
\def\t{\hbox{$\Cal{T}\,$}}
\def\tperp{\hbox{$\Cal T$}^\perp}
\def\diag(#1,#2){\hbox{{\rm diag}}(#1,#2)}
\def\Re{\mathop{\hbox{\rm Re}}\nolimits}
\def\Im{\mathop{\hbox{\rm Im}}\nolimits}

\def \mat(#1;#2;#3;#4)%
{\left(\matrix#1&#2\\#3&#4\endmatrix\right)}

\def \D{\partial}
\def \S{\hbox{$\Cal S$}}

\newcount\secnum \secnum=0       
\newcount\subsecnum              

\def\section#1{\advance\secnum by 1 \subsecnum=0%
            \head{\the\secnum. #1}\endhead }

\def\subsection#1{\advance\subsecnum by 1%
       \subhead\nofrills{\the\subsecnum. #1}\endsubhead}
\newcount\firstpageno
\newcount\tocpageno

\def\Preg{P_{\!{}_0}}
\def\simR{\mathrel{\hbox{$\sim{}_{\!\!\!\!{{}_R}}$}}}
\def \Dm{{\Cal D}_-}
\def \Fs{{\Cal F}_{\!\!\!\!{}_{\roman {scat}}}}
\def \IFs{{\Cal I \!F}_{\!\!\!\!{}_{\roman {scat}}}}
\def \Hp{{\Cal H}_+}
\def \Hm{{\Cal H}_-}
\def \Gp{{\Cal G}_+}
\def \Dmd{\Dm^{\!\!^{\raise 0.75pt\hbox{$\scriptscriptstyle\roman {disc}$}}}}
\def \Dmc{\Dm^{\!\!^{\raise 0.75pt\hbox{$\scriptscriptstyle\roman {cont}$}}}}


\documentstyle{amsppt}


\def\qedmark{\hbox{\vrule height 4pt width 3pt}}
\def\qedskip{\vrule height 4pt width 0pt depth 1pc}
\def\qed{\nobreak\quad\nobreak{\qedmark\qedskip}}

\def \mat(#1;#2;#3;#4)%
{\left(\matrix#1&#2\\#3&#4\endmatrix\right)}

\def \D{\partial}
\def \S{\hbox{$\Cal S$}}

\NoBlackBoxes

\hoffset = 52pt

\topmatter

\title The Symmetries of Solitons \endtitle

\author Richard S. Palais\endauthor

\address%
Department of Mathematics, %
Brandeis University, Waltham,
Massachusetts 02254\endaddress

\email palais\@math.brandeis.edu\endemail


\abstract 
In this article we will retrace one of the great mathematical adventures of this 
century---the discovery of the soliton and the gradual explanation of its remarkable
properties in terms of hidden symmetries. We will take an historical approach, starting 
with a famous numerical experiment carried out by Fermi, Pasta, and Ulam on one of 
the first electronic computers, and with Zabusky and Kruskal's insightful explanation 
of the surprising results of that experiment (and of a follow-up experiment of their own)
in terms of a new concept they called  ``solitons''. Solitons however raised even more 
questions than they answered. In particular, the evolution equations that govern solitons 
were found to be Hamiltonian and have infinitely many conserved quantities, pointing to the 
existence of many non-obvious symmetries.  We will cover next the elegant approach to solitons 
in terms of the Inverse Scattering Transform and Lax Pairs, and finally explain how those ideas 
led step-by-step to the discovery that Loop Groups, acting by ``Dressing Transformations'', 
give a conceptually satisfying explanation of the secret soliton symmetries. 
\endabstract

\thanks During the preparation of this paper, the author was supported in part by the 
Mathematics Institute and Sonderforschungsbereich 256 of Bonn University.\endthanks

\date June 9, 1997\enddate


\endtopmatter

\document

{
\parindent=0pt
\parfillskip=0pt
\leftskip=1.5pc \rightskip=\leftskip
\baselineskip=13pt

\def\tocsection{%
  \par
  \ifdim\lastskip > 2pt \else \removelastskip \vskip2pt \fi
  \def\tocbreak{F}%
  \noindent }

\def\tocline{%
  \par
  \if F\tocbreak \nobreak \else \def\tocbreak{T}\fi
  \noindent\kern30pt }

\def\thepagenum#1{%
  \tocpageno=\firstpageno
  \advance\tocpageno by-1
  \advance\tocpageno by#1
  \hfil
  \the\tocpageno
  \unskip\par }

\specialhead Contents\endspecialhead

\firstpageno = \pageno

\tocsection 1. Introduction \thepagenum{2}

\tocsection 2. Review of Classical Mechanics \thepagenum{4}
\tocline Newton's Equations \thepagenum{4}
\tocline The Lagrangian Viewpoint \thepagenum{5}
\tocline Noether's Principle \thepagenum{6}
\tocline The Hamiltonian Viewpoint \thepagenum{7}
\tocline Symplectic Manifolds \thepagenum{8}
\tocline Examples of Classical Mechanical Systems \thepagenum{12}
\tocline Physics Near Equilibrium \thepagenum{13}
\tocline Ergodicity and Thermalization \thepagenum{15}

\tocsection 3. Origins of Soliton Theory \thepagenum{17}
\tocline The Fermi-Pasta-Ulam Experiments \thepagenum{18}
\tocline The Kruskal-Zabusky Experiments \thepagenum{21}
\tocline A First Look at KdV \thepagenum{24}
\tocline ``Steepening'' and ``Breaking'' \thepagenum{26}
\tocline Dispersion \thepagenum{27}
\tocline Split-stepping KdV \thepagenum{28}
\tocline A Symplectic Structure for KdV \thepagenum{29}

\tocsection 4. The Inverse Scattering Method \thepagenum{32}
\tocline Lax Equations: KdV as an Isospectral Flow \thepagenum{33}
\tocline The Scattering Data and its Evolution \thepagenum{35}
\tocline The Inverse Scattering Transform \thepagenum{38}
\tocline An Explicit Formula for KdV Multi-Solitons \thepagenum{39}
\tocline The KdV Hierarchy \thepagenum{40}

\tocsection 5. The ZS-AKNS Scheme  \thepagenum{42}
\tocline Flat Connections and the Lax Equation, ZCC \thepagenum{42}
\tocline Some ZS-AKNS Examples\thepagenum{43}
\tocline The Uses of Solitons \thepagenum{44}
\tocline Nonlinear Schr\"odinger as a Hamiltonian Flow \thepagenum{47}
\tocline The Nonlinear Schr\"odinger Hierarchy \thepagenum{49}

\tocsection 6. ZS-AKNS Direct Scattering Theory \thepagenum{52}
\tocline Statements of Results \thepagenum{52}
\tocline Outline of Proofs \thepagenum{54}

\tocsection 7. Loop Groups, Dressing Actions, and Inverse Scattering \thepagenum{58}
\tocline Secret Sources of Soliton Symmetries \thepagenum{58}
\tocline Terng-Uhlenbeck Factoring and the Dressing Action \thepagenum{59}
\tocline The Inverse Scattering Transform \thepagenum{60}
\tocline ZS-AKNS Scattering Coordinates \thepagenum{60}

\tocsection References \thepagenum{63}

}

\section{Introduction}

\noindent
  In the past several decades, two major themes have dominated
developments in the theory of dynamical systems. On the one hand
there has been a remarkable and rapid development in the theory of
so-called ``chaotic'' systems, with a gradual clarification of the
nature and origins of the surprising properties from which these
systems get their name. Here what cries out to be explained is how
a system that is deterministic can nevertheless exhibit behavior
that appears erratic and unpredictable.

  In this article I will be discussing a second class of
systems---equally  puzzling, but for almost the opposite reason.
For these so-called ``integrable systems'', the challenge
is to explain the striking predictability, regularities, and
quasi-periodicities exhibited by their solutions, a behavior
particularly apparent for a special class of solutions, called
``solitons''. The latter exhibit a ``particle-like'' behavior
that gives them their name; for example they have geometric shapes
that show a remarkable degree of survivability under conditions
that one might normally expect to destroy such features.

  Such conservation of geometric features is known to be
intimately bound up with notions of symmetry---in fact,
when suitably formalized, a famous theorem of E. Noether
states that conserved quantities correspond to one-parameter
groups of automorphisms of the dynamical system---and therein
lies a puzzle. These systems do not have manifestly obvious
symmetries to account for these anomalous conservation laws,
and to fully understand their surprising behavior we must
search for the secret sources of their hidden symmetries.
This article  will be about that search, and about the
many mathematical treasures it has so far revealed.

  A major problem for anyone attempting an exposition of ``soliton
mathematics'' or ``integrable systems'' is the vast extent of its
literature. The theory had its origins in the 1960's, and so can
be considered relatively recent. But early research in the subject
revealed mysterious new mathematical phenomena that quickly
attracted the attention and stimulated the curiosity of many
mathematicians throughout the world. As these researchers took up
the intriguing challenge of understanding these new phenomena, an
initial trickle of papers soon grew to a torrent, and the eventual
working out of the details of the theory resulted from a concerted
effort by hundreds of mathematicians whose results are spread
over a still growing bibliography of many thousands of papers.

  Attempting to cover the subject in sufficient detail to mention 
all these con\-trib\-utions---or even most of the important
contributions---would require hundreds of pages. I have neither 
the time nor the expertise to undertake such a task, and instead 
I have tried to provide a guided tour through what I consider some
of the major highlights of the subject. But the reader should realize 
that any attempt to compress such a massive subject in so few pages must 
be an exercise in selectivity that will in large measure reflect personal 
taste and biases of the author rather than some objective measure of 
importance.

   Another disclaimer: as we proceed I will try to present some
of the remarkable story of how the subject began and developed.
I say ``story'' rather than ``history'' because my report will be
anecdotal in nature. I will try to be accurate, but I do not pretend 
to have done careful historical research. It is particularly important 
to keep in mind that during most of the development of the theory
of integrable systems there was a very large and active group
of mathematicians working on the subject in the former Soviet
Union. Since communication of results between this group and
the group of western mathematicians working in the field was
slower than that within each group, even more than usual there
were frequent cases in which similar advances were made nearly
simultaneously in one group and the other. Statements made in
this article to the effect that some person discovered a certain
fact should not be interpreted as claiming that person had
priority or sole priority in the discovery.

  There have been a number of fine volumes written
that make a serious effort to encompass the bulk of soliton
theory, giving careful historical and bibliographic references.
I hope my abbreviated account will stimulate readers to consult
these more complete sources, several of which are listed in the
references ([AC], [FT], [N], [NMPZ]).

   The organization of this article will be in part historical.
We will start with some surprising numerical experiments of
Fermi-Pasta-Ulam and of Zabusky-Kruskal that were the origins
of soliton theory. We will next consider the remarkable Inverse 
Scattering Transform and the related concept of Lax Pairs, first 
in the original context of the Korteweg-de Vries (KdV) equation, 
and then for the more general hierarchies of integrable systems 
introduced by Zakharov and Shabat and by Ablowitz, Kaup, Newell, and 
Segur (ZS-AKNS). We will trace how developments that grew out of the
ZS-AKNS approach eventually led to a synthesis that explains most of 
the phenomena of soliton theory from a unified viewpoint. In particular, 
it uncovers the source of the hidden symmetries of solitons, explaining 
both the existence of so many commuting constants of the motion and also
the characteristic phenomenon of B\"acklund Transformations. This
synthesis had its origins in the idea of ``dressing transformations'',
and in explaining it I will follow the recent approach of Chuu-lian
Terng and Karen Uhlenbeck. I would like to express my sincere thanks to 
Chuu-lian for putting up with my countless requests that she interrupt 
her own work in order to explain to me some detail of this approach. 
Without these many hours of help, it would not have been possible 
for me to complete this article.

  This article is a revised version of notes from a series of 
Rudolf-Lipschitz Lectures that I delivered at Bonn University in 
January and February of 1997. I would like to thank the Mathematisches
Institut of Universit\"at Bonn and its Sonderforschungsbereich 256 for 
honoring me with the invitation to give that lecture series, and to 
thank the lively and high-level audience who, by their interest, 
stimulated me to write up my rough notes. 

  My thanks to Bob Palais for pointing out a problem in my original
discussion of split-stepping---and for helping me to re-write it

  And special thanks to barbara n beeton for an exceptional job of
proof-reading. The many changes she suggested have substantially
improved readability.

\section {Review of Classical Mechanics}

\noindent
In this section we will review Classical Mechanics, in both
the Lagrangian and Hamiltonian formulations. This is intended
mainly to establish notational conventions, not as an exposition
for novices. We shall also review the basic geometry of symplectic
manifolds.

\subsection {Newton's Equations}

\noindent
  Let $\C$ be a Riemannian manifold (``configuration space'') and
$\Pi: T\C \to \C$ its tangent bundle.
 A vector field $X$ on $T\C$ is called a second order ODE on
$\C$ if $D\Pi(X_v) = v$ for all $v$ in $T\C$. If $\gamma$ is a
solution curve of $X$ and $\sigma = \Pi(\gamma)$ is its projection
onto $\C$ then, by the chain rule,
$\sigma^\prime(t) = D\Pi(\gamma^\prime(t)) = D\Pi(X_{\gamma(t)}) =
\gamma(t)$,
i.e., $\gamma$ is
the velocity field of its projection. An easy argument shows
conversely that if this is true for all solutions of a vector
field $X$ on $T\C$ then $X$ is a second order ODE on $\C$.
For this reason we shall say that a smooth curve $\sigma(t)$
in $\C$ satisfies the second order ODE $X$ if $\sigma^\prime$
is a solution curve of $X$.

   Given coordinates $x_1,\ldots,x_n$ for $\C$ in $O$, we define
associated ``canonical'' coordinates
$q_1,\ldots,q_n ,\dot q_1,\ldots,\dot q_n$ in $\Pi^{-1}(O)$
by $q_i = x_i\circ\Pi$ and $\dot q_i = dx_i$.
Let $\sigma :[a,b] \to \C$
be a smooth curve in $\C$, $\sigma^\prime : [a,b] \to T\C$
its velocity. If we define $x_i(t) = x_i(\sigma(t))$ and
$q_i(t) = q_i(\sigma^\prime(t)) = x_i(t)$, then
$\dot q_i(t) := \dot q_i(\sigma^\prime(t)) = dx_i(\sigma^\prime(t))
= {dx_i(t)\over dt} = {dq_i(t)\over dt}$.
It follows that a vector field $X$ on $\C$ is a second order ODE
if and only if in each canonical coordinate system it has the form
$X = \sum_i(\dot q_i{\D/\D q_i} + F(q_i,\dot q_i){\D/\D \dot q_i})$,
or equivalently the condition for $\sigma^\prime$ to be a
solution of $X$ is that $dq_i(t)/dt = \dot q_i(t)$,
$d \dot q_i(t)/dt = F_i(q_i(t),\dot q_i(t))$ (so
$d^2 x_i(t)/dt^2 = F_i(x_i(t), d x_i(t)/dt)$, explaining why
it is called a second order ODE)\null.

  The classic example of a second order ODE on $\C$ is the
vector field $X$ generating the geodesic flow on $T\C$---for
each $v$ in $T\C$ the solution curve of $X$ with initial
condition $v$ is $\sigma^\prime$ where $\sigma(t) = \exp(tv)$
is the unique geodesic on $\C$ with $\sigma^\prime(0) = v$.
In local coordinates, $x_i(\sigma(t))$ satisfy the system:
$${d^2 x_i\over dt^2} =
   -\Gamma^i_{jk}(x) {dx_j\over dt} {dx_k\over dt}$$
(where the $\Gamma^i_{jk}$ are the Christoffel symbols).
What we shall call {\it Newton's Equations\/} (NE) for $\C$
is a second order ODE $X^U$ for $\C$ that is a slight
generalization of the geodesic flow and is determined by a
smooth real-valued function $U$ on $\C$ called the potential
energy function:
$${d^2 x_i\over dt^2} =
   -\Gamma^i_{jk}(x) {dx_j\over dt} {dx_k\over dt}
   -{\D U\over \D x_i}.\leqno{(NE)}$$
[Here is an intrinsic, geometric description of (NE).
The gradient of $U$, $\grad U$ is a vector field on
$\C$, and we call $-\grad U$ the force. If $\sigma(t)$ is any
smooth curve in $\C$, and $v(t)$ is any tangent vector field along
$\sigma$ (i.e., a lifting of $\sigma$ to $T\C$), then the
Levi-Civita connection allows us to covariantly differentiate $v$
along $\sigma$ to produce another vector field $Dv/dt$ along
$\sigma$.  In particular, if for $v$ we take the velocity field
$\sigma^\prime(t)$, we can interpret $D\sigma^\prime/dt$ as the
acceleration of $\sigma$, and the curve $\sigma$ satisfies
Newton's Equations (for the potential $U$) if and only if
$D\sigma^\prime/dt = -\grad U$.]

\subsection {The Lagrangian Viewpoint}

\noindent
  We define the kinetic energy function $K$ on $T\C$ by
$K(v) = {1\over2}\norm|v|^2$, and we also consider the potential
energy as a function on $T\C$ by $U(v) = U(\Pi(v))$. Their
difference $\L = K - U$ is called the Lagrangian function
on $T\C$, and if $\sigma :[a,b] \to \C$ is any smooth curve
in $\C$ we define its action
$A(\sigma)=\int_a^b \L(\sigma^\prime(t))\,dt$.
In canonical coordinates as above,
$\L(q,\dot q) = \halb\sum_{ij}g_{ij} \dot q_i \dot q_j - U(q)$,
so if we write $x_i(t) = x_i(\sigma(t))$, then
$q_i(\sigma^\prime(t)) = x_i(t)$,
$\dot q_i(\sigma^\prime(t)) = dx_i/dt$, and therefore
$$A(\sigma)= \int_a^b \L(q(t),\dot q(t))\, dt=
  \int_a^b \halb \sum_{ij}g_{ij}(x(t))
   {dx_i\over dt} {dx_j\over dt} - U(x(t)) \,dt.$$

  Let $\sigma_\es:[a,b] \to \C$ be a smooth one-parameter family of
curves defined for $\es$ near zero, and with $\sigma_{{}_0} = \sigma$.
If we define $\delta\sigma = ({d\over d\es})_{{}_{\es =0}} \sigma_\es$
(a vector field along $\sigma$) then it is easy to see that
$({d\over d\es})_{{}_{\es =0}} A(\sigma_\es)$ depends only on $\sigma$
and $\delta\sigma$, and we denote it by $DA_\sigma(\delta\sigma)$.
Define
$q_i(t,\es) = q_i(\sigma_\es^\prime(t)) = x_i(\sigma_\es(t))$,
$\delta q_i(t) = \D q_i(t,0)/\D \es$,
$\dot q_i(t,\es) = \dot q_i(\sigma_\es^\prime(t))$
and
$\delta \dot q_i(t) = \D \dot q_i(t,0)/\D \es$.
Then clearly
$\dot q_i(t,\es) = \D q_i(t,\es)/\D t$,
so, by equality of cross derivatives,
$\delta \dot q_i(t) = {d\over dt} \delta q_i$.

  It is now easy to compute $DA_\sigma(\delta\sigma)$. In fact,
differentiating under the integral sign, using the chain rule, and
integrating by parts gives:

$$\eqalign{
   DA_\sigma(\delta\sigma) &=
  \int_a^b \sum_i\left(
   {\D\L\over\D q_i} \delta q_i +
   {\D\L\over\D \dot q_i} \delta \dot q_i\right)\,dt\cr
                            &=
   \int_a^b \sum_i\left(
   {\D\L\over\D q_i}  -
   {d\over dt}{\D\L\over\D \dot q_i} \right)\,\delta q_i\,dt
    + \left[ \sum_i{\D\L\over\D \dot q_i} \delta q_i  \right]_a^b\cr
                           &=
   \int_a^b \sum_i\left(
   {\D\L\over\D q_i}  -
   {d\over dt}{\D\L\over\D \dot q_i} \right)\,\delta q_i\,dt
    + \left[ \IP<\sigma^\prime(t),\delta\sigma(t)>  \right]_a^b\cr
}$$

  The curve $\sigma$ is called a critical point of the action
functional $A$ if $DA_\sigma(\delta\sigma)$ vanishes for all
variations $\delta\sigma$ vanishing at the endpoints $a$ and $b$,
or equivalently if the Euler-Lagrange equations
${\D\L\over\D q_i} - {d\over dt}{\D\L\over\D \dot q_i} = 0$
are satisfied. Substituting in the expression for $\L(q,\dot q)$
above, and recalling the definition of the Christoffel symbols,
it is easy to check that $\sigma$ is a critical point of the
action functional if and only if it satisfies Newton's Equations.

  It follows that if $\sigma$ is a solution of Newton's Equations
then for any variation $\delta \sigma$, not necessarily vanishing
at the endpoints,
$$ DA_\sigma(\delta\sigma) =
  \left[\, \IP<\sigma^\prime(t),\delta\sigma(t)> \, \right]_a^b.$$

   As a first application, consider the variation of $\sigma$
defined by $\sigma_\es(t) = \sigma(t+\es)$. Clearly
$\delta\sigma(t) = \sigma^\prime(t)$ and
$A(\sigma_\es) = \int_{a+\es}^{b+\es}\L(\sigma^\prime)\, dt$,
so the definition of $ DA_\sigma(\delta\sigma)$ gives
$DA_\sigma(\delta\sigma) = [\L(\sigma^\prime(t))]_a^b$,
while the above general formula for
$DA_\sigma(\delta\sigma)$ when $\sigma$ satisfies (NE) gives
$DA_\sigma(\delta\sigma)=
[\,\norm|\sigma^\prime(t)|^2\,]_a^b = [2K(\sigma^\prime(t))]_a^b$.

  If we define the {\it Hamiltonian\/} or total energy function $H$
on $T\C$ by $H = 2K - \L = 2K - (K - U) = K + U$, then it follows
that [$H(\sigma^\prime)]_a^b = 0$, or in other words $H$ is
constant along $\sigma^\prime$ whenever $\sigma$ is a solution of
Newton's Equations. Now a function $F$ on $T\C$ that is constant along
$\sigma^\prime$  whenever $\sigma: [a,b] \to \C$ satisfies (NE) is called
a constant of the motion for Newton's Equations, so we have proved:
\proclaim {Conservation of Energy Theorem} The Hamiltonian $H = K + U$
is a constant of the motion for Newton's Equations.\endproclaim

\noindent
[Here is a more direct proof.
$K(\sigma^\prime) = {1\over 2}g(\sigma^\prime,\sigma^\prime)$,
where $g$ is the metric tensor. By definition of the Levi-Civita
connection, $Dg/dt =0$, and (NE) says
$D\sigma^\prime/dt = -\grad U$, so
$d K(\sigma^\prime)/dt = g(-\grad U,\sigma^\prime) = -dU/dt$.]

\subsection {Noether's Principle}

\noindent
   A diffeomorphism $\phi$ of $\C$ induces a diffeomorphism
$D\phi$ of $T\C$, and we call $\phi$ a symmetry of Newton's
Equations if $D\phi$ preserves $\L$, i.e., if $\L\circ D\phi = \L$.
In particular, any isometry of $\C$ that preserves $U$ is a
symmetry of (NE). We note that if $\phi$ is a symmetry of (NE)
and $\sigma$ is any smooth path in $\C$ then
$A(\phi\circ\sigma) = A(\sigma)$, and it follows that $\phi$
permutes the critical points of $A$. Thus if $\sigma$ is a solution
of (NE) then so is $\phi\circ\sigma$. A vector field $Y$ is called an
infinitesimal symmetry of Newton's equations if it generates a
one-parameter group of symmetries of Newton's equations, so in
particular any Killing vector field that is tangent to the level
surfaces of $U$ is an infinitesimal symmetry of Newton's Equations.

   Suppose that $Y$ is any vector field on $\C$ generating a one-parameter
group of diffeomorphisms $\phi_t$ of $\C$. We associate to $Y$ a function
$\hat Y$ on $T\C$, called its conjugate momentum function, by
$\hat Y(v) = \IP<v,Y_{{}_{\Pi(v)}}>$. If $\sigma$ is any smooth
path in $\C$, then we can  generate a variation of $\sigma$ defined by
$\sigma_\es(t) = \phi_\es(\sigma(t))$. Then by definition,
$\delta\sigma(t) = Y_{\sigma(t)}$ so, by the above general formula,
if $\sigma$ is a solution of Newton's Equations then
$DA_\sigma(\delta\sigma) = [\hat Y(\sigma^\prime(t))\,]_a^b$
Now suppose $Y$ is an infinitesimal symmetry of Newton's Equations.
Then since $A(\sigma_\es) = A(\phi_\es\circ\sigma) = A(\sigma)$,
$DA_\sigma(\delta\sigma)$ is zero by definition, hence
$[\hat Y(\sigma^\prime(t))\,]_a^b = 0$, i.e., $\hat Y$ is constant
along $\sigma^\prime$. This proves:
\proclaim {E. Noether's Principle} The conjugate momentum of an
infinitesimal symmetry is a constant of the motion.\endproclaim

  The conjugate momentum to the vector field $\D/\D q_i$ is denoted
by $P_i$; $P_i = \sum_j g_{ij}\dot q_j = {\D \L \over \D \dot q_i}$,
and it follows from the non-degeneracy of the inner-product that
we can use
$q_1,\ldots,q_n ,P_1,\ldots,P_n$ as coordinates in $\Pi^{-1}(O)$.
The fact that Newton's Equations are equivalent to the Euler-Lagrange
equations says that in these coordinates
Newton's Equations take the form:
${dq_i \over dt} = \dot q_i$, ${dP_i\over dt} = {\D\L\over \D q_i}$
(i.e , $X^U =
\sum_i(\dot q_i {\D \over \D q_i}+ {\D\L\over \D q_i}{\D \over \D P_i})$ ).
Since $\sum_i P_i \dot q_i = 2K$,
$H = \sum_i P_i \dot q_i - \L$, so
$dH = \sum_i
(\dot q_i dP_i + P_i d\dot q_i
  - {\D \L\over \D q_i}dq_i
   - {\D \L\over \D \dot q_i}d\dot q_i) =
   \dot q_i dP_i - {\D \L\over \D q_i}dq_i$,
or in other words, ${\D H\over \D q_i} = -{\D \L\over \D q_i}$
and ${\D H\over \D P_i} = \dot q_i$. Thus Newton's Equations take
the very simple and symmetric form (called Hamilton's Equations)
${dq_i \over dt} = {\D H\over \D P_i}$,
${dP_i\over dt} = -{\D H\over \D q_i}$.
Equivalently, the vector field $X^U$ has the form
$X^U =
\sum_i({\D H\over \D P_i} {\D \over \D q_i}-
{\D H\over \D q_i}{\D \over \D P_i})$.

\goodbreak
\subsection {The Hamiltonian Viewpoint}

\noindent
  So far we have looked at the dynamics of Newton's Equations
on the tangent bundle $T\C$ of the configuration space. We will
refer to this as the Lagrangian viewpoint. Since $\C$ is Riemannian,
there is a canonical bundle isomorphism $L : T\C \to T^*\C$
of $T\C$ with the cotangent bundle, which in this setting is
called the Legendre transformation. Explicitly, $L(v)(u) = \IP<u,v>$.
The Hamiltonian  viewpoint towards particle mechanics consists
in moving the dynamics over to $T^*\C$ via the Legendre
transformation. Remarkably, the transferred dynamics preserves
the natural symplectic structure on $T^*\C$, and this fact is
the basis for powerful tools for better analyzing the situation.
The functions $\L\circ L^{-1}$ and $H\circ L^{-1}$ are still
called the Lagrangian and Hamiltonian function respectively and
will still be denoted by $\L$ and $H$. By further such abuse of
notation we will denote the vector field $DL(X^U)$ on $T^*\C$ by $X^U$.

  Just as with the tangent bundle, coordinates $x_1,\ldots, x_n$
for $\C$ in $O$ define natural coordinates
$q_1,\ldots, q_n,p_1,\ldots, p_n$
for the cotangent bundle in $\Pi^{-1}O$. Namely, $q_i = x_i \circ \Pi$
as before, while the $p_i$ are defined by
$p_i(\ell) = \ell(\D/\D x_i)$. It is immediate from the definitions
that $q_i \circ L = q_i$ while $p_i \circ L = P_i$, so it follows from
the calculation above that the vector field $X^U$ (i.e., $DL(X^U)$)
on $T^*\C$ describing the dynamics of Newton's Equations is
$X^U =
\sum_i({\D H\over \D p_i} {\D \over \D q_i}-
{\D H\over \D q_i}{\D \over \D p_i})$

  There is a natural $1$-form $\omega$ on $T^*\C$; namely if
$\ell$ is a cotangent vector of $\C$, then $\omega_\ell = D\Pi^*(\ell)$,
or in other words, for $Y$ a tangent vector to $T^*C$ at $\ell$,
$\omega_\ell(Y) = \ell(D\Pi(Y))$, where $\Pi: T^*\C \to \C$ is the
bundle projection. (We note that $\omega$ does {\it not\/} involve
the Riemannian metric, and in fact is natural in the sense that if
$\phi$ is any diffeomorphism of $\C$ 
and $\Phi = (D\phi)^*$ is the induced diffeomorphism of $T^*\C$
then $\Phi^*(\omega) = \omega$.)
We define the natural $2$-form $\Omega$ on $T^*\C$ by
$\Omega = d\omega$, so $\Omega$ is exact and hence closed,
i.e., $d\Omega = 0$.

 It is then easy to check that
$\omega = \sum_i p_i \, dq_i$ and hence
$\Omega = \sum_i dp_i\wedge dq_i$.
An immediate consequence of this is that $\Omega$ is non-degenerate,
i.e., the map $v \mapsto i_v\Omega$ is an isomorphism of the
tangent bundle of $T^*\C$ with its cotangent bundle. (Here
$i_v\Omega(u) = \Omega(v,u)$.) In fact, if
$v = \sum_i (A_i {\D \over \D q_i} + B_i {\D \over \D p_i})$ then
$i_v\Omega = \sum_i (A_i dp_i - B_i dq_i )$.
In particular
$i_{{}_{X^U}}\Omega = \sum_i
({\D H\over \D p_i} dp_i + {\D H\over \D q_i} dq_i ) = dH$.

   {\it Any\/} coordinates $q_1,\ldots, q_n,p_1,\ldots, p_n$
for $T^*\C$ are called ``canonical coordinates'' provided
$\Omega = \sum_i dp_i\wedge dq_i$. It follows that the
``equations of motion'' for solutions of Newton's Equations take
the Hamiltonian form:
${dp_i \over dt} = -{\D H\over \D q_i}$,
${dq_i \over dt} = {\D H\over \D p_i}$, for any such coordinates.
If $H$ happens not to involve a particular $q_i$ explicitly,
i.e., if $H$ is invariant under the one parameter group of
translations $q_i \mapsto q_i + \es$, then this $q_i$ is called a
{\it cyclic variable\/}, and its ``conjugate momentum'' $p_i$
is clearly a constant of the motion since
${dp_i\over dt} = -{\D H\over\D q_i} = 0$.
If we can find canonical coordinates
$q_1,\ldots, q_n,p_1,\ldots, p_n$
such that {\it all\/} of the $q_i$ are cyclic then we call
these variables {\it action-angle variables\/}, and when
such coordinates exist we say that the Hamiltonian system
is {\it completely integrable\/}. The solutions of a
completely integrable system are very easy to describe in
action-angle variables. Note that we have $H = H(p_1,\ldots,p_n)$.
For each $c$ in $\reals^n$ we have a submanifold
$\Sigma_c = \{\ell \in T^*\C \mid p_i(\ell) = c_i\}$, and since the
$p_i$ are all constants of the motion, these are invariant submanifolds
of the flow. Moreover these submanifolds foliate $T^*\C$, and on each
of them $q_1,\ldots, q_n$ are local coordinates. If we define
$\omega_i(c) = {\partial H \over \partial p_i}(c)$, then on $\Sigma_c$
Hamilton's Equations reduce to ${d q_i \over dt} = \omega_i(c)$, so
on $\Sigma_c$ the coordinates $q_i(t)$ of a solution curve are given
by $q_i(t) = q_i(0) + \omega_i(c) t$. Frequently the surfaces
$\Sigma_c$ are compact, in which case it is easy to show that each
connected component must be an $n$-dimensional torus. Moreover in
practice we can usually determine the $q_i$ to be the angular
coordinates for the $n$ circles whose product defines the torus
structure---which helps explain the terminology action-angle variables.

 Later we will look in more detail at the problem of determining
whether a Hamiltonian system is completely integrable.

\subsection {Symplectic Manifolds}

\noindent
  The cotangent bundle of a manifold is the model for what is
called a symplectic manifold. Namely, a symplectic manifold is a
smooth manifold $P$ together with a closed non-degenerate $2$-form
$\Omega$ on $P$. If $F: P \to \reals$ is a smooth real-valued function
on $P$ then there is a uniquely determined vector field $X$ on $P$
such that $i_X \Omega = dF$, and we call $X$ the symplectic gradient of $F$
and denote it by $\sgrad F$. Thus we can state our observation above
by saying that the vector field $X^U$ on $T^*\C$ is the symplectic
gradient of the Hamiltonian function: $X^U = \sgrad H$.

  By an important theorem of Darboux, ([Ar], Chapter 8) in the 
neighborhood of any point of $P$ there exist ``canonical coordinates''
$q_1,\ldots,q_n,p_1,\ldots,p_n$ in which $\Omega$ has the form
$\sum_i dp_i\wedge dq_i$, and in these coordinates
$\sgrad H = \sum_i({\D H\over \D p_i} {\D \over \D q_i}-
{\D H\over \D q_i}{\D \over \D p_i})$, or equivalently the
solution curves of $\sgrad H$ satisfy Hamilton's equations
${dp_i \over dt} = -{\D H\over \D q_i}$,
${dq_i \over dt} = {\D H\over \D p_i}$.

  Before considering Poisson brackets on symplectic manifolds,
we first make a short digression to review Lie derivatives.
Recall that if $X$ is a smooth vector field on a smooth manifold $M$,
generating a flow $\phi_t$, and if $T$ is any smooth tensor
field on $M$, then the Lie derivative of $T$ with respect to $X$
is the tensor field $\LD(X?T) = {d\over dt}|_{{}_{t=0}} \phi_t^*(T)$.
If $\LD(X?T) = 0$ then we shall say  that ``$X$ preserves $T$'', for
this is the necessary and sufficient condition that the flow
$\phi_t$ preserve $T$, i.e., that $\phi_t^*(T) = T$ for all $t$.
There is a famous formula of Cartan for the Lie derivative operator
$\LD(X?)$ restricted to differential forms, identifying it with
the anti-commutator of the exterior derivative operator $d$
and the interior product operator $i_X$:
$$\LD(X?) = d i_X + i_X d.$$
If $\theta$ is a closed $p$-form this gives
$\LD(X?\theta) = d(i_X \theta)$, so $X$ preserves $\theta$
if and only if the $(p-1)$-form $i_X \theta$ is closed.
In particular this demonstrates the important fact
that a vector field $X$ on a symplectic manifold $P$
is symplectic (i.e., preserves the symplectic form, $\Omega$)
if and only if $i_X\Omega$ is a closed $1$-form
(and hence, at least locally, the differential of a smooth function).
The well known identity $\LD([X,Y]?) = [\LD(X?),\LD(Y?)]$ implies
that the space of symplectic vector fields on $P$ is a Lie algebra,
which we can think of as the Lie algebra of the group of symplectic
diffeomorphisms of $P$. It is an interesting and useful fact that
the space of Hamiltonian vector fields on $P$, i.e., those for which
$i_X\Omega$ is an exact form, $dF$, is not only a linear  subspace, but
is even a Lie subalgebra of the  symplectic vector fields, and moreover
the commutator subalgebra of the symplectic vector fields is included
in the Hamiltonian vector fields. To demonstrate this we shall show that
if $i_X\Omega$ and $i_Y\Omega$ are closed forms, then
$i_{[X,Y]} \Omega$ is not only closed but even exact, and in fact it
is the differential of the function $\Omega(Y,X)$. First, using the fact
that Lie derivation satisfies a Leibnitz formula with respect to any
natural bilinear operation on tensors (so in particular with
respect to the interior product),
$\LD(X?(i_Y\Omega)) = i_{{}_{(\LD(X?Y))}}\Omega + i_Y (\LD(X?\Omega))$.
Thus, since $\LD(X?Y) = [X,Y]$ and $\LD(X?\Omega)=0$,
$\LD(X?(i_Y\Omega)) = i_{{}_{[X,Y]}}\Omega$. Finally, since
$d (i_Y\Omega) = 0$, Cartan's formula for $\LD(X?(i_Y\Omega))$ gives
$i_{{}_{[X,Y]}}\Omega = d i_X(i_Y\Omega) = d(\Omega(Y,X))$.

\noindent
{\bf Remark.} It is possible to prove Cartan's Formula by an ugly,
brute force calculation of both sides, but there is also an elegant,
no-sweat proof that I first learned from S.~S.~Chern (when I proudly showed
him my version of the ugly proof). There is an important involutory
automorphism $\omega \mapsto \bar \omega$ of the algebra $A$ of 
differential forms on a manifold. Namely, it is the identity on 
forms of even degree and is minus the identity on forms of odd degree.
A linear map $\partial: A \to A$ is called an {\it anti-derivation\/} if
$\partial(\lambda\omega) = 
    \partial\lambda \wedge\omega + \bar\lambda\wedge\partial\omega$.
It is of course well-known that the exterior derivative, $d$, is an 
anti-derivation (of degree $+1$) and an easy check shows that the interior
product $i_X$ is an anti derivation (of degree $-1$). Moreover, the 
anti-commutator of two anti-derivations is clearly a derivation, so that
$\LD(X?)$ and $d i_X + i_X d$ are both derivations of $A$, and hence to 
prove they are equal it suffices to check that they agree on a set of
generators of $A$. But $A$ is generated by forms of degree zero 
(i.e., functions) and the differentials of functions, and it is
obvious that $\LD(X?)$ and $d i_X + i_X d$ agree on these.
\smallskip

   We shall also have to deal with symplectic structures on infinite
dimensional manifolds. In this case we  still require that $\Omega$ is a
closed form and we also still require that $\Omega$ is {\it weakly\/}
non-degenerate,  meaning that for each point $p$ of $P$, the map
$v \mapsto i_v\Omega$ of $TP_p$ to $TP_p^*$ is injective. In finite
dimensions this of course implies that $\Omega$ is strongly
non-degenerate---meaning  that the latter map is in fact an
isomorphism---but that is rarely  the case in infinite dimensions,
so we will {\it not\/} assume it.  Thus, if $F$ is a smooth function
on $P$, it does not automatically  follow that there is a symplectic
gradient vector field $\sgrad F$ on $P$ satisfying
$\Omega((\sgrad F)_p, v) = dF_p(v)$ for all $v$ in
$TP_p$---this must be proved separately. However, if a symplectic
gradient does exist, then weak non-degeneracy shows that it is unique.
In the infinite dimensional setting
we call a function $F:P \to \reals$ a {\it Hamiltonian\/}
function if it has a symplectic gradient, and vector fields of the
form $\sgrad F$ will be called Hamiltonian vector fields. Obviously
the space of Hamiltonian functions is linear, and in fact the
formula $d (FG) = F dG + G dF$ shows that it is even an algebra, and that
$\sgrad(FG) = F\sgrad G + G \sgrad F$. We shall call a vector field $X$
on $P$ symplectic if the $1$-form $i_X \Omega$ is a closed but not
necessarily exact, for as we have seen, this is the condition
for the flow generated  by $X$ to preserve $\Omega$.

  Of course if $P$ is a vector space the distinction between
Hamiltonian and symplectic disappears: if $i_X\Omega$ is closed,
then $H(p) = \int_0^1 \Omega_{tp}(X_{tp}, p)\, dt$ defines a
Hamiltonian function with $\sgrad H = X$. Moreover, in this case it
is usually straightforward to check if $i_X\Omega$ {\it is\/} closed.
Given $u,v$ in $P$, consider them as constant vector fields on $P$,
so that $[u,v] = 0$. Then the formula
$d \theta (u,v) = u(\theta(v)) - v(\theta(u)) -\theta([u,v])$
for the exterior derivative of a \hbox{$1$-form} shows that symmetry of
${d\over dt}\big | _{t = 0} \Omega(X_{p+tu},v)$ in $u$ and $v$
is necessary and sufficient for $i_X\Omega$ to be
closed (and hence exact). In case $\Omega$ is a constant form
(i.e., $\Omega_p(u,v)$ is independent of $p$) then
${d\over dt}\big | _{t = 0} \Omega(X_{p+tu},v) = \Omega((DX_p)(u),v)$,
where $(DX)_p(u)={d\over dt}\big | _{t = 0} X_{p+tu}$
is the differential of $X$  at $p$. Since $\Omega$
is skew-symmetric in $u$ and $v$, this shows that if $\Omega$ is constant
then $X$ is Hamiltonian if and only if $(DX)_p$ is ``skew-adjoint'' with
respect to $\Omega$.

  If two smooth real-valued functions $F_1$ and $F_2$ on a symplectic
manifold $P$ are Hamiltonian, i.e., if they have symplectic gradients
$\sgrad F_1$ and $\sgrad F_2$, then they determine a third function 
on $P$, called  their {\it Poisson bracket\/}, defined by:
$$\PB (F_1,F_2) = \Omega(\sgrad F_2, \sgrad F_1).$$
The formula
$i_{{}_{[X,Y]}}\Omega = d(\Omega(Y,X))$
shows that the Poisson bracket is also a Hamiltonian
function, and in fact
$$\sgrad \PB (F_1,F_2) = [\sgrad F_1, \sgrad F_2].$$
What this formula says is that Hamiltonian functions
$F:P \to \reals$ are not only a commutative and associative algebra
under pointwise product, but also a Lie algebra under Poisson bracket,
and $F \mapsto \sgrad F$
is a Lie algebra homomorphism of this Lie algebra onto the
Lie algebra of Hamiltonian vector fields on $P$. In particular,
we see that the Poisson bracket satisfies the Jacobi identity,
$$\PB(\PB (F_1,F_2),F_3) + \PB(\PB (F_2,F_3),F_1) +
\PB(\PB (F_3,F_2),F_2) = 0,$$
and the Leibnitz Rule $\sgrad(FG) = F\sgrad G + G \sgrad F$ gives:
$$\PB(F_1, F_2 F_3) = \PB (F_1,F_2) F_3 + F_2 \PB(F_1,F_3),$$
which we will also call the Leibnitz Rule.

\smallskip
\noindent
{\bf Remark.} A {\it Poisson structure\/} for a smooth manifold
is defined to be a Lie algebra structure $\PB(F,G)$ on the algebra
of smooth functions that satisfies the Leibnitz Rule.
\smallskip

Since
$\PB (F_1,F_2) = \Omega(\sgrad F_2, \sgrad F_1) =
dF_2(\sgrad F_1) = \sgrad F_1 (F_2)$,
we can interpret the Poisson bracket of $F_1$ and $F_2$ as
the rate of change of $F_2$ along the solution curves of the
vector field $\sgrad F_1$. If we are considering some fixed Hamiltonian
system ${dx\over dt} = \sgrad H_x$ on $P$, then we can write this as
${dF\over dt} = \PB (H,F)$, and we see that {\it the vanishing of the
Poisson bracket $\PB (H,F)$ is the necessary and sufficient condition
for $F$ to be a constant of the motion\/}. By the Jacobi Identity, 
a corollary to this observation is that the Poisson Bracket of two
constants of the motion is also a constant of the motion. And since
$\PB (H,H) = 0$, $H$ itself is always a constant of the motion. 
(This is a proof of conservation of energy from the Hamiltonian 
point of view,  and below we  will also see how to prove Noether's 
Theorem in the Hamiltonian framework.)

Since the Poisson bracket is skew-symmetric, $\PB (F_1,F_2)$ is zero
if and only if $\PB (F_2,F_1)$ is zero, and in this case we say that
$F_1$ and $F_2$ are {\it in involution\/}. More generally $k$
Hamiltonian functions  $F_1,\ldots, F_k$ are said to be in involution
if all of the Poisson brackets $\PB(F_i,F_j)$ vanish. Note that since
$\sgrad \PB (F_i,F_j) = [\sgrad F_i, \sgrad F_j]$, if the $F_i$
are in involution then the vector fields $\sgrad F_i$
commute, i.e., $[\sgrad F_i, \sgrad F_j] = 0$, or equivalently the
flows they generate commute. In particular we see that
{\it if $F_1,\ldots, F_n$ are in involution and if each $\sgrad F_i$
generates a one parameter group of diffeomorphisms $\phi_t^i$ of $P$
then $(t_1,\ldots,t_n) \mapsto
\phi_{t_1}^1 \circ \phi_{t_2}^2\circ\ldots\circ \phi_{t_n}^n$
defines a symplectic action of the abelian group $\reals^n$ on $P$\/}.

  Suppose $P$ is a symplectic manifold of dimension $2n$ and that 
there exist $n$ functions $F_i$ such that the $dF_i$ are everywhere
linearly independent. If the functions $F_i$ are in involution with 
each other and with a function $H$, then the so-called Arnold-Liouville 
Theorem ([Ar], Chapter 10) states that the Hamiltonian system $\sgrad H$ 
is completely integrable in the sense mentioned earlier, i.e., there exist 
action-angle variables $q_1,\ldots,q_n,p_1,\ldots,p_n$ . In fact, complete
integrability of a $2n$ dimensional Hamiltonian system is often {\it defined\/} 
as the existence of $n$ functionally independent constants of the motion in 
involution. 

  This leads naturally to two interesting problems: finding ways 
to construct symplectic manifolds with lots functions in involution, and
determining whether a given Hamiltonian system is completely integrable.
In the late 1970's M.~Adler [Ad], B.~Kostant [Kos], and W.~Symes [Sy] 
independently and nearly simultaneously found a beautiful approach to the 
first question using certain special splittings of Lie algebras. For 
excellent surveys of finite dimensional completely integrable systems 
see [AdM] and [Pe]\null. The Adler-Kostant-Symes Theorem is explained in
detail in both of these references, and we shall not discuss it further 
here, except to note that it is closely related to an earlier method of 
Peter Lax [La1], that will be one of our main tools in later sections, 
and that, as Adler's paper showed, the Adler-Kostant-Symes Theorem also 
applies to infinite dimensional systems. In fact Adler's paper applied 
the method to the KdV equation, and later many other PDE were treated by 
the A-K-S approach in [Dr], [DS], [RS], [Se1], [Se2], and [Te2].

  As for the second problem,there is no magic test to check if a given 
system is completely integrable, and the  principal technique is to try 
to show that it can be manufactured using the  Adler-Kostant-Symes method. 
In fact, one often hears it said that ``all known  completely integrable 
systems arise in this way''.

 If a symplectic structure $\Omega$ is
``exact''---i.e., if $\Omega = d\omega$ for some $1$-form $\omega$
on $P$ (as we saw was the case for a cotangent bundle)
and if a vector field $X$ not only preserves $\Omega$, but
even preserves $\omega$, then Cartan's formula gives
$0 = \LD(X?\omega) = di_X\omega + i_X \Omega$, so if we define
$X^\omega = -i_X\omega = -\omega(X)$, then $\sgrad (X^\omega) =X$.
If $Y$ is a second such vector field on $P$, then
a computation completely analogous to that for $i_{{}_{[X,Y]}}\Omega$
above (replacing $\Omega$ by $\omega$) gives
$[X,Y]^\omega = \omega([Y,X]) = i_{[Y,X]}\omega = i_Y d(i_X \omega)
= -dX^\omega(Y) = -dX^\omega(\sgrad Y^\omega) = \PB(X^\omega,Y^\omega)$.
Thus $X \mapsto X^\omega$ is a Lie algebra homomorphism
inverse to $F \mapsto \sgrad F$ from the Lie algebra of vector
fields preserving $\omega$ to the Lie algebra of Hamiltonian
functions under Poisson bracket.

 In particular going back to Newton's Equations on our configuration
space $\C$, we see that if $X$ is a Killing vector field on $\C$ such
that $XU = 0$ then $\omega(X)$ is a constant of the motion for
Newton's Equations. It is easy to see that $\omega(X)$ is just
the conjugate momentum of $X$, so this gives a proof of Noether's
Principle in the Hamiltonian framework.
\goodbreak
\subsection {Examples of Classical Mechanical Systems}

\noindent
  While any choice of potential function $U$ on any Riemannian
manifold $\C$ defines a ``Classical Mechanical System'', in some
generalized sense, this name is often reserved for certain more
special cases that arise from physical considerations.

  One important and interesting class of examples describes the
motion of rigid bodies or ``tops'' with no external forces acting.
Here the configuration space $\C$ is the rotation group $\SOn(3)$,
while the metric tensor (also called the Inertia tensor in this case)
is any left-invariant metric on $\C$, and $U = 0$. We refer the
reader to any book on Classical Mechanics (e.g., [AbM], [Ar]) for a
discussion of these example, but be warned, that the full theory
is covered in a multi-volume treatise [KS]\null. An excellent recent book
is [Au]\null.

  A second important class of examples, usually referred to as
``particle mechanics'', describes the motion under mutual forces
of $N$ particles in the Euclidean space $\reals^k$ (where
usually $k = 1,2,$ or $3$). In this  case $\C = (\reals^k)^N$,
a point $x = (x_1,\ldots, x_N)$ of $\C$ representing the positions of
$N$ particles. For an important subclass, the force on each particle
is the sum of forces exerted on it by the remaining particles. In
this case the potential $U$ is a function of the distances
$r_{ij} = \norm|x_i - x_j|$ separating the particles.
It follows that the Lie group $G$ of Euclidean motions of $\reals^k$
is a group of symmetries, so the conjugate momenta of
the Lie algebra of $G$ give $k$ linear momenta (from the translations)
and $k(k-1)/2$ angular momentum (from the rotations)
that are conserved quantities.

   A simple but important example from particle mechanics is the
``harmonic oscillator''.  Here $k = N = 1$, so $\C = \reals$,
the metric on $T\C = \reals \times \reals$ is given by
$\norm|(x,v)|^2 = m v^2$ (where $m$ is the mass of the oscillator)
and $U(x) = \halb k x^2$, where $k>0$ is the so-called
spring constant of the oscillator. This models a particle that
is in equilibrium at the origin, but which experiences a Hooke's Law
linear ``restoring  force'' of magnitude $-kx$ directed towards the
origin when it is  at the point $x$ in $\C$. Newton's Equation of
motion is
$m\ddot x = -kx$, and the solutions are of the form
$x(t) = A\cos(\omega (t-t_0))$, where the angular
frequency $\omega$ is $\sqrt{k/m}$. The Hamiltonian formulation
of the harmonic oscillator is given in terms of canonical variables
$q = x$ and $p = m (dx/dt)$ by $H(q,p) = \halb(p^2/m + kq^2)$. Note
that $P = \halb(p^2 + mk q^2)$ and $Q = \arctan(p/q\sqrt{mk})$
define action-angle variables for the harmonic oscillator.

  Only notationally more complicated is the case of $N$ uncoupled
harmonic oscillators, with masses $m_1,\ldots,m_N$ and spring constant
$k_1,\ldots,k_N$. Now $\C = \reals^N$, the metric on
$T\C = \reals^N \times \reals^N$ is given by
$\norm|(x,v)|^2 = \sum_i m_i v_i^2$, and the potential function
is $U(x) = \halb \sum_i k_i x_i^2$. Newton's Equations are
$m_i\ddot x_i = -k_i x_i$ with the solutions
$x_i(t) = A_i\cos(\omega_i (t-t_0^i))$,
where $\omega_i = \sqrt{k_i/m_i}$.  The Hamiltonian for this example is
$H(q,p) = \sum_i\halb(p_i^2/m_i + kq_i^2)$. Note that not only
is the total Hamiltonian, $H$, a constant of the motion, but so
also are the $N$ partial Hamiltonians,
$H_i(q,p) = \halb(p_i^2/m_i + kq_i^2)$---i.e., the sum of the
kinetic plus potential energy of each individual oscillator
is preserved during the time evolution of any solution.
In this case we get one pair of action-angle variables from the
action-angle variables for each of the individual harmonic
oscillators, so it is again completely integrable.

   A seemingly more complicated example is the case of $N$
{\it coupled\/} harmonic oscillators. Starting from the
previous example, we imagine adding Hooke's Law springs
with spring constants $K_{ij}$ joining the $i$-th and $j$-th
particles. The force on the $i$-th particle is now
$F_i = - k_i x_i - K_{ij}(x_i - x_j)$, so we can take as
our potential function
$U(x) = \halb \sum k_i x_i^2 + \halb \sum_{ij} K_{ij} (x_i - x_j)^2$.
Notice that this is clearly a positive definite quadratic form,
so without loss of generality we can consider the somewhat more
general potential function $U(x) = \halb \sum_{ij} k_{ij} x_i x_j$,
where $k_{ij}$ is a positive definite symmetric matrix. Newton's Equations
are now $m_i\ddot x_i = -\sum_j k_{ij} x_j$.
Because of the off-diagonal elements of $ k_{ij}$ (the so-called
``coupling constants'') Newton's Equations no longer have separated
variables, and integrating  them appears much more difficult. This is
of course an illusion; all that is required to reduce this case to the
case of uncoupled harmonic oscillators is to diagonalize the quadratic
form that gives the potential energy, i.e., find an orthonormal basis
$e_i,\ldots,e_n$ such that if $y = y_1 e_i + \ldots + y_n e_n$ then
$U(y) = \halb \sum_i \lambda_i y_i^2$.
The solutions of Newton's Equations are now all of the form
$\sum_i A_i \cos(\sqrt{\lambda_i}(t_i - t_0^i))e_i$. Solutions for
which one $A_i$ is non-zero and all the others are zero are referred
as ``normal modes'' of the coupled harmonic oscillator system.
Since the coupled harmonic oscillator system is just the uncoupled
system in disguise, we see that it also is completely integrable.
Moreover, when we express a solution $x(t)$ of Newton's Equations
as a sum of normal modes, then not only is the kinetic energy plus
the potential energy of $x(t)$ a constant of the motion, but also
the kinetic plus the potential energy of each of these normal modes
is also a constant of the motion.

   There are two properties of the coupled harmonic oscillators
that make it an exceptionally important model system. First, it
is exactly and explicitly solvable, and secondly, as we shall see
in the next  section, it is an excellent first approximation to what
happens in an arbitrary system near a so-called ``vacuum solution''.
i.e., a stable equilibrium.

\subsection {Physics Near Equilibrium}

\noindent
  Physical systems are normally close to equilibrium, so it is
important to analyze well what happens in the phase space of a
physical system in the near neighborhood of an equilibrium point.

  We shall assume that our system is described as above by a
potential $U$ on a configuration space $\C$. By an ``equilibrium
point'' we mean a point $p$ of $\C$ that is not just a critical
point of $U$, but in fact a non-degenerate local minimum. Since
$U$ is only determined up to an additive constant, we can assume
that $U(p) = 0$. Since $\grad U$ vanishes at $p$, it is clear
that $\sigma(t) = p$ is a solution of Newton's Equations, and
physicists sometimes refer to such a solution as a ``vacuum
solution''.

  By a famous result of Marston Morse, we can find local
coordinates $y_1,\ldots, y_n$, in a neighborhood $O$ of $p$
and centered at $p$ such that
$U(q) = \sum_i y_i(q)^2$, so that
$N(\es) = \{q \in O \mid U(q) < \es\}$
is a neighborhood basis for $p$.  It follows that a vacuum
solution is stable, i.e., a solution of
Newton's Equations with initial conditions sufficiently close to
those of a vacuum solution will remain close to the vacuum
solution for  all time. To be precise, suppose
$\gamma(t)$ is a solution of Newton's Equations such that
$\gamma(0)$ is in $N(\halb\es)$ and $K(\gamma^\prime(0)) < \halb \es$.
Then $U(\gamma(0)) + K(\gamma^\prime(0)) < \es$,
so that, by conservation of total energy,
$U(\gamma(t)) + K(\gamma^\prime(t)) < \es$ for all $t$, and since
$K$ is non-negative, $U(\gamma(t)) < \es$ for all $t$,
i.e., the solution $\gamma(t)$ remains inside $N(\es)$.

  But we can be much more precise about the nature of these
solutions that are near the vacuum. To simplify the exposition
somewhat we will make the (inessential) assumption that the
metric on $\C$ is flat---as it usually is in particle mechanics.
Then we can choose orthogonal coordinates $x_1,\ldots,x_n$
centered at $p$ that simultaneously diagonalizes both the
kinetic energy and the Hessian matrix of $U$ at $p$, and the
assumption that $p$ is a non-degenerate local minimum just means
that the diagonal elements,  $k_i$, of the Hessian are positive.
(The diagonal elements, $m_i$, of the kinetic energy are of course
also positive, and have the interpretations of masses)
Thus, by Taylor's Theorem,
$U(x) = \halb \sum_j k_j x_j^2 +
{1\over 6}\sum_{jkl} a_{jkl}(x) x_j x_k x_l$,
where the functions $a_{jkl}(x)$ are smooth and symmetric
in their last two indices, and Newton's Equations take the form:
$$m_i{d^2 x_i(t)\over dt^2} =  -k_i x_i -
    \sum_{jk}a_{ijk}(x)x_k x_l + O(\norm|x|^3)  .$$
(For later reference, we note that if we adopt a Hamiltonian viewpoint,
and move to the cotangent bundle using the Legendre transform,
then in the canonical symplectic coordinates associated to $x_1,\ldots,x_n$,
the kinetic energy is $K$ is given by $K = {1\over 2}\sum_i{p_i^2\over m_i}$,
the potential energy is $U = \halb \sum_j k_j q_j^2 +
{1\over 6}\sum_{jkl} a_{jkl}(q) q_j q_k q_l$, and the Hamiltonian
is $H = K + U$.)

  The system of uncoupled harmonic oscillators obtained by
dropping the nonlinear terms is called the ``linearized system''
(at the given equilibrium $p$) and its normal modes are
referred to by physicists as the ``degrees of freedom''
of the system.

  An obvious question is, ``To what extent do solutions of
the linearized system approximate those of the full system?''.
One answer is easy, and no surprise---Gronwal's Inequality
implies that, as the initial position tends to $p$ and the initial
velocity tends to zero, a solution of the linearized equation
approximates that of the full equation better and better, and
for a longer period of time.

  A more subtle, but also more interesting question is, ``How will
the kinetic and potential energy of a solution become distributed,
on average, among the various degrees of freedom of the full system?''.
It is not difficult to give a precise formulation of this question.
The kinetic energy in the $i$-th mode is clearly
$K_i = \halb {p_i^2\over m_i}$,
and it is natural to assign to the $i$-th mode the potential energy
$U_i = \halb  k_i q_i^2 + {1\over 6}\sum_{kl} a_{ikl}(q) q_i q_k q_l$.
Then $H_i = K_i + U_i$ is that part of the total energy in the
$i$-th mode, and the total energy $H$ is just the sum of these $H_i$.
We know that for the linearized system each of the $H_i$
is a constant of the motion; that is, $H_i$ is constant
along any solution of Newton's Equations.
But it is easy to see that cannot be true for the full system,
and energy will in general flow between the normal modes because
of the nonlinear coupling between them. The question is, will the
``average behavior'' of the $H_i$ and $K_i$ have some predictable
relationship over large time intervals. 

   To make the concept of
``average'' precise, given any function $F:T\C \to \reals$,
define its ``time average'', $\bar F$ along a given solution $x(t)$ by:
$\bar F = \lim_{T \to \infty}{1\over T}\int_{-T}^T F({\dot x}(t)) \,dt$.
Then, what can we say about the time averages of the above partial
energy functions and their relations to each other. Of course a
first question is whether the limit defining the time average really
exists, and this is already a non-trivial point. Fortunately,
as we shall see in the next section, it is answered by the ``Individual
Ergodic Theorem'' of G.~D.~Birkhoff, [Bi], according to which the time
average will exist for ``almost all'' initial conditions.

   Starting in the late Nineteenth Century, physicists such as
Maxwell, Boltzmann, and Gibbs developed a very sophisticated
theory of statistical mechanics that gave convincing explanations
for (and good predictions of) the behavior of large assemblages
of molecules. The theoretical foundations for this theory were
based on just such time averages and their {\it hypothesized\/}
equality with another kind of average that is easier to investigate,
so-called ``space averages'', or ``microcanonical averages''.
As we will see, the space average of the kinetic energy in each normal 
mode is the same---a fact referred to as ``equipartition of energy''.
This important fact is the very basis for the definition of temperature
in statistical mechanics. Namely, for a system near equilibrium,
if the absolute temperature is $T$, then the average kinetic energy in
each degree of freedom is $kT\over 2$, where $k$ is the so-called
Boltzmann constant.

  But it is the {\it time\/} averages of the kinetic energy that should
really  determine the temperature, and if energy equipartition holds for
time averages, and if the system is experimentally started in one  of
its normal modes and is then followed in time, one should see an
equilibration take place, in which the kinetic energy should gradually
flow out of the original single mode in which it was concentrated and
become equally divided (on average) among all the various degrees of
freedom of the system. Because of the above relation between temperature
and equipartition of energy, this hypothesized equilibration process is
referred to as ``thermalization''. Intuitively speaking, this refers
to the transformation of the large scale motion of the system in a
single mode into ``heat'', i.e., lots of tiny fluctuating bits of energy
of amount  $kT\over 2$ in each of the many degrees of freedom.

  It should now be clear why physicists placed so much emphasis on
proving the supposed equality of the time average and the microcanonical
average, but mathematically this proved to be a highly intractible
problem. There were heuristic proofs, based on vague physical reasoning,
and also semi-rigorous arguments based on so-called ``ergodic hypotheses''.
The latter were assumptions to the effect that the solution curves would
wander on an energy surface in a sufficiently space filling way (ergodic
comes from the Greek word for energy). Unfortunately these ergodicity
assumptions were vague and in certain cases topologically impossible, and
it was only with the development of measure theory that von Neumann and
Birkhoff were able to state the precise condition (``metric transitivity'')
under which one could prove that time and space averages must necessarily
coincide.

   Nevertheless, physicists were morally convinced of the correctness
of the time-average based concept of thermalization; so much so that when
Fermi, Pasta, and Ulam undertook the numerical experiments that we will
consider later, they stated that their goal was not so much to discover
{\it if\/}  there would be be thermalization, but rather to discover
experimentally what the rate of approach to thermalization would be!

  For those readers who are interested, we will provide more of the
mathematical details concerning equipartition of energy in the next
section.

\subsection {Ergodicity and Thermalization}

   Let $P$ is a symplectic manifold (say of dimension $2n$)
with symplectic $2$-form $\Omega$, and let $H$ denote a
Hamiltonian function on $P$, generating a symplectic flow
$\phi_t$, that is the infinitesimal generator of $\phi_t$
is $\sgrad H$, the symplectic gradient of $H$.
As we have seen, this implies that the flow $\phi_t$
preserves the symplectic structure, and also that $H$
is a ``constant of the motion'', meaning  that it is constant
along every orbit, $\phi_t(p)$, or equivalently, that the constant
energy hypersurfaces $\Sg{c}$ (defined by $H = c$) are invariant under
the flow.  In classical examples, the Hamiltonian is usually bounded
below and  proper (so that all the $\Sg{c}$ are compact) and we shall
assume  this in what follows. Since $H$ is only defined up to an
additive constant we can assume the minimum value of $H$ is zero.

   The $2n$-form $\Omega^n$ defines a measure $d\mu$ on $P$ (the
 \hbox{Liouville} measure), and this is of course invariant under the
flow. We can factor $\Omega^n$ as $\Omega^n = \lambda \wedge dH$,
and the $2n-1$ form $\lambda$ is uniquely determined modulo the
ideal generated by $dH$, so it induces a unique measure on each
energy hypersurface $\Sg{c}$. We will denote these measures by
$d\sigma$, and they are of course likewise invariant under the flow.
Since $\Sg{c}$ is compact, its total measure, $\sigma(c)$, is finite,
and so, for any integrable function $f$ on $\Sg{c}$, we can define
its spatial average by
$\hat f = \sigma(c)^{-1} \int _{\Sg{c}} f(x)\, d\sigma(x)$.
(This is the quantity called the ``microcanonical average''
of $f$ in statistical mechanics.)
We note that these measure $d\sigma$ are canonically determined in
terms of the Liouville form, $\Omega^n$, and the Hamiltonian function
$H$, so {\it if $\psi$ is any diffeomorphism of $P$ that preserves
$\Omega^n$ and preserves $H$, then $\psi$ will also preserve the $d\sigma$
and hence all microcanonical averages, i.e., if $g = f\circ \psi$,
then $\hat g = \hat f$\/}.

  We return now to the question of ``equipartition of energy''.
We assume that we have canonical variables
$(p_1,\ldots,p_n,q_1,\ldots,q_n)$ in $P$ in which $H$ takes the
classical form $H = K + U$ where $K = {1\over 2}\sum_i{p_i^2\over m_i}$,
and $U$ is a function of the $q_i$ with a non-degenerate local minimum,
zero, at the origin. (It follows that for small $c$ the energy surfaces
$\Sg{c}$ are not only compact, but are in fact topologically
spheres.)  Since the $p$'s and $q$'s are canonical, $\Omega$ has the
standard Darboux form $\sum_i dp_i \wedge dq_i$, and so the
Liouville $2n$-form is just
$dp_1\wedge dq_1\wedge\ldots \wedge dp_n\wedge dq_n$, giving
Lebesgue measure as the Liouville measure in these coordinates.
Our goal is to prove that if $K_i ={p_i^2\over m_i}$,
then the microcanonical averages $\hat K_i$, $i = 1,\ldots,n$
(over any fixed energy surface $\Sg{c}$) are all the same.
Without loss of generality we can assume that $i = 1$ and $j = 2$,
and by the remark above it will suffice to find
a diffeomorphism $\psi$ that preserves $H = K + U$
and the Liouville form such that $K_2 = K_1\circ \psi$.
In fact, define
$$\psi(p_1,p_2,p_3,\ldots,p_n,q_1,\ldots,q_n) =
    (\alpha p_2,\alpha^{-1}p_1,p_3,\ldots,p_n,q_1,\ldots,q_n),$$
where $\alpha =\sqrt{{m_1/m_2}}$. Now, while $\psi$ is clearly
{\it not\/} symplectic, it just as clearly does preserve the
Liouville form. Moreover a trivial calculation shows that
$K_2 = K_1\circ \psi$ and $K_1 = K_2\circ \psi$, while
$K_i = K_i\circ \psi$ for $i > 2$. Since $K = \sum_i K_i$,
$K = K\circ \psi$. Since $U$ is a function of the $q$'s and not the
$p$'s, $U = U\circ \psi$, so $H = H\circ \psi$ also, and this completes
the proof that $\hat K_2 = \hat K_1$.

  There is an important corollary of the above proof. Suppose that
we can write the potential energy $U$ as the sum of $n$ functions
$U_i$, and let us define $H_i = K_i + U_i$. You should think of $U_i$
as representing the ``potential energy in the $i$-th normal mode'', 
and similarly $H_i$ represents the part of the total energy that is 
``in'' the $i$-th normal mode. In applications where the potential $U$ 
describes an interaction between identical particles, these partial
potentials will satisfy 
$U_1(q_1,q_2,\ldots,q_n) = U_2(q_2,q_1,\ldots,q_n)$, and similarly
for other pairs of indices. (For the example of the preceding section,
we note that these conditions will be satisfied if the ``spring constants''
$k_i$ are all equal and if the functions $a_{ijk}$ are symmetric in all
three indices.) We remark that, in particular, these conditions are 
satisfied for the Fermi-Pasta-Ulam Lattice that we will consider shortly.
If we now redefine $\psi$ above to simply 
interchange $q_i$ and $q_j$, then the same argument as
before shows that $\hat U_i = \hat U_j$, and so of course
we also have $\hat H_i = \hat H_j$. In words, for such systems 
not only kinetic energy per mode, but also potential and total 
energies per mode are ``equi-partitioned'', in the sense 
that their microcanonical averages are equal.

   Next recall that for $p$ in $\Sg{c}$
we define the time average of $f$ on the orbit of $p$ by:
$$\bar f(p) = \lim_{T \to \infty}{1\over T} \int^T_{-T} f(\phi_t(p))\, dt,$$
provided the limit exists.   G.~D.~Birkhoff's Individual Ergodic Theorem
([Bi]) states that $\bar f(p)$ is defined for
almost all $p$ in $\Sg{c}$, and then clearly $\bar f$
is invariant under the flow. It is  moreover again an integrable
function on $\Sg{c}$ with the same spatial average as $f$ itself.
It is then easily seen that the following four conditions are equivalent:

\roster
\item "1)" For every integrable function $f$ on $\Sg{c}$, its time
     average $\bar f$ is constant (and hence equal to its spatial average).

\item "2)" Every measurable subset of $\Sg{c}$ that is invariant
   under the flow either has measure zero or has measure $\sigma(c)$.

\item "3)" If an integrable functions on $\Sg{c}$ is constant on each
     orbit of the flow then it is constant (almost everywhere) on $\Sg{c}$.

\item "4)" Given two subsets $E_1$ and $E_2$ of $\Sg{c}$ having positive
         measure, some translate $\phi_t(E_1)$ of $E_1$ meets $E_2$ in a
         set of positive measure.
\endroster

\noindent
and if these equivalent conditions are satisfied, then the flow is said
to be {\it ergodic\/} or {\it metrically transitive\/} on $\Sg{c}$.

By choosing $f$ to be the characteristic function of an open set $O$, we
see from 1) that ergodicity implies that the motion has a ``stochastic''
nature---that is, the fraction of time that an orbit spends in $O$ is equal
to the measure of $O$ (so in particular almost all orbits are dense in
$\Sg{c}$). This implies that (apart from $\Sg{c}$ itself) there cannot
exist any stable fixed point, periodic orbit, or more general stable
invariant set. To put it somewhat more informally, orbits on an
ergodic $\Sg{c}$ cannot exhibit any simple asymptotic behavior.

 Note that any function of a constant of the motion will again
be a constant of the motion---and in particular any function of $H$
is a constant of the motion. There may of course be constants of
the motion that are functionally independent of $H$. But if the
flow is ergodic on every energy surface, then it follows from 3) that
any constant of the motion, will be constant on each level set
of $H$---which is just to say that it is a function of $H$.
This shows that Hamiltonian systems with many independent constants
of the motion (and in particular completely integrable systems) are
in some sense at the opposite extreme from ergodic systems.

  So what is the status of the old belief that a ``generic'' (in some 
suitable sense) Hamiltonian system should be ergodic on each energy surface?  
On the one hand, Fermi [Fe] proved a result that points in this direction.
And there is a famous result of Oxtoby and Ulam ([OU]) to the 
effect that in the set of all measure preserving {\it homeomorphisms\/} 
of an energy surface, those that are metrically transitive are 
generic in the sense of category. But the measure 
preserving {\it diffeomorphisms\/} of an energy surface are
themselves only a set of first category in the measure preserving
homeomorphisms, so the Oxtoby-Ulam theorem is not particularly relevant
to this question. In fact, the KAM (Kolmagorov-Arnold-Moser) Theorem ([Ar],
Appendix 8) shows that any Hamiltonian flow that is sufficiently close to a
completely  integrable system in a suitable $C^k$ topology will have a set of
invariant  tori of positive Liouville measure, and so cannot be ergodic.
Indeed, proving  rigorously that any particular Hamiltonian system is 
ergodic is  quite difficult. For some examples of such theorems 
see [AA]\null.

\section {Origins of Soliton Theory}

  Perhaps the single most important event leading up to the explosive
growth of soliton mathematics in the last decades were some seemingly
innocuous numerical experiments, carried out by Enrico Fermi, John Pasta,
and Stanislaw Ulam in 1954--55, on the Los Alamos MANIAC computer.
(Originally published as Los Alamos Report LA1940 (1955) and reprinted
in [FPU]).

\subsection {The Fermi-Pasta-Ulam Experiments}

\noindent
The following quotation is taken from Stanislaw Ulam's
autobiography, ``Adventures of a Mathematician''.
\smallskip

{\narrower
  Computers were brand new; in fact the Los Alamos
Maniac was barely finished\dots.As soon as the machines were
finished, Fermi, with his great common sense and intuition,
recognized immediately their importance for the study
of problems in theoretical physics, astrophysics, and
classical physics. We discussed this at length and decided to
formulate a problem simple to state, but such that a solution
would require a lengthy computation which could not be done with
pencil and paper or with existing mechanical computers\dots$\!\!.$[W]e
found a typical one\dots the consideration of an elastic string
with two fixed ends, subject not only to the the usual elastic
force of stress proportional to strain, but having, in addition,
a physically correct nonlinear term\dots. The question was to
find out how\dots the entire motion would eventually thermalize\dots$\!\!.$

  John Pasta, a recently arrived physicist, assisted us in
the task of flow diagramming, programming, and running the
problem on the Maniac\dots$\!\!.$

   The problem turned out to be felicitously chosen. The results
were entirely different qualitatively from what even Fermi, with
his great knowledge of wave motion had expected.
\smallskip
\par}

\noindent
What Fermi, Pasta, and Ulam (FPU) were trying to do was to verify
numerically a basic article of faith of statistical mechanics; namely
the belief that if a mechanical system has many degrees of freedom and
is close to a stable equilibrium, then a generic nonlinear interaction
will ``thermalize'' the energy of the system, i.e., cause the
energy to become equidistributed among the normal modes of the
corresponding linearized system. In fact, Fermi believed he had
demonstrated this fact in [Fe]\null. Equipartition of energy among
the normal modes is known to be closely related to the ergodic
properties of such a system, and in fact FPU state their goal
as follows: ``The ergodic behavior of such systems was studied
with the primary aim of establishing, experimentally, the rate of
approach to the equipartition of energy among the various degrees
of freedom of the system.''

  FPU make it clear that the problem that they want to simulate is
the vibrations of a ``one-dimensional continuum'' or ``string''
with fixed end-points and nonlinear elastic restoring forces,
but that ``for the purposes of numerical work this continuum
is replaced by a finite number of points $\ldots$ so that the
PDE describing the motion of the string is replaced by a finite
number of ODE''\null. To rephrase this in the current jargon, FPU
study a one-dimensional lattice of $N$ oscillators with
nearest neighbor interactions and zero boundary conditions.
(For their computations, FPU take $N=64$.)

  We imagine the original string to be stretched along the
$x$-axis from $0$ to its length $\ell$. The $N$ oscillators
have equilibrium positions $p_i = ih$, $i = 0, \ldots, N-1$,
where $h = \ell/(N-1)$ is the lattice spacing, so their
positions at time $t$ are $X_i(t) = p_i + x_i(t)$,
(where the $x_i$ represent the displacements of the oscillators
from equilibrium). The force attracting any oscillator to one of
its neighbors is taken as $k(\delta + \alpha \delta^2)$,
$\delta$ denoting the ``strain'', i.e., the
deviation of the distance separating these
two oscillators from their equilibrium separation $h$.
(Note that when $\alpha = 0$ this is just a linear Hooke's
law force with spring constant $k$.)
The force acting on the $i$-th oscillator due to its right neighbor is
$F(x)_i^+ = k[(x_{i+1} - x_i) + \alpha((x_{i+1} - x_i)^2]$
while the force acting on the $i$-th oscillator due to its left neighbor is
$F(x)_i^- = k[(x_{i-1} - x_i) - \alpha((x_{i-1} - x_i)^2]$.
Thus the total force acting on the $i$-th oscillator will be the
sum of these two forces, namely:
$F(x)_i = k(x_{i+1} + x_{i-1} -2x_i)[1 + \alpha(x_{i+1} - x_{i-1})],$
and assuming that all of the oscillators have the same mass, $m$,
Newton's equations of motion read:
$$m\ddot x_i = k(x_{i+1} + x_{i-1} -2x_i)[1 + \alpha(x_{i+1} - x_{i-1})],$$
with the boundary conditions $x_0(t) = x_{N-1}(t) = 0$. In addition,
FPU looked at motions of the lattice that start from rest, i.e.,
they assumed that $\dot x_i(0) = 0$, so the motion of the lattice
is completely specified by giving the $N - 2$ initial displacements
$x_i(0)$, $i = 1, \ldots, N-2$. We shall call this the FPU
initial value problem (with initial condition $x_i(0)$).

  It will be convenient to rewrite Newton's equations in
terms of parameters that refers more directly to
the original string that we are trying to model. Namely,
if $\rho$ denotes the density of the string then $m = \rho h$, while
if $\kappa$ denotes the Young's modulus for the string, (i.e., the
spring constant for a piece of unit length) then $k=\kappa/h$
will be the spring constant for
a piece of length $h$. Defining $c = \sqrt{\kappa/\rho}$ we can now
rewrite Newton's equations as:
$$\ddot x_i = c^2 \left({x_{i+1} + x_{i-1} -2x_i \over h^2}\right)
      [1 + \alpha(x_{i+1} - x_{i-1})]. \leqno (FPU)$$
and in this form we shall refer to them as the FPU Lattice Equations.
We can now ``pass to the continuum limit'', i.e., by letting $N$
tend to infinity (so $h$ tends to zero) we can attempt to derive
a PDE for the function $u(x,t)$ that measures the displacement at
time $t$ of the particle of string with equilibrium position $x$.
We shall leave the nonlinear case for later, and here restrict
our attention to the linear case, $\alpha = 0$. If we take
$x = p_i$,  then by definition $u(x,t) = x_i(t)$ and since
$p_i + h = p_{i+1}$ while $p_i - h = p_{i-1}$, with $\alpha = 0$
the latter form of Newton's equations gives:
$$ u_{tt}(x,t) = c^2{u(x + h,t) + u(x - h,t) - 2u(x,t)\over h^2}.$$
By Taylor's formula:
$$f(x\pm h)= f(x) \pm h f^{'}(x) + {h^2\over 2!} f^{''}(x)
\pm  {h^3\over 3!} f^{'''}(x) + {h^4\over 4!} f^{''''}(x) + O(h^5),$$
and taking $f(x) = u(x,t)$ this gives:
$${u(x + h,t) + u(x - h,t) - 2u(x,t)\over h^2}
       = u_{xx}(x,t)  +
\Bigr ({{\scriptstyle h^2}\over {\scriptstyle 12}}\Bigr)  u_{xxxx}(x,t)
+ O(h^4),$$ so letting $h \to 0$, we find $u_{tt} = c^2 u_{xx}$, i.e.,
$u$ satisfies the linear wave equation, with propagation speed $c$,
(and of course the boundary conditions $u(0,t) = u(\ell,t) = 0$,
and initial conditions $u_t(x,0) =0$, $u(x,0) = u_0(x)$).

  This is surely one of the most famous initial value problems of
mathematical physics, and nearly every mathematician sees a
derivation of both the d'Alembert and Fourier version of its
solution early  in their careers.  For each positive integer $k$
there is a normal  mode or ``standing wave'' solution:
$$u_k(x,t) =
\cos\left({k \pi c t\over \ell}\right)\sin \left({k \pi x\over\ell}\right)$$
and the solution to the initial value problem is
$u(x,t) = \sum_{k=1}^\infty a_k u_k(x,t)$
where the $a_k$ are the Fourier coefficients of $u_0$:
$$a_k =
{2\over l}\int_0^\ell u_0(x) \sin \left({k \pi x\over\ell}\right)\,dx.$$

  Replacing $x$ by $p_j = j h$ in $u_k(x,t)$  (and using $\ell = (N-1)h$)
we get functions
$$ \xi_j^{(k)}(t) =
\cos\left({k \pi c t\over (N-1)h}\right)
\sin \left({k j \pi \over N-1}\right),$$
and it is natural to conjecture that these will be the normal modes for
the FPU initial value problem (with $\alpha = 0$ of course). This is
easily checked using the addition formula for the sine function.
It follows that, in the linearized case, the solution to the
FPU initial value problem with initial conditions $x_i(0)$ is given
explicitly by $x_j(t) = \sum_{k=1}^{N-2} a_k \xi_j^{(k)}(t)$, where the
Fourier coefficients $a_k$ are determined from the formula:
$$ a_k = \sum_{j=1}^{N-2} x_j(0) \sin \left({k j \pi \over N-1}\right).$$

  Of course, when $\alpha$ is zero and the interactions are linear
we are in effect dealing with $N-2$ uncoupled harmonic oscillators
(the above normal modes) and there is no thermalization. On the
contrary, the sum of the kinetic and potential energy of each
of the normal modes is a constant of the motion!

  But if $\alpha$ is small but non-zero, FPU expected (on the basis of
then generally accepted statistical mechanics arguments) that the energy
would  gradually shift between modes so as to eventually roughly equalize
the total of potential and kinetic energy in each of the $N-2$ normal
modes $\xi^{(k)}$. To test this they started the lattice in the
fundamental mode $\xi^{(1)}$, with various values of $\alpha$, and
integrated Newton's equations numerically for a long time interval,
interrupting the evolution from time to time to compute the total
of kinetic plus potential energy in each mode. What did they find?
Here is a quotation from their report:\par
\smallskip

{\narrower
Let us say here that the results of our computations show features
which were, from the beginning, surprising to us. Instead of a gradual,
continuous flow of energy from the first mode to the higher modes,
all of the problems showed an entirely different behavior. Starting
in one problem with a quadratic force and a pure sine wave as the initial
position of the string, we did indeed observe initially a gradual increase
of energy in the higher modes as predicted (e.g., by Rayleigh in an
infinitesimal analysis). Mode $2$ starts increasing first, followed by mode
$3$, and so on. Later on, however this gradual sharing of energy among the
successive modes ceases. Instead, it is one or the other mode that
predominates. For example, mode $2$ decides, as it were, to increase rather
rapidly at the cost of the others. At one time it has more energy than all
the others put together! Then mode $3$ undertakes this r\^ole. It is
only the first few modes which exchange energy among themselves, and
they do this in a rather regular fashion. Finally, at a later time,
mode $1$ comes back to within one percent of its initial value, so that the
system seems to be almost periodic.
\par}
\smallskip

  There is no question that Fermi, Pasta, and Ulam realized they had
stumbled onto something big. In his autobiography [Ul], Ulam devotes
several pages to a discussion of this collaboration. Here is a
little of what he says:\par
\smallskip

{\narrower
I know that Fermi considered this to be, as he said,
``a minor discovery.'' And when he was invited a year later to give
the Gibbs Lecture (a great honorary event at the annual American
Mathematical Society meeting), he intended to talk about it. He became
ill before the meeting, and his lecture never took place\dots.

  The results were truly amazing. There were many attempts to find the
reasons for this periodic and regular behavior, which was to be the
starting point of what is now a large literature on nonlinear vibrations.
Martin Kruskal, a physicist in Princeton, and Norman Zabusky, a
mathematician at Bell Labs wrote papers about it. Later, Peter Lax
contributed signally to the theory. \par}
\smallskip

\noindent
  Unfortunately, Fermi died in 1955, even before the paper cited
above was published. It was to have been the first in a series
of papers, but with Fermi's passing it fell to others to follow
up on the striking results of the Fermi-Pasta-Ulam experiments.

  The MANIAC computer, on which FPU carried out their remarkable
research, was designed to carry out some computations needed for
the design of the first hydrogen bombs, and of course it was a
marvel for its day. But it is worth noting that it was very weak by
today's standards---not just when compared with current supercomputers,
but even when compared with modest desktop machines. At a conference
held in 1977 Pasta recalled, ``The program was of course punched
on cards. A DO loop was executed by the operator feeding in the deck of cards
over and over again until the loop was completed!''

\subsection {The Kruskal-Zabusky Experiments}

\noindent
  Following the FPU experiments, there were many attempts
to explain the surprising quasi-periodicity of solutions
of the FPU Lattice Equations. However it was not until ten
years later that Martin Kruskal and Norman Zabusky took the
crucial steps that led to an eventual understanding of this
behavior [ZK]\null.											

  In fact, they made two significant advances. First they demonstrated
that, in a continuum limit, certain solutions of the FPU Lattice Equations
could be described in terms of solutions of the so-called Korteweg-de Vries
(or KdV) equation. And secondly, by investigating the initial value
problem for the KdV equation numerically on a computer, they discovered
that its solutions had remarkable behavior that was related to, but
if anything even more surprising and unexpected than the anomalous
behavior of the FPU lattice that they had set out to understand.

  Finding a good continuum limit for the nonlinear FPU lattice is a
lot more sophisticated than one might at first expect after the easy
time we had with the linear case. In fact the approach to the limit
has to be handled with considerable skill to avoid inconsistent results,
and it involves several non-obvious steps.

  Let us return to the FPU Lattice Equations
$$\ddot x_i = c^2 \left({x_{i+1} + x_{i-1} -2x_i \over h^2}\right)
      [1 + \alpha(x_{i+1} - x_{i-1})], \leqno (FPU)$$
and as before we let $u(x,t)$ denote the function measuring the displacement
at time $t$ of the particle of string with equilibrium position $x$, so
if $x = p_i$ then, by definition, $x_i(t) = u(x,t)$,
$x_{i+1}(t) = u(x+h,t)$, and $x_{i-1}(t) = u(x-h,t)$.
Of course $\ddot x_i =u_{tt}(x,t)$ and, as noted earlier, Taylor's Theorem
with remainder gives
$$\eqalign{{x_{i+1} + x_{i-1} -2x_i\over h^2} &=
  {u(x + h,t) + u(x - h,t) - 2u(x,t)\over h^2}\cr
       &= u_{xx}(x,t)  +
\Bigr ({{\scriptstyle h^2}\over {\scriptstyle 12}}\Bigr) u_{xxxx}(x,t) +
O(h^4).}$$
By a similar computation
$$\alpha(x_{i+1} - x_{i-1}) =
   (2\alpha h) u_x(x,t) +
\Bigr ({{\scriptstyle \alpha h^3}\over {\scriptstyle 3}}\Bigr) u_{xxx}(x,t) +
O(h^5),$$
so substitution in (FPU) gives
$$\Bigr ({{\scriptstyle 1}\over {\scriptstyle c^2}}\Bigr) u_{tt} - u_{xx} =
    (2\alpha h) u_x u_{xx} +
\Bigr ({{\scriptstyle h^2}\over {\scriptstyle 12}}\Bigr) u_{xxxx} + O(h^4) .$$
As a first attempt to derive a continuum description for
the FPU lattice in the nonlinear case, it is tempting to just
let $h$ approach zero and assume that $2\alpha h$
converges to a limit $\eps$. This would give the PDE
$$ u_{tt} = c^2 (1 + \eps u_x)u_{xx}$$
as our continuum limit for the FPU Lattice equations and the
nonlinear generalization of the wave equation. But this leads to
a serious problem. This equation is familiar in applied
mathematics---it was studied by Rayleigh in the last century---and
it is easy to see from examples that its solutions develop discontinuities
(shocks) after a time on the order of $(\eps c)^{-1}$,
which is considerably shorter than the time scale of the
almost periods observed in the Fermi-Pasta-Ulam experiments. It was Zabusky
who realized that the correct approach was to retain the term of order $h^2$
and study the equation
$${\Bigr ({{\scriptstyle 1}\over {\scriptstyle c^2}}\Bigr) u_{tt} - u_{xx} =
    (2\alpha h) u_x u_{xx} +
\Bigr ({{\scriptstyle h^2}\over {\scriptstyle 12}}\Bigr) u_{xxxx}.}\leqno
(ZK)$$
If we differentiate this equation with respect to $x$ and  make
the substitution $v = u_x$, we see that it reduces to the more
familiar Boussinesq equation
$$\Bigr ({{\scriptstyle 1}\over {\scriptstyle c^2}}\Bigr)v_{tt} =
v_{xx} + \alpha h {\partial (v^2) \over \partial x^2} +
\Bigr ({{\scriptstyle h^2}\over {\scriptstyle 12}}\Bigr) v_{xxxx},$$
(The effect of the fourth order term is to add dispersion to the
equation, and this smoothes out incipient shocks before they can
develop.)

  It is important to realize that, since $h \ne 0$, (ZK) cannot logically
be considered a true continuum limit of the FPU lattice. It should rather
be regarded as an asymptotic approximation to the lattice model that
works for small lattice spacing $h$ (and hence large $N$). Nevertheless,
we shall now see how to pass from (ZK) to a true continuum description
of the FPU lattice.

  The next step is to notice that, with $\alpha$ and $h$ small,
solutions of (ZK) should behave qualitatively like
solutions of the linear wave equation $u_{tt} = c^2 u_{xx}$, and
increasingly so as $\alpha$ and $h$ tend to zero.
Now the general solution of the linear wave equation is of course
$u(x,t) = f(x+ct) + g(x-ct)$, i.e., the sum of an arbitrary left
moving traveling wave and an arbitrary right moving traveling wave,
both moving with speed $c$. Recall that it is customary to simplify
the analysis in the linear case by treating each kind of wave
separately, and we would like to do the same here. That is, we would
like to look for solutions $u(x,t)$ that behave more and more like
(say) right moving traveling waves of velocity $c$---and for longer and
longer periods of time---as $\alpha$ and $h$ tend to zero.

   It is not difficult to make precise sense out of this requirement.
Suppose that $y(\xi,\tau)$  is a smooth function of two real variables
such that the map $\tau \mapsto y(\cdot,\tau)$ is uniformly continuous
from $\reals$ into the bounded functions on $\reals$ with the
sup norm---i.e., given $\eps > 0$ there is a positive $\delta$
such that $|\tau - \tau_0| < \delta$ implies
$|y(\xi,\tau) - y(\xi, \tau_0)| < \eps$.
Then for $|t - t_0| < T = \delta/(\alpha h c)$ we have
$|\alpha h c t - \alpha h c t_0| < \delta$, so
$|y(x - ct,\alpha h c t) - y(x - ct, \alpha h c t_0)| < \eps$.
In other words, the function $u(x,t) = y(x - ct, \alpha h c t)$
is uniformly approximated by the traveling wave
$u^0(x,t) = y(x - ct, \alpha h c t_0)$
on the interval $|t - t_0| < T $ (and of course $T \to \infty$
as $\alpha$ and $h$ tend to zero). To restate this a little
more picturesquely, $u(x,t) = y(x - ct, \alpha h c t)$ is
approximately a traveling wave whose shape gradually changes
in time. Notice that if $y(\xi,\tau)$ is periodic or almost
periodic in $\tau$, the gradually changing shape of the
approximate traveling wave will also be periodic or almost periodic.

  To apply this observation, we define new variables $\xi = x - ct$
and $\tau = (\alpha h) ct$. Then by the chain rule,
$\partial^k/\partial x^k = \partial^k/\partial \xi^k$,
$\partial/\partial t =
 -c(\partial/\partial \xi - (\alpha h)  \partial/\partial \tau)$,
and $\partial^2/\partial t^2 =
 c^2(\partial^2/\partial \xi^2
 - (2\alpha h)  \partial^2/\partial \xi\partial\tau)
 +(\alpha h)^2 \partial^2/\partial \tau^2)$.

\noindent
Thus in these new coordinates the wave operator transforms to:
$${1\over c^2}{\partial^2\over \partial t^2} -
{\partial^2\over\partial x^2} =
  -2\alpha h {\partial^2 \over \partial \xi \partial \tau} +
(\alpha h)^2 {\partial^2 \over \partial \tau^2},$$
so substituting $u(x,t) = y(\xi,\tau)$ in (ZK) (and dividing by
$-2\alpha h$) gives:
$$y_{\xi\tau} -
\Bigl({\scriptstyle{\alpha h}\over \scriptstyle 2 }\Bigr) y_{\tau\tau} =
- y_\xi y_{\xi \xi} -
\Bigr ({{\scriptstyle h}\over {\scriptstyle 24\alpha}}\Bigr)
y_{\xi\xi\xi\xi},$$
and, at last, we are prepared to pass to the continuum limit.
We assume that $\alpha$ and $h$ tend to zero at the same rate, i.e.,
that as $h$ tends to zero, the quotient $h/\alpha$ tends
to a positive limit, and we define
$\delta = \lim_{h\to 0}\sqrt{h/(24\alpha)}$.
Then $\alpha h = O(h^2)$, so letting $h$ approach zero gives
$y_{\xi \tau} + y_\xi y_{\xi \xi} + \delta^2 y_{\xi \xi \xi \xi} = 0$.
Finally, making the substitution $v = y_\xi$ we arrive at the KdV
equation:
$${v_\tau + v v_\xi + \delta^2 v_{\xi \xi \xi} = 0.}
\leqno (KdV)$$
{\bf Remark.}
Note that if we re-scale the
independent variables by $\tau \to \beta \tau$ and $\xi \to \gamma \xi$,
then the KdV equation becomes:

$$v_\tau +
\Bigr ({{\scriptstyle \beta}\over {\scriptstyle \gamma}}\Bigr) v v_\xi +
\Bigr ({{\scriptstyle \beta}\over {\scriptstyle \gamma^3}}\Bigr)
\delta^2 v_{\xi \xi \xi} = 0,$$
so by appropriate choice of $\beta$ and $\gamma$ we can obtain
any equation of the form
$v_\tau + \lambda v v_\xi+ \mu v_{\xi\xi\xi} = 0$, and any such
equation is referred to as ``the KdV equation''. A commonly used
choice that is convenient for many purposes is
$v_\tau +6 v v_\xi +  v_{\xi\xi\xi} = 0$, although the form
$v_\tau -6 v v_\xi +  v_{\xi\xi\xi} = 0$ (obtained by replacing
$v$ by $-v$) is equally common. We will use both these forms.
\smallskip

  Let us recapitulate the relationship between the FPU Lattice
and the KdV equation. Given a solution $x_i(t)$ of the FPU Lattice
we get a function $u(x,t)$ by interpolation---i.e., $u(ih,t) = x_i(t)$,
$i = 0,\ldots,N$. For small lattice spacing $h$ and nonlinearity
parameter $\alpha$ there will be solutions $x_i(t)$ so that the
corresponding $u(x,t)$ will be an approximate right moving
traveling wave with slowly varying shape, i.e., it will be of the
form $u(x,t) = y(x-ct, \alpha h c t)$ for some smooth function
$y(\xi,\tau)$, and the function $v(\xi,\tau) = y_\xi(\xi,\tau)$
will satisfy the KdV equation
$v_\tau + v v_\xi +  \delta^2 v_{\xi \xi \xi}= 0$, where
$\delta^2 = h/(24\alpha)$.

 Having found this relationship between the FPU Lattice and the
KdV equation, Kruskal and Zabusky made some numerical experiments,
solving the KdV initial value problem for various initial data.
Before discussing the remarkable results that came out of these
experiments, it will be helpful to recall some of the early
history of this equation.

\subsection {A First Look at KdV}

\noindent
  Korteweg and de Vries derived their equation in 1895
to settle a debate that had been going on since 1844, when the
naturalist and naval architect John Scott Russell, in an oft-quoted
paper [Ru], reported  an experience a decade earlier in which he
followed the bow wave of a  barge that had suddenly stopped in a
canal. This ``solitary wave'', some thirty feet long and a foot high,
moved along the channel at about eight miles per hour, maintaining
its shape and speed for over a mile as Russell raced after it on
horseback. Russell became fascinated with this phenomenon, and
made extensive further experiments with such waves in a wave tank
of his own devising, eventually deriving a (correct) formula for
their speed as a function of height. The mathematicians Airy and
Stokes made calculations which appeared to show that any such wave
would be unstable and not persist for as long as Russell claimed.
However, later work by Boussinesq (1872), Rayleigh (1876) and finally
the Korteweg-de Vries paper in 1895 [KdV] pointed out errors in the
analysis of Airy and Stokes and vindicated Russell's conclusions.

   The KdV equation is now accepted as controlling the dynamics of
waves moving to the right in a shallow channel. Of course, Korteweg
and de Vries did the obvious and looked for traveling-wave solutions
for their equation by making the Ansatz $v(x,t) = f(x-ct)$. When
this is substituted  in the standard form of the KdV equation it
gives $-cf^\prime + 6f f^\prime + f^{\prime\prime\prime} = 0$. If we
add the boundary conditions that $f$ should vanish at infinity, then
a fairly routine analysis leads to the one parameter family of
traveling wave solutions $v(x,t) = 2 a^2 \sech^2(a(x - 4a^2 t))$,
now referred to as the one-soliton solutions of KdV. (These are
of course the solitary waves of Russell.)  Note that the amplitude
$2 a^2$ is exactly half the speed $4 a^2$, so that taller waves
move faster than their shorter brethren.

  Now, back to Zabusky and Kruskal. For numerical reasons, they
chose to deal with the case of periodic boundary conditions---in
effect studying the KdV equation
$u_t + u u_x + \delta^2 u_{xxx} = 0$ (which they label (1) )
on the circle instead of on the line. For their published report,
they chose $\delta = 0.022$ and used the initial condition
$u(x,0) = \cos(\pi x)$. Here is an extract from their report
(containing the first use of the term ``soliton'') in which
they describe their observations:
\smallskip

{\narrower\noindent
(I) Initially the first two terms of Eq. (1) dominate and
the classical overtaking phenomenon occurs; that is $u$ steepens
in regions where it has negative slope. (II) Second, after $u$
has steepened sufficiently, the third term becomes important and
serves to prevent the formation of a discontinuity. Instead,
oscillations of small wavelength (of order $\delta$) develop on
the left of the front. The amplitudes of the oscillations grow,
and finally {\it each\/} oscillation achieves an almost steady
amplitude (that increases linearly from left to right) and has the
shape of an individual solitary-wave of (1). (III) Finally, each
``solitary wave pulse'' or {\it soliton\/} begins to move uniformly
at a rate (relative to the background value of $u$ from which the
pulse rises) which is linearly proportional to its amplitude.
Thus, the solitons spread apart. Because of the periodicity, two or
more solitons eventually overlap spatially and interact nonlinearly.
Shortly after the interaction they reappear virtually unaffected
in size or shape. In other words, solitons ``pass through'' one
another without losing their identity.$\,${\it Here we have a nonlinear
physical process in which interacting localized pulses do not scatter
irreversibly\/}.  \par}
\smallskip

\noindent
(If you are not sure what Zabusky and Kruskal mean here by ``the classical 
overtaking phenomenon'', it will be explained in the next section.)

  Zabusky and Kruskal go on to describe a second interesting
observation, a recurrence property of the solitons that goes a
long way towards accounting for the surprising recurrence observed
in the FPU Lattice. Let us explain again, but in somewhat different
terms, the reason why the recurrence in the FPU Lattice is so
surprising. The lattice is made up of a great many identical
oscillators. Initially the relative phases of these oscillators
are highly correlated by the imposed cosine initial condition.
If the interactions are linear ($\alpha = 0$), then the oscillators
are harmonic and their relative phases remain constant. But, when $\alpha$
is positive, the anharmonic forces between the oscillators cause their
phases to start drifting relative to each other in an apparently
uncorrelated manner. The expected time before the phases of all of
the oscillators  will be simultaneously close to their initial phases
is enormous, and increases rapidly with the total number $N$.
But, from the point of view of the KdV solitons,
an entirely different picture appears. As mentioned in the above
paragraph, if $\delta$ is put equal to zero in the KdV equation, it
reduces to the so-called inviscid Burgers' Equation, which exhibits
steepening and breaking of a negatively sloped wave front in a
finite time $T_B$. (For the above initial conditions, the
breaking time, $T_B$, can be computed theoretically to be $1/\pi$.)
However, when $\delta >0$, just before breaking would occur,
a small number of solitons emerge (eight in the case of the above
initial wave shape, $\cos(\pi x)$) {\it and this number depends only
on the initial wave shape, not on the number of oscillators\/}.
The expected time for their respective centers of gravity to all
eventually  ``focus'' at approximately the same point of the circle
is of course much smaller than the expected time for the much larger
number of oscillators to all return approximately to their original phases.
In fact, the recurrence time $T_R$ for the solitons turns out to be
approximately equal to $30.4 T_B$, and at this time the wave shape
$u(x,T_R)$ is uniformly very close  to the initial wave form $u(x,0) =
\cos(\pi x)$. There is a second  (somewhat weaker) focusing at  time $t = 2
T_R$, etc. (Note that these times are measured in units of the ``slow time'',
$\tau$, at which the shape of the FPU traveling wave evolves, not in
the ``fast time'', $t$, at which the traveling wave moves.) In effect,
the KdV solitons are providing a hidden correlation between the relative
phases of the FPU oscillators!

   Notice that, as Zabusky and Kruskal emphasize, it is the persistence
or shape conservation of the solitons that provides the explanation of
recurrence. If the shapes of the solitons were not preserved when they
interacted, there would be no way for them to all get back together
and approximately reconstitute the initial condition at some later
time. Here in their own words is how they bring in solitons to
account for the fact that thermalization was not observed in the FPU
experiment:
\smallskip

{\narrower\noindent
Furthermore, because the solitons are remarkably stable entities,
preserving their identities throughout numerous interactions, one
would expect this system to exhibit thermalization (complete
energy sharing among the corresponding linear normal modes)
only after extremely long times, if ever. \par}
\smallskip

\noindent
But this explanation, elegant as it may be, only pushes the basic question
back a step. A full understanding of FPU recurrence requires that we
comprehend the reasons behind the remarkable new phenomenon of solitonic
behavior, and in particular {\it why\/} solitons preserve their shape.
In fact, it was quickly recognized that the soliton was itself
a vital new feature of nonlinear dynamics, so that
understanding it better and discovering other nonlinear wave equations
that had soliton solutions became a primary focus for research in both
pure and applied mathematics. The mystery of the FPU Lattice recurrence
soon came to be regarded as an important but fortuitous spark that ignited
this larger effort.

\goodbreak
\smallskip
  The next few short sections explain some elementary but important facts about 
one-dimensional wave equations. If you know about shock development, and how dispersion 
smooths shocks, you can skip these sections without loss of continuity.
\goodbreak

\subsection {``Steepening'' and ``Breaking''}

\noindent
   Several times already we have referred to the phenomenon of ``steepening and
breaking of negatively sloped wave-fronts'' for certain wave equations. If you
have never seen this explained it probably sounds suggestive but also a little 
mysterious. In fact something very simple is going on that we will now explain.

  Let us start with the most elementary of all one-dimensional wave equations,
the linear advection equation (or forward wave equation), $u_t + c u_x = 0$. 
If we think of the graph of $x \mapsto u(x,t)$ as representing the profile of 
a wave at time $t$, then this equation describes a special evolutionary behavior 
of the wave profile in time. In fact, if $u_0(x) = u(x,0)$ is the ``initial'' 
shape of the wave, then the unique solution of the equation with this initial
condition is the so-called ``traveling wave'' $u(x,t) = u_0(x - ct)$, 
i.e., just the initial wave profile translating rigidly to the right at a
uniform velocity $c$. In other words, we can construct the wave profile at
time $t$ by translating each point on the graph of $u_0(x)$ horizontally
by an amount $ct$. As we shall now see, this has a remarkable generalization.

  We shall be interested in the non-viscous Burgers' equation, $u_t + u u_x = 0$,
but it is just as easy to treat the more general equation $u_t + f(u) u_x = 0$,
where $f:\reals \to \reals$ is some smooth function. Let me call this simply the
nonlinear advection equation or NLA. 

\proclaim {Proposition}  Let $u(x,t)$ be a smooth solution of the nonlinear advection 
equation $u_t + f(u) u_x = 0$ for $x \in \reals$ and $t \in [0,t_0]$, and with 
initial condition $u_0(x) = u(x,0)$. Then for $t < t_0$ the graph of $x \mapsto u(x,t)$
can be constructed from the graph of $u_0$ by translating each point $(x,u_0(x))$
horizontally by an amount $f(u_0(x)) t$.\endproclaim
\proof The proof is by the ``method of characteristics'', i.e., we look for curves
$(x(s),t(s))$ along which $u(x,t)$ must be a constant (say $c$), because $u$ satisfies NLA. 
If we differentiate $u(x(s),t(s)) = c$ with respect to $s$, then the chain rule gives
$u_x(x(s),t(s)) x^\prime(s) + u_t(x(s),t(s)) t^\prime(s) = 0$, and hence
$d x/dt =  x^\prime(s)/t^\prime(s) = - u_t(x(s),t(s))/u_x(x(s),t(s))$, and now
substitution from NLA gives:
$$dx/dt = f(u(x(s),t(s))) = f(c),$$ 
so the characteristic curves are straight lines, whose slope is $f(c)$, where $c$ is
the constant value the solution $u$ has along that line. In particular, if we
take the straight line with slope $f(u_0(x))$ starting from the point $(x,0)$,
then $u(x,t)$ will have the constant value $u_0(x)$ along this line, a fact
that is equivalent to the conclusion of the Proposition.
\qed

  It is now easy to explain steepening and breaking. We assume that the 
function $f$ is monotonically increasing and that $u_0(x)$ has negative
slope (i.e., is strictly decreasing) on some interval $I$. If we follow
the part of the wave profile that is initially over the interval $I$, we 
see from the Proposition that the higher part (to the left) will move faster
than the lower part (to the right), and so gradually overtake it. 
The result is that the wave ``bunches up'' and its slope increases---this is 
steepening---and eventually there will be a first time $T_B$ when the graph
has a vertical tangent---this is breaking. Clearly the solution cannot be
continued past $t = T_B$, since for $t > T_B$ the Proposition would give a
multi-valued graph for $u(x,t)$. It is an easy exercise to show that the
breaking time $T_B$ is given by $|\min(u_0^\prime(x))|^{-1}$.

  This explains the first part of the above quotation from Zabusky and Kruskal,
namely, ``Initially the first two terms of Eq. (1) dominate and
the classical overtaking phenomenon occurs; that is $u$ steepens
in regions where it has negative slope.'' But what about their next comment:
``Second, after $u$ has steepened sufficiently, the third term becomes important 
and serves to prevent the formation of a discontinuity.''? To explain this we have 
to take up the matter of dispersion.

\subsection {Dispersion}

\noindent
  Let us next consider {\it linear\/} wave equations of the form
$u_t + P\left ({\partial\over \partial x}\right ) u = 0$, where $P$ is a polynomial.
Recall that a solution $u(x,t)$ of the form $e^{i(kx-\omega t)}$ is 
called a plane-wave solution; $k$ is called the {\it wave number\/}
(waves per unit length) and $\omega$ the (angular) {\it frequency\/}. Rewriting
this in the form $e^{ik(x-(\omega/k) t)}$, we recognize that this is
a traveling wave of velocity $\omega\over k$. If we substitute this $u(x,t)$
into our wave equation, we get a formula determining a unique frequency 
$\omega(k)$ associated to any wave number $k$, which we can write in the
form ${\omega(k)\over k} = {1\over ik}P(ik)$. This is called the 
``dispersion relation'' for this wave equation. Note that it expresses
the velocity for the plane-wave solution with wave number $k$. For example,
$P\left({\partial\over \partial x}\right) = c {\partial\over \partial x}$ gives
the linear advection equation $u_t + c u_x = 0$, which has the dispersion
relation ${\omega(k)\over k} = c$, showing of course that all plane-wave
solutions travel at the same velocity $c$, and we
say that we have trivial dispersion in this case. On the other hand, if we take 
$P\left ({\partial\over \partial x}\right ) = 
\left({\partial\over \partial x}\right)^3$, then our wave equation is $u_t + u_{xxx} = 0$,
which is the KdV equation without its nonlinear term, and we have the
non-trivial dispersion relation ${\omega(k)\over k} = -k^2$. In this case, 
plane waves of large wave-number (and hence high frequency) are traveling 
much faster than low-frequency waves. The effect of this 
is to ``broaden a wave-packet''. That is, suppose our initial condition is $u_0(x)$. 
We can use the Fourier Transform to write $u_0$ in the form
$u_0(x) = \int \hat u_0(k) e^{ikx}\,dk$, and then, by superposition,
the solution to our wave equation will be
$u(x,t) = \int \hat u_0(k) e^{ik(x - (\omega(k)/k)t)}\,dk$.
Suppose for example our initial wave form is a highly peaked Gaussian.
Then in the case of the linear advection equation all the Fourier
modes travel together at the same speed and the Gaussian lump remains 
highly peaked over time. On the other hand, for the linearized
KdV equation the various Fourier modes all travel at different
velocities, so after a short time they start cancelling each other
by destructive interference, and the originally sharp Gaussian
quickly broadens. This is what Zabusky and Kruskal are referring to
when they say that ``\dots the third term becomes important and serves to 
prevent the formation of a discontinuity.'' Just before breaking or
shock-formation, the broadening effects of dispersion start to cancel
the peaking effects of steepening. Indeed, careful analysis shows that 
in some sense, what gives KdV solitons their special properties of 
stability and longevity is a fine balance between the yin effects of
dispersion and the yang effects of steepening.

\subsection {Split-stepping KdV}

\noindent

There is an interesting question that is suggested by our
analysis in the last two sections. In the KdV equation,
$u_t = -6uu_x - u_{xxx}$, 
if we drop the nonlinear term,
we have a constant coefficient linear PDE  whose initial value
problem can be solved explicitly by the Fourier Transform. 
On the other hand, if we ignore the linear third-order term, then
we are left with the inviscid Burgers' equation, whose initial value
problem can be solved numerically by a variety of methods. 
(It can also be solved in implicit form analytically, for short times,
by the method of characteristics, 
$$u=u_o(x-6ut)$$
but the solution is not conveniently represented on a fixed numerical grid.)
Can we somehow combine the methods for solving each of the two parts into 
an efficient numerical method for solving the full KdV initial value problem?

  In fact we can, and indeed there is a very general technique that
applies to such situations. In the pure mathematics community it
is usually referred to as the Trotter Product Formula, while in
the applied mathematics and numerical analysis communities it is
called split-stepping. Let me state it in the context of
ordinary differential equations. Suppose that $Y$ and $Z$ are two
smooth vector fields on $\reals^n$, and we know how to solve each of
the differential equations $dx/dt = Y(x)$ and $dx/dt = Z(x)$, meaning
that we know both of the flows $\phi_t$ and $\psi_t$ on $\reals^n$
generated by $X$ and $Y$ respectively. The Trotter Product Formula is
a method for constructing the flow $\theta_t$ generated by $Y + Z$
out of $\phi$ and $\psi$; namely, letting $\Delta t = {t \over n}$,
$\theta_t = \lim_{n\to\infty} (\phi_{\Delta t}\psi_{\Delta t})^{n}$.
The intuition behind the formula is simple. Think of approximating
the solution of $dx/dt = Y(x) + Z(x)$ by Euler's Method. If we are
currently at a point $p_0$, to propagate one more time step $\Delta t$
we go to the point $p_0 + \Delta t(Y(p_0) + Z(p_0))$.
Using the split-step approach on the other hand,
we first take an Euler step in the $Y(p_0)$ direction,
going to $p_1 = p_0 + \Delta t Y(p_0)$, then take a second
Euler step, but now from $p_1$ and in the $Z(p_1)$ direction, going
to $p_2 = p_1 + \Delta t Z(p_1)$.  If $Y$ and $Z$ are
constant vector fields then this gives exactly the same final result
as the simple full Euler step with $Y + Z$, while for continuous
$Y$ and $Z$ and small time step $\Delta t$ it is a good enough
approximation that the above limit is valid.

  The situation is more delicate for flows on infinite dimensional
manifolds, nevertheless it was shown  by F.~Tappert in [Ta] that the
the Cauchy Problem for KdV can be solved numerically  by using
split-stepping to combine solution methods for $u_t = -6uu_x$ and
$u_t = - u_{xxx}$. In addition to providing a perspective on an evolution
equation's relation to its component parts, split-stepping allows one to 
modify a code from solving KdV to the Kuramoto-Sivashinsky equation 
($u_t+uu_x=-u_{xx}-u_{xxxx}$), or study the joint zero-diffusion-dispersion
limits KdV-Burgers' equation ($u_t+6uu_x=\nu u_{xx}+\epsilon u_{xxxx}$),
by merely changing one line of code in the Fourier module.

  Tappert uses an interesting variant, known as Strang splitting,
which was first suggested in [St] to solve multi-dimensional
hyperbolic problems by split-stepping one-dimensional problems.
The advantage of splitting comes from the greatly reduced effort required
to solve the smaller bandwidth linear systems which arise when implicit 
schemes are necessary to maintain stability.
In addition, Strang demonstrated that second-order accuracy of the
component methods need not be compromised by the assymmetry of the
splitting, as long as the pattern 
$\phi_{\Delta t \over 2} \psi_{\Delta t \over 2} 
 \psi_{\Delta t \over 2} \phi_{\Delta t \over 2}$ 
is used, to account for possible non-commutativity of $Y$ and $Z$.
(This may be seen by multiplying the respective exponential series.)
No higher order analogue of Strang splitting is available. Serendipitously, 
when output is not required, several steps of Strang splitting
require only marginal additional effort: 
$(\phi_{\Delta t \over 2} \psi_{\Delta t \over 2} 
 \psi_{\Delta t \over 2} \phi_{\Delta t \over 2})^n=
(\phi_{\Delta t \over 2} \psi_{\Delta t} 
(\phi_{\Delta t} \psi_{\Delta t} )^{n-1}
 \phi_{\Delta t \over 2} $

\subsection {A Symplectic Structure for KdV}

\noindent
  The FPU Lattice is a classical finite dimensional mechanical system,
and as such it has a natural Hamiltonian formulation. However its
relation to KdV is rather complex---and KdV is a PDE rather than a
finite  dimensional system of ODE---so it is not clear that it too
can be viewed as a Hamiltonian system. We shall now see how this can be
done  in a simple and natural way. Moreover, when interpreted as the
infinite dimensional analogue of a Hamiltonian system, KdV turns
out to have a key property one would expect from any  generalization
to infinite dimensions of the concept of complete integrability in the
Liouville sense; namely the existence of infinitely many functionally
independent constants of the motion that are in involution. (Later, in
discussing the inverse scattering  method, we will indicate how complete
integrability was proved in a more precise sense  by Fadeev and Zakharov[ZF];
they demonstrated that the ``scattering data'' for the KdV equation obey
the characteristic  Poisson bracket relations for the action-angle
variables of a  completely integrable system.)

  In 1971, Gardiner and Zakharov independently showed how to
interpret KdV as a Hamiltonian system, starting from a Poisson
bracket approach, and from this beginning Poisson brackets have played
a significantly more important r\^ole in the infinite dimensional theory
of Hamiltonian systems than they did in the more classical finite
dimensional theory, and in recent years this has led to a whole theory
of so-called Poisson manifolds and Poisson Lie groups. However, we will
start with the more classical approach to Hamiltonian systems, defining  
a symplectic structure for
KdV first and then obtain the Poisson bracket structure as a derived
concept (cf.\ Abraham and Marsden [AbM]). Thus, we will first exhibit a
symplectic structure $\Omega$ for the phase space $P$ of the KdV equation
and a Hamiltonian function, $H:P\to \reals$, such that the KdV equation
takes the form $\dot u = (\sgrad H)_u$.

  For simplicity, we shall take as our phase space $P$ the Schwartz
space, $\S(\reals)$, of rapidly decreasing functions $u: \reals \to \reals$,
although a much larger space would be possible. (In [BS] it is proved that
KdV defines a global flow on the Sobolev space $H^4(\reals)$
of functions $u: \reals \to \reals$ with derivatives of order up to
$4$ in $L^2$, and it is not hard to see that $P$ is an invariant
subspace of this flow. See also [Ka1], [Ka2].)
For $u,v$ in $P$ we will denote
their $L^2$ inner product $\int_{-\infty}^\infty u(x) v(x) \, dx$ by
$\IP<u,v>$ and we define
$$\Omega(u,v) =
{1\over 2}\int_{-\infty}^\infty
(v(x)\Int u(x) - u(x)\Int v(x)) \, dx,$$
where $\Int u(x) = \int_{-\infty}^x u(y)\, dy$ denotes the indefinite
integral of $u$. (For the periodic KdV equation we take $P$ to be all smooth
periodic functions of period $2\pi$ and replace the
$\int_{-\infty}^\infty$ by  $\int_0^{2\pi}$.)

  We denote by $\D$ the derivative operator, $u \mapsto u^\prime$,
so $\D\Int u = u$, and $\int_{-\infty}^\infty \D u =0$ for functions $u$
that vanish at infinity.
We will also write $\Der_k(u)$ for $\partial^k u$,
but for small $k$ we shall also use $u = \Der_0(u)$, $u_x = \Der_1(u)$,
$u_{xx} = \Der_2(u)$, etc.

There is a simple but important relation
connecting $\Omega$, $\D $, and the $L^2$ inner product; namely:
$$\Omega(\D u,v) =\IP<u,v>.$$
This is an immediate consequence of three obvious identities:
$\D (u\Int v) = (\D u)\,\Int v + u\, v$,
$\int_{-\infty}^\infty \D (u\Int v) = 0$, and
$\Omega(\D u,v) = (1/2)\int_{-\infty}^\infty (v\, u - (\D u)\Int v)$.

  One important consequence of this is the weak
non-degeneracy of $\Omega$. For, if $i_v\Omega$ is zero, then
in particular
$\IP<u,v> =
\Omega(\D u,v)=-\Omega(v,\D u)= -(i_v\Omega)(\D u) = 0$ for all $u$,
so $v = 0$.

  $\Omega$ is clearly a skew-bilinear form on $P$. Since $P$ is
a vector space, we can as usual identify $P$ with its tangent space
at every point, and then $\Omega$ becomes a ``constant'' $2$-form
on $P$. Since it is constant, of course $d\Omega = 0$.
(Below we will exhibit an explicit $1$-form $\omega$ on $P$
such that $d\omega = \Omega$.) Thus $\Omega$ is a symplectic
form for $P$, and henceforth we will consider $P$ to be a
symplectic manifold.

  A second consequence of $\Omega(\D u,v) =\IP<u,v>$ is that if
$F: P \to \reals$ is a smooth function (or ``functional'') on $P$
that has a gradient $\grad F$ with respect to the flat Riemannian
structure on $P$ defined by the $L^2$ inner product, then the
symplectic gradient of $F$ also exists and
is given by $(\sgrad F)_u = \D ((\grad F)_u)$.
Recall that $dF$, the differential of $F$, is the $1$-form on $P$
defined by
$$dF_u(v) = {d\over d\eps}\biggr |_{\eps = 0}F(u+\eps v),$$
and the gradient of $F$ is the vector field dual to
$dF$ with respect to the $L^2$ inner product (if such a vector
field indeed exists), i.e., it is
 characterized by $(dF)_u(v) = \IP<(\grad F)_u,v>$. Since
$\IP<(\grad F)_u,v> = \Omega((\D (\grad F)_u),v)$, it then follows
that $(\sgrad F)_u$ also exists and equals $\D ((\grad F)_u)$.

  We shall only consider functions $F : P \to \reals$ of the
type normally considered in the Calculus of Variations, i.e., of the
form:
$$F(u) = \int_{-\infty}^\infty \tilde F(u,u_x, u_{xx},\ldots)\, dx,$$
where $\tilde F :\reals^{k+1} \to \reals$ is a polynomial function
without a constant term. Then the usual integration by parts argument
of the Calculus of Variations shows that such an $F$ has a gradient, given by:
$$(\grad F)_u = {\partial \tilde F\over \partial u} -\,
							\D\biggl({\partial \tilde F\over \partial u_x}\biggr) +\, 
 \D^2\biggl({\partial \tilde F\over \partial u_{xx}}\biggl) -
                \ldots$$

\smallskip\noindent
{\bf Remark.} The above formula is written using the standard but
somewhat illogical conventions of the Calculus of Variations and
needs a little interpretation. $\tilde F$ is a function of variables
$y = (y_0, y_1, y_2,\ldots y_k)$, and for example
$\partial \tilde F/ \partial u_{xx}$
really means the function on $\reals$ whose value at $x$ is
$\partial \tilde F/ \partial y_2$ evaluated at
$y = (\Der_0(u)(x), \Der_1(u)(x),\Der_2(u)(x),\ldots \Der_k(u)(x))$.

\smallskip\noindent
From what we saw above, the symplectic
gradient of such an $F$ exists and is given by:
$$(\sgrad F)_u = \,\D \biggl({\partial \tilde F\over \partial u}\biggr) -\,
                \D ^2\biggl({\partial \tilde F\over \partial u_x}\biggr) +\, \D
^3\biggl({\partial \tilde F\over \partial u_{xx}}\biggr) -
                \ldots$$
Thus every such $F$ is a Hamiltonian function on $P$, defining
the Hamiltonian flow
$\dot u = (\sgrad F)_u$, where $u(t)$ denotes a smooth curve
in $P$. If instead of $u(t)(x)$ we write $u(x,t)$, this symbolic
ODE in the manifold $P$ becomes the PDE:
$$u_t =  \D \biggl({\partial \tilde F\over \partial u}\biggr) -\,
             \D ^2\biggl({\partial \tilde F\over \partial u_x}\biggr) +\, \D
^3\biggl({\partial \tilde F\over \partial u_{xx}}\biggr) -
                \ldots$$
and in particular if we take
$\tilde F(u,u_x) = -u^3 + u_x^2/2$ , then we get the KdV equation in
standard form:
$u_t = \D (-3u^2) - \D ^2(u_x) = -6u\,u_x - u_{xxx}$.

\noindent
{\bf Remark.} The formula defining $\Omega$ can be motivated as follows.
Define linear functionals
$p_x$ and $q_x$ on $P$ by $q_x(u) = u(x)$ and $p_x(u) = \Int u(x)$.
(Think of these as providing ``continuous coordinates''
for $P$.) These give rise to differential $1$-forms $dp_x$
and $dq_x$ on $P$. Of course, since $p_x$ and $q_x$ are linear,
at every point $u$ of $P$, we have $dp_x = p_x$ and $dq_x = q_x$.
Then $\Omega$ can now be written in the suggestive form
$\Omega = \sum_x dp_x\wedge dq_x$, where $\sum_x$ is shorthand for
$\int_{-\infty}^\infty$. This suggests that we define a $1$-form
$\omega$ on $P$ by $\omega = \sum_x p_x\, dq_x$, i.e.,
$\omega_w(u) = \int_{-\infty}^\infty \Int w(x)u(x)\,dx$.
Consider this as a function $f(w)$ on $P$ and let us compute
its directional derivative at $w$ in the direction $v$,
$(vf)(w) = {d\over d\eps} |_{\eps = 0}f(w+\eps v)$.
We clearly get $v(\omega(u))=
\int_{-\infty}^\infty \Int v(x)u(x)\,dx$.
Since $u$ and $v$ are constant vector fields, their bracket
$[u,v]$ is zero, and we calculate
$d\omega(u,v) = v(\omega(u)) - u(\omega(v)) = \Omega(u,v)$,
as expected.
\smallskip

 We now again specialize to the phase space $P$ for the KdV equation,
namely the Schwartz space $\S(\reals)$ with its $L^2$ inner product
$\IP<u,v>$ and symplectic form $\Omega(u,v)$, related by
$\Omega(\D u,v) = \IP<u,v>$. Then, since $\sgrad F = \D (\grad F)$,
we obtain the formula
$$\PB (F_1,F_2) = \Omega(\sgrad F_2, \sgrad F_1) =
   \Omega(\D\grad F_2, \D\grad F_1) =
    \IP<\grad F_2, \D (\grad F_1)>$$
for Poisson brackets in terms of the Riemannian structure for $P$,
and in particular we see that $F_1$ and $F_2$ are in involution
if and only if the two vector fields $\grad F_1$ and $\D\grad F_2$
on $P$ are everywhere orthogonal.

\section {The Inverse Scattering Method}

\noindent
  In 1967, in what would prove to be one of the most cited mathematical
papers in history, [GGKM], Clifford Gardner, John Greene, Martin
Kruskal, and Robert Miura introduced an ingenious method, called the
Inverse Scattering Transform (IST), for solving the KdV equation. In
the years that followed, the IST changed applied mathematics like
no other tool since the Fourier Transform (to which it is closely
related) and it soon became clear that it was the key to understanding
the remarkable properties of soliton equations.

  Before starting to explain the IST, we recall the basic philosophy
of using ``transforms'' to solve ODE\null. Suppose we are interested in
some evolution equation $\dot x = X(x)$ on a smooth manifold $M$.
That is, $X$ is a smooth vector field on $M$ that generates a
flow $\phi_t$ on $M$. Usually our goal is to understand the
dynamical properties of this flow---and perhaps get an explicit
``formula'' for $\phi_t(x)$, at least for some initial conditions $x$.
A {\it transform\/} is a diffeomorphism $T$ of $M$ onto some other
manifold $N$, mapping the vector $X$ onto a vector field $Y = DT(X)$ on
$N$. If $\psi_t$ is the flow generated by $Y$, then clearly
$T(\phi_t(x)) = \psi_t(Tx)$, and it follows that if we understand
$\psi_t$  well, and moreover have explicit methods for computing
$T(x)$ and $T^{-1}(y)$, then we in effect also know all about $\phi_t$.

  It is important to realize that there is usually more at stake
than just finding particular solutions of the original initial value
problem. Essential structural features of the flow that are hidden
from view in the original form of the evolution equation may become
manifest when viewed in the transform space $N$.

  For example, consider the case of a {\it linear\/} evolution equation
$\dot x = X(x)$ on some vector space $M$. We can formally
``solve'' such an equation in the form $x(t) = \exp(tX) x(0)$.
However, explicit evaluation of the linear operator
$\exp(tX)$ is not generally feasible, nor does the
formula provide much insight into the structure of the flow.
But suppose we can find a linear diffeomorphism $T:M \to N$ so that
the linear operator $Y = T X T^{-1}$ is diagonal in some ``basis''
(discrete or continuous) $\{w_\alpha\}$ for $N$, say
$Y w_\alpha = \lambda_\alpha w_\alpha$.
Then $\exp(tY) w_\alpha = e^{{\lambda_\alpha}t} w_\alpha$, hence if
$y(0) = \sum_\alpha \hat y_\alpha w_\alpha$ then the solution
to the initial value problem $\dot y = Y(y)$ with initial value $y(0)$
is $y(t) = \sum_\alpha (e^{{\lambda_\alpha}t}\hat y_\alpha) w_\alpha$.
Not only do we have an explicit formula for $\psi_t$, but we see
the important structural fact that the flow is just a direct sum
(or integral) of uncoupled one-dimensional flows, something not
obvious when viewing the original flow.

  This is precisely why the Fourier transform is such a powerful tool
for analyzing constant coefficient linear PDE---it simultaneously
diagonalizes all such operators! Since the Fourier transform is an
excellent model for understanding the more complex IST, let us quickly
review it in our current context. It will be convenient to complexify
$P$  temporarily, i.e., regard our phase space as the complex
vector space of complex-valued Schwartz functions on $\reals$.
Then the Fourier Transform, $v \mapsto  w = \FT(v)$, is a linear
diffeomorphism of $P$ with  itself, defined by
$w(\alpha) =
{1\over \sqrt{2\pi}}\int_{-\infty}^\infty v(x) e^{-i \alpha x}\, dx$,
and the Inverse Fourier Transform, $w \mapsto v =\IFT(w)$ is given by
$v(x) =
{1\over \sqrt{2\pi}}
   \int_{-\infty}^\infty w(\alpha) e^{i \alpha x}\, d\alpha$.

Given any $n+1$-tuple of real numbers $a = (a_0,\ldots a_n)$, we
let $F_a(y)$ denote the polynomial
$a_0 y + a_1 y^3 + \ldots + a_n y^{2n+1}$, and $F_a(\D)$ the constant
coefficient linear differential operator
$a_0 \D + a_1 \D^3 + \ldots + a_n \D^{2n+1}$.
Note that $F_a(\D)$ is a vector field on $P$.
In fact, if we put
$\tilde H_a(\Der_0(v),\ldots, \Der_n(v)) =
{1\over 2}\sum_{j=0}^n a_j({\Der_j(v)})^2$, and define the
corresponding functional
$H_a(v) =
\int_{-\infty}^\infty \tilde H_a(\Der_0(v),\ldots, \Der_n(v))\, dx$,
then clearly $F_a(\D) = \sgrad H_a$.
It is trivial that if $b= (b_0,\ldots b_m)$
is some other $m+1$-tuple of real numbers then
$[F_a(\D), F_b(\D)] = 0$, i.e., all these differential operators
(or vector fields) on $P$ commute, and it is easy to check
directly that $\PB(H_a,H_b) = 0$, i.e., that
the corresponding Hamiltonian functions Poisson commute.

  The transform, $G_a$, of the vector field $F_a(\D)$
under the Fourier Transform is easy to compute:
$G_a(w)(\alpha) = F_a(i\alpha) w(\alpha)$, or in words,
the partial differential operator $F_a(\D)$ is transformed
by $\FT$ into multiplication by the function $F_a(i\alpha)$.
In ``physicist language'', this shows that the $G_a$ are
all diagonal in the continuous  basis for $P$ given by the
evaluations $w \mapsto w(\alpha)$.

  Before going on to consider the Scattering Transform we should
mention another classical and elementary transform---one linearizing
Burgers' Equation, $v_t = v_{xx} -2v v_x$. The transform,
$\CH$ mapping $v$ to $w$, is $w = \exp(-\Int{v})$,
and the inverse transform $\ICH$ that
recovers $v$ from $w$ is $v = -\D \log(w)= -\D w /w$.
Clearly $w$ must be positive for this to be defined, and it is
easily checked that if $w$ is a positive solution of
the linear heat conduction (or diffusion) equation
$w_t = w_{xx}$ then $v$ satisfies Burgers' Equation.
So if we start with any positive integrable function
$w(x,0)$, we can use the Fourier Transform method to find
$w(x,t)$ satisfying the heat equation, and then
$v(x,t) = -w_x(x,t)/w(x,t)$ will give a solution of
Burgers' Equation. ($\CH$ is usually referred to as the Cole-Hopf
Transform, but the fact that it linearizes Burgers' Equation was
actually pointed out by Forsyth in 1906, four decades before Cole
and Hopf each independently rediscovered it.)

\subsection {Lax Equations: KdV as an Isospectral Flow}

\noindent
  In discussing the Inverse Scattering Transform it will be
useful to have available an interesting reinterpretation of
the KdV equation as formulated by Peter Lax. Namely, if $u(x,t)$
is a solution of the KdV equation, and we consider the one-parameter
family $L(t)$ of self-adjoint operators on $L^2(\reals)$ that are given by
the Schr\"odinger operators with potentials $u(t)(x) = u(x,t)$
(i.e., $L(t) \psi (x) = -{d^2\over dx^2} \psi(x) + u(x,t) \psi(x)$),
then these operators are isospectral, and in fact unitarily equivalent.
That is, there is a smooth one parameter family $U(t)$ of unitary
operators on $L^2(\reals)$ such that $U(0) = I$ and
$L(t) = U(t) L(0) U(t)^{-1}$.

   By the way, in the following it will be convenient to take KdV in the form
$u_t -6 u u_x +  u_{xxx} = 0$.

  Suppose we have a smooth one-parameter family $U(t)$ of unitary
transformations of a Hilbert space $H$ with $U(0) = I$.
$U_t(t)$, the derivative of $U(t)$,
is a tangent vector at $U(t)$ of the group
${\Cal U}(H)$ of unitary transformations of $H$, so
$B(t) = U_t(t) U(t)^{-1} = U_t(t) U(t)^*$ is a tangent vector
to ${\Cal U(H)}$ at the identity, $I$.
Differentiating $U U^* = I$ gives $U_t U^* + U U_t^* = 0$,
and since $U_t = B U$ and $U_t^* = U^* B^*$,
$0 = BUU^* + UU^*B^*$, so $B^* = -B$, i.e., $B(t)$ is a family of
skew-adjoint operators on $H$. Conversely, a smooth map
$t \mapsto B(t)$ of $\reals$ into the skew-adjoint operators defines
a time-dependent right invariant vector field $X_U(t) = B(t) U$
on ${\Cal U}(H)$ and so (at least in finite dimensions)
a smooth curve $U(t)$ of unitary operators starting from $I$
such that $U_t(t) = B(t)U(t)$.

  Now suppose that $L(0)$ is a self-adjoint operator on $H$, and define
a family of conjugate operators $L(t)$ by $L(t) = U(t) L(0) U(t)^{-1}$,
so $L(0) = U(t)^* L(t) U(t)$. Differentiating the latter with respect
to $t$, $0 = U_t^*LU + U^*L_tU + U^*LU_t = U^*(-BL + L_t + LB)U$.
Hence, writing $[B,L] = BL - LB$ as usual for the commutator of $B$
and $L$, we see that $L(t)$ satisfies the so-called {\it Lax Equation\/},
$L_t = [B,L]$.

  Given a smooth family of skew-adjoint operators $B(t)$, the Lax Equation
is a time-dependent linear ODE in the vector space ${\Cal S}$ of
self-adjoint operators on $H$, whose special form expresses the
fact that the evolution is by unitary conjugation. Indeed, since the
commutator   of a skew-adjoint operator and a self-adjoint operator is
again self-adjoint, $B(t)$ defines a time-dependent vector field, $Y$,
on ${\Cal S}$ by $Y(t)(L) = [B(t),L]$. Clearly a smooth curve $L(t)$ in
${\Cal S}$ satisfies the Lax Equation if and only it is a solution curve
of $Y$.  By uniqueness of solutions of linear ODE, the  solution $L(t)$ of
this ODE with initial condition $L(0)$ must be the  one-parameter family
$U(t) L(0) U(t)^{-1}$ constructed above.

   Given any $\psi(0)$ in $H$, define $\psi(t) = U(t)\psi(0)$.
Since $U(t) L(0) =  L(t) U(t)$, it follows that if $\psi(0)$ is
an  eigenvector of $L(0)$ belonging to the eigenvalue $\lambda$,
then $\psi(t)$ is an eigenvalue of $L(t)$ belonging
to the {\it same\/} eigenvalue $\lambda$. Differentiating
the relation defining $\psi(t)$ gives $\psi_t = B \psi(t)$, so we may
consider $\psi(t)$ to be defined as the solution of this linear ODE
with initial value $\psi(0)$. Since this is one of the main ways in
which we will use Lax Equations, we will restate it as what we shall
call the:

\proclaim {Isospectral Principle} Let $L(t)$ and $B(t)$ be smooth
one-parameter families of self-adjoint and skew-adjoint operators
respectively on a Hilbert space $H$, satisfying the Lax Equation
$L_t = [B,L]$, and let $\psi(t)$ be a curve in $H$ that is a solution
of the time-dependent linear ODE $\psi_t = B\psi$.
If the initial value, $\psi(0)$, is an eigenvector of $L(0)$ belonging 
to an eigenvalue $\lambda$, then $\psi(t)$ is an eigenvector of $L(t)$
belonging to the same eigenvalue $\lambda$.
\endproclaim

\noindent
{\bf Remark.} There is a more general (but less precise) version
of the Isospectral Principle that follows by an almost identical
argument. Let $V$ be any topological vector space and $B(t)$ a
family of linear operators on $V$ such that the evolution equation
$U_t = BU$ is well-defined. This means that for each $\psi(0)$ in
$V$ there should exist a unique solution to the time-dependent linear
ODE $\psi_t(t) = B(t) \psi(t)$. The evolution operator $U(t)$ is of
course then defined by $U(t)\psi(0) = \psi(t)$, so $U_t = BU$. Then
clearly the conclusion of the Isospectral Principle still holds. That
is to say, if a smooth family of linear operators $L(t)$ on $V$
satisfies the Lax Equation $L_t =[B,L]$, then $U(t)L(0) =  L(t)U(t)$,
so if $L(0) \psi(0) = \lambda \psi(0)$ then
$L(t) \psi(t) = \lambda \psi(t)$.
\smallskip

  We now apply the above with $H = L^2(\reals)$. We will see
that if $u$ satisfies KdV, then the family of Schr\"odinger operators
$L(t)$ on $H$ defined above satisfies the Lax Equation $L_t = [B,L]$,
where
$$B(t)\psi(x) =
   -4\psi_{xxx}(x) + 3\left(u(x,t)\psi_x(x) + (u(x,t)\psi(x))_x\right),$$
or more succinctly, $B = -4\D^3 +3(u\D + \D u)$.
Here and in the sequel it is convenient to use the same symbol both
for an element $w$ of the Schwartz space, $\S(\reals)$, and for the
bounded self-adjoint multiplication operator $v \mapsto wv$ on $H$. Since
$H$ is infinite  dimensional and our operators $B$ and $L$ are unbounded on
$H$, some care is needed for a rigorous treatment. But this is relatively
easy. Note that all the operators involved have the Schwartz space as
a common dense domain, so we can use the preceding remark taking
$V = \S(\reals)$ (we omit details).

  Note that since $\D$ is skew-adjoint, so is any odd power, and in
particular $4\D^3$ is skew-adjoint. Also, the multiplication operator
$u$ is self-adjoint, while the anti-commutator of a self-adjoint and a
skew-adjoint operator is skew-adjoint, so $u\D + \D u$ and hence $B$ is
indeed skew-adjoint.

  Since clearly $L_t = u_t$, while $u_t - 6uu_x + u_{xxx} = 0$ by
assumption, to prove that $L_t = [B,L]$ we must check that
$[B,L] = 6uu_x -u_{xxx}$. Now
$[B,L] = 4[\D^3,\D^2] - 4[\D^3,u]
         -3[u \D,\D^2] +3[u \D,u]
         -3[\D u,\D^2] +3[\D u,u]$,
and it easy to compute the six commutators relations
$[\D^3,\D^2] =0$,
$[\D^3,u] = u_{xxx}+3u_{xx}\D +3u_x\D^2$,
$[u\D,\D^2] =-u_{xx}\D-2u_x\D^2$,
$[u\D,u] =uu_x$,
$[\D u,\D^2]= -3u_{xx}\D -2u_x\D^2-u_{xxx}$, and
$[\D u,u] = -uu_x$,
from which the desired expression for $[B,L]$ is immediate.

  Let us now apply the Isospectral Principle to this example.

\proclaim {KdV Isospectrality Theorem}
  Suppose $u(x,t)$ is a solution of the KdV equation,
$$u_t -6u u_x + u_{xxx} = 0,$$
whose initial value $u(x,0)$ is in the Schwartz space $\S(\reals)$,
and that $\psi(x)$ is an eigenfunction of the
Schr\"odinger Equation with potential $u(x,0)$ and eigenvalue $\lambda$:
$$ -{d^2\over dx^2} \psi(x) + u(x,0) \psi(x) = \lambda \psi(x).$$
Let $\psi(x,t)$ be the solution of the evolution equation $\psi_t = B\psi$,
i.e.,
$${\D \psi\over \D t} = -4{\D^3 \psi\over \D x^3} +
3\Bigl(u(x,t){\D \psi\over \D x}(x,t) +
   {\D \over \D x}\left(u(x,t)\psi(x,t)\right)\Bigr)$$
with the initial value $\psi(x,0) = \psi(x)$. Then $\psi(x,t)$ is an
eigenfunction for the Schr\"odinger Equation with potential $u(x,t)$
and the same eigenvalue $\lambda$:
$$ -\psi_{xx}(x,t) + u(x,t) \psi(x,t) = \lambda \psi(x,t),$$
and moreover, if $\psi(x)$ is in $L^2$, then the $L^2$ norm of
$\psi(\cdot,t)$ is independent of $t$. Finally, $\psi(x,t)$ also satisfies
the first-order evolution equation
$$\psi_t - (4\lambda + 2u)\psi_x + u_x\psi = 0.$$
\endproclaim

\proof Except for the final statement this is an immediate
application of the Isospectrality Principle. Differentiating the
eigenvalue equation for $\psi(x,t)$ with respect to $x$ gives
$\psi_{xxx} = u_x\psi + (u-\lambda)\psi_x$, and substituting this
into the assumed evolution equation for $\psi$ gives the asserted
first-order equation for $\psi$. \qed

  By the way, it should be emphasized that the essential point
is that when a potential evolves via KdV then the corresponding
Schr\"odinger operators are isospectral, and this is already clearly
stated in [GGKM]\null.	 Lax's contribution was to explain the mechanism
behind this remarkable fact, and to formulate it in a way that
was easy to generalize. In fact, almost all generalizations of
the phenomena first recognized in KdV have used the Lax Equation
as a jumping off place.

\subsection {The Scattering Data and its Evolution}

\noindent
  We now fix a ``potential function'' $u$ in the Schwartz space
$\S(\reals)$ and look more closely at the space $E_\lambda(u)$
of $\lambda$ eigenfunctions of the Schr\"odinger operator with this
potential. By definition,  $E_\lambda(u)$ is just the kernel of the
linear operator
$ L^u(\psi) = -{d^2 \psi\over dx^2} + u\psi - \lambda \psi$
acting on the space $C^\infty(\reals)$, and by the elementary theory
of second-order linear ODE it is, for each choice of $\lambda$,
a two-dimensional linear subspace of $C^\infty(\reals)$.
Using the special form of $L^u$ we can describe $E_\lambda(u)$
more precisely. We will ignore the case $\lambda = 0$, and consider
the case of positive and negative $\lambda$ separately.

   Suppose $\lambda = - \kappa^2$, $\kappa > 0$. Note that any
$\psi$ in $E_\lambda(u)$ will clearly be of the form
$\psi(x) = a e^{\kappa x} + b e^{-\kappa x}$
in any interval on which $u$ vanishes identically. Thus if
$u$ has compact support, say $u(x) = 0$ for $|x| > M$,
then we can find a basis
$\psi_{\lambda,-\infty}^+,\psi_{\lambda,-\infty}^-$ for $E_\lambda(u)$
such that for $x < -M$,  $\psi_{\lambda,-\infty}^\pm(x) = e^{\pm\kappa x}$,
or equivalently
$\psi_{\lambda,-\infty}^+(x)e^{-\kappa x} = 1$ and
$\psi_{\lambda,-\infty}^-(x)e^{\kappa x} = 1$ for $x < -M$.
Similarly there is a second basis
$\psi_{\lambda,\infty}^+,\psi_{\lambda,\infty}^-$ for $E_\lambda(u)$
such that
$\psi_{\lambda,\infty}^+(x)e^{-\kappa x} = 1$ and
$\psi_{\lambda,\infty}^-(x)e^{\kappa x} = 1$ for $x > M$.
When $u$ does not have compact support but is only
rapidly decreasing then it can be shown that there still exist
two bases $\psi_{\lambda,-\infty}^+,\psi_{\lambda,-\infty}^-$
and $\psi_{\lambda,\infty}^+,\psi_{\lambda,\infty}^-$ for $E_\lambda(u)$
such that
$\lim_{x \to -\infty}\psi_{\lambda,-\infty}^+(x)e^{-\kappa x} = 1$ and
$\lim_{x \to -\infty}\psi_{\lambda,-\infty}^-(x)e^{\kappa x} = 1$,
while
$\lim_{x \to \infty}\psi_{\lambda,\infty}^+(x)e^{-\kappa x} = 1$ and
$\lim_{x \to \infty}\psi_{\lambda,\infty}^-(x)e^{\kappa x} = 1$.
(A more descriptive way of writing these limits is
$\psi_{\lambda,-\infty}^+(x) \sim e^{\kappa x}$ and
$\psi_{\lambda,-\infty}^-(x) \sim e^{-\kappa x}$ as $x \to -\infty$,
while
$\psi_{\lambda,\infty}^+(x) \sim e^{\kappa x}$ and
$\psi_{\lambda,\infty}^-(x) \sim e^{-\kappa x}$ as $x \to \infty$.)
Let us define functions $f(\lambda)$ and $c(\lambda)$ by
$\psi_{\lambda,-\infty}^+ =
  f(\lambda)\psi_{\lambda,\infty}^+  + c(\lambda) \psi_{\lambda,\infty}^-$.
Using these bases it is easy to detect when $\lambda$ is a so-called
``discrete eigenvalue'' of $L^u$, i.e., when
$E_\lambda(u)$ contains a non-zero element $\psi$ of $L^2(\reals)$.
We can assume $\psi$ has $L^2$ norm one, and since
$\psi_{\lambda,-\infty}^-$ blows up at $-\infty$ while
$\psi_{\lambda,\infty}^+$ blows up at $\infty$, $\psi$ must
be both a multiple of $\psi_{\lambda,-\infty}^+$
and of $\psi_{\lambda,\infty}^-$, and since $\psi \ne 0$
it follows that $f(\lambda) = 0$. Conversely, if $f(\lambda) = 0$ then
$\psi_{\lambda,-\infty}^+ = c(\lambda) \psi_{\lambda,\infty}^-$
decays exponentially both at $\infty$ and $-\infty$ and so
we can normalize it to get an element of $E_\lambda(u)$ with $L^2$
norm one. Thus the discrete eigenvalues of $L^u$ are precisely
the roots of the function $f$.

  It follows from standard arguments of Sturm-Liouville theory
that in fact $L^u$ has only finitely many discrete eigenvalues,
$\lambda_1,\ldots,\lambda_N$, with corresponding $L^2$ normalized
eigenfunctions $\psi_1,\ldots,\psi_N$, and these determine
so-called ``normalization constants'' $c_1,\ldots,c_N$ by
$\psi_n = c_n \psi_{\lambda_n,\infty}^-$, i.e., if we write
$\lambda_n = -\kappa_n^2$, then $c_n$ is characterized by
$\psi_n(x) \sim c_n e^{-\kappa_n x}$ as $x \to \infty$.
We note that the $\psi_n$ and hence the normalization constants
$c_n$ are only determined up to sign, but we will only use
$c_n^2$ in the Inverse Scattering Transform.

  For $\lambda = k^2$, $k>0$ there are similar considerations.
In this case if $u(x)$ vanishes for $|x| > M$ then any element
of $E_\lambda(u)$ will be of the form $a e^{ikx} + b e^{-ikx}$
for $x < -M$ and also of the form $c e^{ikx} + d e^{-ikx}$ for
$x > M$. If $u$ is only rapidly decaying then we can still find
bases $\psi_{\lambda,-\infty}^+,\psi_{\lambda,-\infty}^-$
and $\psi_{\lambda,\infty}^+,\psi_{\lambda,\infty}^-$ for $E_\lambda(u)$
such that
$\psi_{\lambda,-\infty}^+(x) \sim e^{ikx}$ and
$\psi_{\lambda,-\infty}^-(x) \sim e^{-ikx}$ as $x \to -\infty$,
while
$\psi_{\lambda,\infty}^+(x) \sim e^{ikx}$ and
$\psi_{\lambda,\infty}^-(x) \sim e^{-ikx}$ as $x \to \infty$. Then
$\psi_{\lambda,-\infty}^- =
   \alpha \psi_{\lambda,\infty}^- + \beta \psi_{\lambda,\infty}^+$, where
$\alpha$ can be shown to be non-zero. Dividing by $\alpha$ we get a
particular eigenfunction $\psi_k$, called the Jost solution,
with the special asymptotic behavior
$\psi_k(x) \sim a(k) e^{-ikx}$ as $x \to -\infty$
and $\psi_k(x) \sim e^{-ikx} + b(k) e^{ikx}$ as $x \to \infty$.

  The functions $a(k)$ and $b(k)$ are called the transmission coefficient
and reflection coefficient respectively, and $b(k)$ together with the
above normalizing constants $c_1,\ldots c_n$ make up the
``Scattering Data'', $\S(u)$ for $u$.

  While it is perhaps intuitively clear that the bases
$\psi_{\lambda,\pm\infty}^\pm$ must exist, to supply the
asymptotic arguments required for a rigorous proof of the
crucial theorem on the time evolution of the Scattering
Data it is essential to give them precise definitions, and
we do this next.

  First consider the simpler problem of the first order
ODE $L^u\psi = {d\psi\over dx} - u\psi$. If we make
the substitution $\psi = e^{\lambda x} \phi$, then
the eigenvalue equation $L^u(\psi) = \lambda \psi$ becomes
${d\phi\over dx} = u\phi$, so (assuming $u$ depends on a
parameter $t$) we have
$\phi(x,t) = \exp \bigl(\int_{-\infty}^x u(\xi,t)\,d\xi\bigr)$.
Note that
$\lim _{x \to -\infty} \phi(x,t) = 1$
while
$$\lim _{x \to \infty} \phi(x,t) =
\exp \Bigl(\int_0^\infty u(\xi,t)\,d\xi\Bigr) = c(t),$$
so if $\psi(x,t)$ is an eigenfunction of $L^u$,
$\psi(x,t) \sim c(t) e^{\lambda x}$
(i.e., $\lim _{x \to \infty}\psi(x,t) e^{-\lambda x} = c(t)$)
and since $u(x,t)$ is rapidly decaying we can moreover
differentiate under the integral sign to obtain
$\psi_t(x,t) \sim c^\prime(t) e^{\lambda x}$.
One can not differentiate asymptotic relations in general
of course, and since we will need a similar relation for
eigenfunctions of Schr\"odinger operators we must make a
short detour to justify it by an argument similar to the above.

  If we now make the substitution $\psi = \phi e^{-\kappa x}$ in
the eigenvalue equation $\psi_{xx} = \kappa^2 \psi + u\psi$, then
we get after simplifications $\phi_xx - 2\kappa\phi_x = u\phi$,
or $\D(\D - 2\kappa)\phi = u\phi$. Recall the method of solving
the inhomogeneous equation $\D(\D - 2\kappa)\phi = f$ by
``variation of parameters''. Since $1$ and $e^{2\kappa x}$
form a basis for the solutions of the homogeneous equation,
we look for a solution of the form
$\phi = \Theta_1 + \Theta_2 e^{2\kappa x}$, and to make the
system determined we add the relation
$\Theta_1^\prime + \Theta_2^\prime e^{2\kappa x} = 0$.
This leads to the equations
$\Theta_1^\prime = -{f\over 2\kappa}$ and
$\Theta_2^\prime = {f\over 2\kappa} e^{2\kappa x}$
so
$\phi = -{1\over 2\kappa} \int_0^x f(\xi)\, d\xi
       + {e^{2\kappa x}\over 2\kappa}
          \int_0^x f(\xi)e^{-2\kappa x}\, d\xi$.
If we now take $f = u\phi$ (and use $\phi e^{-\kappa x} = \psi$)
then we get the relation
$\phi(x,t) = {1\over 2\kappa}\int_x^0 u(\xi,t)\phi(\xi,t)\, d\xi
     - {e^{2\kappa x}\over 2\kappa}
          \int_x^0 u(\xi,t) \psi(\xi,t)e^{-\kappa x}\, d\xi$.
Assuming that $-\kappa^2$ is a discrete eigenvalue, and that
$\psi$ has  $L^2$ norm $1$, $u\psi$ will also be in $L^2$ and we
can estimate the second integral using the Schwartz Inequality,
and we see that in fact
$|\int_x^0 u(\xi) \psi(\xi)e^{-\kappa x}\, d\xi| < O(e^{-\kappa x})$,
so the second term is $O(e^{\kappa x})$.
It follows that $\psi(x,t) \sim c(t) e^{\kappa x}$ in the sense that
$\lim _{x \to -\infty}\psi(x,t) e^{-\kappa x} = c(t)$,
where
$c(t) = \phi(-\infty,t) =
{1\over 2\kappa}\int_{-\infty}^0 u(\xi,t)\phi(\xi,t)\,d\xi$.
In other words, the normalizing constant is well defined.
But what is more important, it also follows that
if $u(x,t)$  satisfies KdV, then the normalizing
constant $c(t)$ for a fixed eigenvalue $-\kappa^2$ is a
differentiable function of $t$ and satisfies
$\psi_t(x,t) \sim c^\prime(t) e^{\kappa x}$. This follows
from the fact that we can differentiate the formula for
$c(t)$ under the integral sign because $u$ is rapidly decreasing.
Note that differentiating the relation $\psi e^{\kappa x} = \phi$
gives $\psi_x e^{\kappa x} = \phi_x  - \kappa \psi$. But the formula
for $\phi$ shows that $\phi_x$ converges to zero at $-\infty$,
so $\psi_x(x,t) \sim -\kappa c(t)e^{\kappa x}$.
From the KdV Isospectrality Theorem, we know that if $u(x,t)$
satisfies KdV, then $\psi(x,t)$ satisfies
$\psi_t - (-4\kappa^2 + 2u)\psi_x + u_x\psi = 0$, so
the left hand side times $e^{\kappa x}$ converges
to $c^\prime(t) +4\kappa^2(-\kappa c(t))$ as $x \to \infty$
and hence $c^\prime(t) - 4\kappa^3c(t) = 0$, so
$c(t) = c(0)e^{4\kappa^3 t}$.

  By a parallel argument (which we omit) it follows that the
transmission and reflection coefficients are also well defined
and that the Jost solution $\psi_k(x,t)$ satisfies
$(\psi_k)_t  \sim a_t(k,t) e^{-ikx}$  at $-\infty$
and
$(\psi_k)_t \sim  b_t (k,t) e^{ikx}$ at $\infty$, and
then one can show from the KdV Isospectrality Theorem
that the transmission coefficients are constant, while the
reflection coefficients satisfy $b(k,t) = b(k,0) e^{8ik^3t}$.

\proclaim  {Theorem on Evolution of the Scattering Data}
   Let $u(t) = u(x,t)$ be a smooth curve in $\S(\reals)$ satisfying
the KdV equation $u_t -6u u_x + u_{xxx} = 0$ and assume that the
Schr\"odinger operator with potential $u(t)$ has discrete eigenvalues
$-\kappa_1^2,\ldots,-\kappa_N^2$ whose corresponding normalized
eigenfunctions have normalization constants $c_1(t),\ldots,c_n(t)$.
Let the transmission and reflection coefficients of $u(t)$
be respectively $a(k,t)$ and $b(k,t)$. Then the transmission
coefficients are all constants of the motion, i.e., $a(k,t) = a(k,0)$,
while the Scattering Data, $c_n(t)$ and $b(k,t)$, satisfy:
\roster
\item "\rom{1)}" $c_n(t) = c_n(0)e^{4\kappa_n^3 t}$,
\item "\rom{2)}" $b(k,t) = b(k,0)e^{8ik^3 t}$. 
\endroster
\endproclaim

\noindent
We note a striking (and important) fact: not only do we now have an
explicit and simple formula for the evolution of the scattering
data $\S(u(t))$ when $u(t)$ evolves by the KdV equation, but further
{\bf this formula does not require any knowledge of $u(t)$.}

  The fact that the transmission coefficients $a(k)$ are constants
of the motion while the logarithms of the reflection coefficients,
$b(k)$ vary linearly with time suggest that perhaps they can somehow
be regarded as action-angle variables for the KdV equation, thereby
identifying KdV as a completely integrable system in a precise
sense. While $a(k)$ and $b(k)$ are not themselves canonical
variables, Zakharov and Fadeev in [ZF] showed that certain
functions of $a$ and $b$ did satisfy the Poisson commutation
relations for action-angle variables. Namely, the functions
$p(k) = (k/\pi)\log |a(k)|^2 = (k/\pi)\log [1 + |b(k)|^2]$
and $q(k) = \arg( b(k))$ satisfy
$\PB({p(k)},{q(k^\prime)}) = \delta (k - k^\prime)$ and
$\PB({p(k)},{p(k^\prime)}) = \PB({q(k)},{q(k^\prime)}) = 0$.

  The above formula for the evolution of the Scattering Data is
one of the key ingredients for The Inverse Scattering Method, and
we are finally in a position to describe this elegant algorithm
for solving the Cauchy problem for KdV.
\medskip

    \centerline {\bf The Inverse Scattering Method}
{\narrower
\noindent
To solve the KdV initial value problem
$u_t -6uu_x +u_{xxx} = 0$
with given initial potential $u(x,0)$ in $\S(\reals)$:
\roster
\item "{1)}" Apply the ``Direct Scattering Transform'', i.e., find
the discrete  eigenvalues $-\kappa_1^2,\ldots,-\kappa_N^2$ for the
Schr\"odinger  operator with potential $u(x,0)$ and compute the
Scattering Data,  i.e., the normalizing constants $c_n(0)$ and the
reflection  coefficients $b(k,0)$.

\item"{2)}" Define $c_n(t) = c_n(0)e^{4\kappa_n^3 t}$ and
$b(k,t) = b(k,0)e^{8ik^3 t}$.

\item "{3)}" Use the Inverse Scattering Transform (described
below) to compute $u(t)$ from $c_n(t)$ and $b(k,t)$.
\endroster
\par}

\subsection {The Inverse Scattering Transform}

\noindent
  Recovering the potential $u$ of a Schr\"odinger operator
$L^u$ from the Scattering Data $\S(u)$ was not something
invented for the purpose of solving the KdV initial value
problem. Rather, it was a question of basic importance to
physicists doing Cyclotron experiments, and the theory was
worked out in the mid-1950's by Kay and Moses [KM], Gelfand and 
Levitan [GL], and Marchenko [M].

  Denote the discrete eigenvalues of $u$ by
$-\kappa_1^2,\ldots,-\kappa_N^2$, the normalizing
constants by $c_1,\ldots,c_N$, and the reflection coefficients
by $b(k)$, and define a function
$$B(\xi) = \sum_{n=1}^N c_n^2 e^{-\kappa_n \xi} +
         {1\over 2\pi}\int_{-\infty}^\infty b(k) e^{ik \xi}\,dk.$$

\proclaim {Inverse Scattering Theorem}
The potential $u$ can be recovered using the formula
$u(x) = -2{d\over dx} K(x,x)$,
where $K(x,z)$ is the unique function on $\reals \times \reals$
that is zero for $z < x$ and satisfies the
Gelfand-Levitan-Marchenko Integral Equation:
$$K(x,z) + B(x+z) + \int_{-\infty}^\infty K(x,y)B(y+z) \,dy = 0.$$
\endproclaim

\noindent
(For a proof, see [DJ], Chapter 3, Section 3, or [La3], Chapter II.)
\smallskip

We will demonstrate by example how the Inverse Scattering
Method can now be applied to get explicit solutions of KdV.
But first a couple of general remarks about solving
the Gelfand-Levitan-Marchenko equation. We assume in
the following that $B$ is rapidly decreasing.

   Let $C(\reals\times\reals)$ denote the Banach space of bounded,
continuous  real-valued functions on $\reals\times\reals$ with the
$\sup$ norm. Define
${\Cal F}^B: C(\reals\times \reals) \to C(\reals\times \reals)$  by
the formula
$${\Cal F}^B(K)(x,z) =
-B(x + z) - \int_{-\infty}^\infty K(x,y)B(y+z)\,dy.$$
Then $K$ satisfies
the Gelfand-Levitan-Marchenko equation if and only if it is a fixed-point
of ${\Cal F}^B$. It is clear that ${\Cal F}^B$ is Lipschitz with
constant $\norm|B|_{L^1}$, so if $\norm|B|_{L^1} < 1$ then by the
Banach Contraction Principle the Gelfand-Levitan-Marchenko equation
has a unique solution, and it is the limit of the sequence $K_n$
defined by $K_1(x,z) = -B(x+z)$, $K_{n+1} = {\Cal F}^B(K_n)$.

  Secondly, we note that if the function $B$ is ``separable''
in the sense that it satisfies an identity of the form
$B(x+z) = \sum_{n=1}^N X_n(x) Z_n(z)$, then the
Gelfand-Levitan-Marchenko equation takes the form
$$K(x,z) + \sum_{n=1}^N X_n(x) Z_n(z) +
\sum_{n=1}^N  Z_n(z) \int_x^\infty K(x,y) X_n(y) \,dy = 0.$$
It follows that $K(x,z)$ must have the form
$K(x,z) = \sum_{n=1}^N L_n(x) Z_n(z)$. If we substitute this for
$K$ in the previous equation and define
$a_{nm}(x) = \int_x^\infty Z_m(y) X_n(y)\,dy$ then we
have reduced the problem to solving
$N$ linear equations for the unknown functions $L_n$, namely:
$L_n(x) + X_n(x) + \sum_{m=1}^N a_{nm}(x)L_m(x) = 0,$
or $X_n(x) + \sum_{m=1}^N A_{nm}(x)L_m(x) = 0,$
where $A_{nm}(x) = \delta_{nm} + a_{nm}(x)$.
Thus finally we have
$$K(x,x) = -\sum_{n=1}^N Z_n(x)
   \sum_{m=1}^N A_{nm}^{-1}(x)X_m(x).$$

\subsection {An Explicit Formula for KdV Multi-Solitons}

\noindent
  A potential $u$ is called ``reflectionless'' if all the
reflection coefficients are zero. Because of the relation
$b(k,t) = b(k,0)e^{8ik^3 t}$, it follows that if
$u(x,t)$ evolves by KdV and if it is reflectionless
at $t = 0$ then it is reflectionless for all $t$. If the
discrete eigenvalues of such a potential are
$-\kappa_1^2,\ldots,-\kappa_N^2$ and the normalizing
constants are $c_1,\ldots,c_N$, then
$B(\xi) = \sum_{n=1}^N c_n^2 e^{-\kappa_n\xi}$, so
$B(x+z) = \sum_{n=1}^N X_n(x) Z_n(z)$,
where $X_n(x) =  c_n^2 e^{-\kappa_nx}$,
and $Z_n(z) = e^{-\kappa_nz}$
and we are in the separable case just considered.
Recall that
$a_{nm}(x) = \int_x^\infty Z_m(y) X_n(y)\,dy
= c_n^2\int_x^\infty  e^{-(\kappa_n +\kappa_m)y}\,dy
=c_n^2 e^{-(\kappa_n +\kappa_m)x}/(\kappa_n +\kappa_m)$,
and that  
$$A_{nm}(x) = \delta_{nm} +a_{nm}(x) =
\delta_{nm} + c_n^2 e^{-(\kappa_n +\kappa_m)x}/(\kappa_n +\kappa_m).$$
Differentiation gives
${d\over dx} A_{nm}(x) = -c_n^2 e^{-(\kappa_n +\kappa_m)x}$,
so by a formula above
$$\eqalign{ K(x,x)
    &= -\sum_{n=1}^N Z_n(x)
         \sum_{m=1}^N A_{nm}^{-1}(x)X_m(x)\cr
    &= \sum_{n=1}^N e^{-\kappa_nx}
        \sum_{m=1}^N A_{nm}^{-1}(x)(-c_m^2 e^{-\kappa_mx})\cr
    &= \sum_{n=1}^N \sum_{m=1}^N A_{nm}^{-1}{d\over dx} A_{mn}(x)\cr
    &= \tr\left (A^{-1}(x) {d\over dx} A(x)\right)\cr
    &= {1\over \det(A(x))} {d\over dx} \det A(x)\cr
    &= {d\over dx} \log \det A(x).\cr
}$$
and so $u(x) = -2{d\over dx} K(x,x) = -2{d^2\over dx^2}\log \det A(x)$.

 If $N=1$ and we put $\kappa = \kappa_1$ it is easy to see that this
formula reduces to our earlier formula for traveling wave solutions
of the KdV equation:
$u(x,t) =-{\kappa^2\over 2} \sech^2(\kappa(x - \kappa^2 t))$.
We can also use it to find  explicit solutions $u(x,t)$ for $N=2$.
Let $g_i(x,t) = \exp(\kappa_i^3 t - \kappa_i x)$,
and set
$A = {(\kappa_1 - \kappa_2)^2 \over (\kappa_1 + \kappa_2)^2}$, then
$$u(x,t) = -2{\kappa_1^2 g_1 + \kappa_2^2 g_2 +
              2(\kappa_1 - \kappa_2)^2 g_1 g_2
              + A g_1 g_2 (\kappa_1^2 g_2 + \kappa_2^2 g_1 )
              \over
              (1+ g_1 + g_2 + A g_1 g_2)^2}.$$

  For general $N$ the solutions $u(x,t)$ that we get this way are
referred to  as the pure $N$-soliton solutions of the KdV equation.
It is not hard to show by an asymptotic analysis that for 
large negative and positive times they behave as a 
superposition of the above traveling wave solutions,
and that after the larger, faster moving waves have all
passed through the slower moving shorter ones and they have
become well-separated, the only trace of their interactions
are certain predictable ``phase-shifts'', i.e., certain constant
translations of the locations of their maxima from where they would
have been had they not interacted. (For details see [L], p.123.)

\subsection {The KdV Hierarchy}

\noindent
  By oversimplifying a bit, one can give a succinct statement
of what makes the KdV equation, $u_t - 6uu_x+u_{xxx}$,
more than just a run-of-the-mill evolution equation; namely it
is equivalent to a Lax equation, $L_t^u = [B,L^u]$, expressing that
the corresponding Schr\"odinger operator
$L^u = -{d^2\over dx^2} + u$
is evolving by unitary equivalence---so that the spectral data
for $L^u$ provides many constants of the motion for KdV,
and in fact enough commuting constants of the motion to make
KdV completely integrable.

   It is natural to ask whether KdV is unique in that respect,
and the answer is a resounding ``No!''.  In his paper introducing
the Lax Equation formulation of KdV, [La1], Peter Lax already
pointed out an important generalization. Recall that
$B = -4\D^3 + 3(u\D +\D u)$. Lax suggested that for each
integer $j$ one should look for an operator of the form
$B_j = \alpha \D^{2j+1} + \sum_{i=1}^j (b_i\D^{2i-1} + \D^{2i-1}b_i)$,
where the operators $b_i$ are to be chosen so as to make the commutator
$[B_j,L^u]$ a zero order operator---that is $[B_j,L^u]$ should be
multiplication by some polynomial, $K_j(u)$, in $u$ and its derivatives.
This requirement imposes $j$ conditions on the $j$ coefficients $b_i$,
and these conditions uniquely determine the $b_i$ as multiplications
by certain polynomials in $u$ and its derivatives.
For example, $B_0 = \D$, and the corresponding Lax Equation
$u_t = K_0(u)$ is $u_t = u_x$, the so-called
Linear Advection Equation. And of course
$B_1$ is just our friend $-4\D^3 + 3(u\D +\D u)$,
whose corresponding Lax Equation is KdV.

  $K_j(u)$ is a polynomial in the derivatives of $u$ up through
order $2j+1$, and the evolution equation $u_t = K_j(u)$ is
referred to as the $j$-th higher order KdV equation. This whole
sequence of flows is known  as ``The KdV Hierarchy'', and the
initial value problem for each of these equations can be solved
using the Inverse Scattering Method in a straightforward
generalization from the KdV case. But even more remarkably:

\proclaim {Theorem} Each of the higher order KdV equations defines
a Hamiltonian flow on $P$. That is, for each positive integer $j$
there is a Hamiltonian function $F_j: P \to \reals$ 
\rom{(}defined by a
polynomial differential operator of order $j$,
$\tilde F(u_{(0)}, \ldots, u_{(j)})$\rom{)} such that
$K_j(u) = (\sgrad F_j)_u$. Moreover, all the functions
$F_j$ are in involution, so that all the higher order
KdV flows commute with each other.\endproclaim

\noindent
The proof can be found in [La3], Chapter I.

  It should be pointed out here that the discovery of the constants
of the motion $F_k$ goes back to the earliest work on KdV as an
integrable system. In fact, it came out of the research in 1966 by
Gardner, Greene, Kruskal, and Miura leading up to their paper
[GGKM] in which the Inverse Scattering Method was introduced.
However, the symplectic structure for the phase space of KdV,
and the fact that these functions were in involution was only
discovered considerably later, in 1971 [G],[ZF]\null.	

   To the best of my knowledge, the higher order KdV equations
are not of independent interest. Nevertheless, the above theorem
suggests a subtle but important change in viewpoint towards the KdV
equation---one that proved important in further generalizing the
Inverse Scattering Method to cover other evolution equations which
are of interest for their own sake. Namely, the key player in the
Inverse Scattering Method should not be seen as the KdV equation
itself, but rather the Schr\"odinger operator $L^u$. If we want
to generalize the Inverse Scattering Method, we should first find
other operators $L$ with a ``good scattering theory'' and then look
among the Lax Equations $L_t = [M,L]$ to find interesting  candidates
for integrable systems that can be solved using  scattering methods.

  In fact, this approach has proved important in investigating
both finite and infinite dimensional Hamiltonian systems, and in
the remainder of this article we will investigate in detail one
such scheme that has not only been arguably the most sucessful
in identifying and solving important evolution equations, but
has moreover a particularly elegant and powerful mathematical
framework that underlies it. This scheme was first introduced by
Zakharov and Shabat [ZS] to study an important special equation (the
so-called Nonlinear Schr\"odinger Equation, or NLS). Soon thereafter,
Ablowitz, Kaup, Newell, and Segur [AKNS] showed that one relatively
minor modification of the Zakharov and Shabat approach recovers the
theory of the KdV equation, while another leads to an Inverse Scattering
Theory analysis for a third very important evolution equation, the
Sine-Gordon Equation (SGE).  AKNS went on to develop the Zakharov and 
Shabat technique into a general method for PDE with values in 
$2 \times 2$-matrix groups, and ZS further generalized it to the
case of $n \times n$-matrix groups. Following current custom, 
we will refer to this method as the ZS-AKNS Scheme.

\section {The ZS-AKNS Scheme}

\subsection {Flat Connections and the Lax Equation, ZCC}

\noindent
   To prepare for the introduction of the ZS-AKNS Scheme, we
must first develop some of the infra-structure on which it
is based. This leads quickly to the central Lax Equation of the
theory, the so-called ``Zero-Curvature Condition'', (or ZCC).

  First we fix a matrix Lie Group $\G$ and denote its Lie algebra
by $\g$. That is, $\G$ is some closed subgroup of the group
$\GLn(n,\Cx)$ of all $n \times n$ complex matrices, and $\g$ is
the set of all $n \times n$ complex matrices, $X$, such that
$\exp(X)$ is in $\G$. If you feel more comfortable working with
a concrete example, think of $\G$  as the group $\SLn(n,\Cx)$ of all
$n \times n$ complex matrices of determinant $1$, and $\g$ as its Lie
algebra $\sln(n,\Cx)$ of all $n \times n$ complex matrices of trace zero.
In fact, for the  original ZS-AKNS Scheme, $\G = \SLn(2,\Cx)$ and
$\g =\sln(2,\Cx)$, and we will carry out most of the later discussion
with these choices,  but for what we will do next the precise nature
of $\G$ is irrelevant.

  Let $\grad$ be a flat connection for the trivial
principal bundle $\reals^2 \times \G$. Then we can write
$\grad = d - \omega$, where $\omega$ is a $1$-form on $\reals^2$
with values in the Lie algebra $\g$. Using coordinates $(x,t)$
for $\reals^2$ we can then write $\omega = A\,dx + B\, dt$
where $A$ and $B$ are smooth maps of $\reals^2$ into $\g$.

   If $X$ is a vector field on $\reals^2$, then the covariant
derivative operator in the direction $X$ is
$\grad_X = \D_X - \omega(X)$, and in particular, the
covariant derivatives in the coordinate directions
$\D\over \D x$ and $\D\over \D t$ are
$\grad_{{\D\over \D x}} = {\D\over \D x} - A$ and
$\grad_{{\D\over \D t}} = {\D\over \D t} - B$.

   Since we are assuming that $\grad$ is flat, it determines
a global parallelism.  If $(x_0,t_0)$ is any point of $\reals^2$
 then we have a map $\psi : \reals^2 \to \G$, where
$\psi(x,t)$ is the parallel translation operator from $(x_0,t_0)$
to $(x,t)$. Considered as a section of our trivial principal bundle,
$\psi$ is covariant constant, i.e., $\grad_X \psi = 0$ for
any tangent vector field $X$. In particular, taking $X$ to be
$\D\over \D x$ and $\D\over \D t$ give the relations
$\psi_x = A\psi$ and $\psi_t = B \psi$.

  There are many equivalent ways to express the flatness of the
connection $\grad$. On the one hand the curvature $2$-form
$d\omega - \omega\wedge\omega$ is zero. Equivalently, the
covariant derivative operators in the $\D\over \D x$ and
$\D\over \D t$ directions commute, i.e.,
$[{\D\over \D x} - A,{\D\over \D t} - B] = 0$, or finally,
equating the cross-derivatives of $\psi$,
$(A\psi)_t =\psi_{xt} = \psi_{tx} = (B\psi)_x$.
Expanding the latter gives
$A_t\psi + A\psi_t = B_x \psi + B \psi_x$ or
$A_t \psi + AB\psi = B_x \psi + B A\psi$, and right multiplying
by $\psi^{-1}$ we arrive at the so-called
``Zero-Curvature Condition'': $A_t - B_x - [A,B] = 0$.
Rewriting this as $-A_t = -B_x + [B, -A]$, and noting that
$[B,{\D\over \D x}] = -B_x$, we see that the Zero-Curvature
Condition has an equivalent formulation as a Lax Equation:
$$\left({\D\over \D x} - A\right)_{\!t} =
    \left[B, {\D\over \D x} - A\right],\leqno{(ZCC)}$$
and it is ZCC that plays the central r\^ole in the
ZS-AKNS Scheme.

  Recall what ZCC is telling us. If we look at $t$ as a parameter,
then the operator ${\D\over \D x} - A(x,t_0)$ is the covariant
derivative in the $x$-direction along the line $t = t_0$, and
the Lax Equation ZCC says that as a function of $t_0$ these
operators are all conjugate. Moreover the operator $\psi(t_0,t_1)$
implementing the conjugation between the time $t_0$ and the time $t_1$
satisfies $\psi_t = B\psi$, which means it is parallel translation
from $(x,t_0)$ to $(x,t_1)$ computed by going ``vertically'' along the
curve $t \mapsto (x,t)$. But since ${\D\over \D x} - A(x,t_0)$ generates
parallel translation along the horizontal curve $x \mapsto (x,t_0)$,
what this amounts to is the statement
that parallel translating horizontally from
$(x_0,t_0)$ to $(x_1,t_0)$ is the same as parallel translation
vertically from $(x_0,t_0)$ to $(x_0,t_1)$ followed by
parallel translation horizontally from $(x_0,t_1)$ to $(x_1,t_1)$
followed by parallel translation vertically from $(x_1,t_1)$ to
$(x_1,t_0)$. Thus, in the case of ZCC, the standard interpretation
of the meaning of a Lax Equation reduces to a special
case of the theorem that if a connection has zero curvature
then the holonomy around a contractible path is trivial.

\subsection {Some ZS-AKNS Examples}

\noindent
  The ZS-AKNS Scheme, is a method for solving the initial value
problem for certain (hierarchies of) evolution equations on a
space of ``potentials'' $P$. In general $P$ will be of the
form $\S(\reals,V)$, where $V$ is some finite dimensional real
or complex vector space, i.e., each potential $u$ will be a map
$x \mapsto u(x)$ of  Schwartz class from $\reals$ into $V$. (A
function $u$ with  values in $V$ is of Schwartz class if, for
each linear functional $\ell$ on $V$, the scalar valued function
$\ell\circ u$ is of Schwartz class, or equivalently if, when we
write $u$ in terms of a fixed basis for $V$, its components are
of Schwartz class.)  The evolution equations in question are of
the form $u_t = F(u)$ where the map $F: P \to P$ is a
``polynomial differential operator''---i.e., it has the form
$F(u) = p(u,u_x,u_{xx},\ldots)$, where $p$ is a polynomial mapping
of $V$ to itself.

  When we say we want to solve the initial value
(or ``Cauchy'') problem for such an equation, we of course
mean that given $u^0 =u(x,0)$ in $P$ we want to find a smooth
map $t \mapsto u(t) = u(x,t)$ of $\reals$ to $P$ with
$u(0) = u^0$ and
$u_t(x,t) = p(u(x,t),u_x(x,t),u_{xx}(x,t),\ldots)$. In essence,
we want to think of $F$ as a vector field on $P$ and construct
the flow $\phi_t$ that it generates. (Of course, if $P$ were
a finite dimensional manifold, then we could construct the
flow $\phi_t$ by solving a system of ODE's, and as we shall
see, the ZS-AKNS Scheme allows us in certain cases to solve the
PDE $u_t = p(u,u_x,u_{xx},\ldots)$ by reducing it to ODE's.)

  The first and crucial step in using the ZS-AKNS Scheme to study
a particular such evolution equation consists in setting up an
interpretation of $A$ and $B$ so that the equation
$u_t = p(u,u_x,u_{xx},\ldots)$ becomes a special  case of ZCC.

   To accomplish this, we first identify $V$ with a subspace of $\g$
(so that $P = \S(\reals,V)$ becomes a subspace of $\S(\reals,\g)$),
and define a map $u \mapsto A(u)$ of $P$ into $C^\infty(\reals,\g)$
of the form $A(u) =$ const${}\,+ \,u$, so that if $u$ depends
parametrically on $t$ then $({\D \over \D x} -A(u))_t = -u_t$.

  Finally (and this is the difficult part) we must define a map
$u \mapsto B(u)$ of $P$ into $C^\infty(\reals,\g)$ so that
$[B(u), {\D \over \D x} -A(u)] = -p(u,u_x,u_{xx},\ldots)$.

  To interpret the latter equation correctly, and in particular to
make sense out of the commutator bracket in a manner consistent
with our earlier interpretation of $A$ and $B$, it is important to
be clear about the interpretation $A(u)$ and $B(u)$ as operators, and
in particular to be precise about the space on which they are operating.
This is just the space $C^\infty(\reals,\gln(2,\Cx))$ of smooth
maps $\psi$ of $\reals$ into the space of all complex $2 \times 2$
matrices. Namely, we identify $A(u)$ with the zero-order differential
operator mapping $\psi$ to $A(u)\psi$, the pointwise matrix product
of $A(u)(x)$ and $\psi(x)$, and similarly with $B(u)$.
(This is a complete analogy with the KdV situation, where in
interpreting the Schr\"odinger operator, we identified our
potential $u$ with the operator of multiplication by $u$.)
Of course $({\D \over \D x}\psi)(x) = \psi_x$.

  We will now illustrate this with three examples: the
KdV equation, the Nonlinear Schr\"odinger Equation (NLS),
and the Sine-Gordon Equation (SGE). In each case $V$ will
be a one-dimensional space that is embedded in
the space of off-diagonal complex matrices 
$\mat(0;b;c;0)$,
and in each case $A(u) = \a\lambda + u$, where $\lambda$
is a complex parameter, and $\a$ is the constant, diagonal,
trace zero matrix $\a = \mat(-i;0;0;i)$.

\medskip
\noindent
{\bf Example 1.} [AKNS] Take $u(x) = \mat(0;q(x);-1;0)$,
and let
$$B(u) = \a\lambda^3 + u\lambda^2 +
   \mat({i\over 2}q;{i\over 2}q_x;0;-{i\over 2}q)\lambda +
   \mat({q_x\over 4};{-q^2\over 2};{q\over 2};-{q_x\over 4}).$$
Then an easy computation shows that ZCC is satisfied if and
only if $q$ satisfies KdV in the form
$q_t = -{1\over 4}(6qq_x + q_{xxx})$.
\goodbreak
\medskip
\noindent
{\bf Example 2.} [ZS] Take $u(x) = \mat(0;q(x);-\bar q(x);0)$,
and let
$$B(u) = \a\lambda^2 + u\lambda +
\mat({i\over 2}|q|^2;{i\over 2}q_x;-{i\over 2}\bar q_x;-{i\over 2}|q|^2).$$
In this case ZCC is satisfied if and only if $q(x,t)$ satisfies
the so-called Nonlinear Schr\"odinger Equation (NLS)
$q_t= {i\over 2}(q_{xx} + 2|q|^2q)$.
\medskip
\noindent
{\bf Example 3.} [AKNS] Take
$u = \mat(0;-{q_x(x)\over 2};{q_x(x)\over 2};0)$, and let
$B(u) = {1\over\lambda}v$ where
$v(x) = {i\over 4}\mat(\cos q(x);\sin q(x);\sin q(x);-\cos q(x))$.
In this case, ZCC is satisfied if and only if
$q$ satisfies the Sine-Gordon Equation (SGE) in the form
$q_{xt} = \sin q$.

  In the following description of the ZS-AKNS Scheme, we will
state definitions and describe constructions in a way that works
for the general ZS-AKNS case---and we will even make occasional
remarks explaining what modifications are necessary
to extend the theory to the more general case of $n \times n$
matrix groups. (For the full details of this latter
generalization the reader should consult [Sa]\null.)
However, working out details in even full
ZS-AKNS generality would involve many
distracting  detours, to discuss various special situations
that are irrelevant to the main ideas. So, for ease and clarity
of exposition, we will carry out most of the further discussion
of the ZS-AKNS Scheme within the framework of the NLS Hierarchy.

\subsection {The Uses of Solitons}

   There are by now dozens of ``soliton equations'', but not only were the three 
examples from the preceding section the first to be discovered, they are 
also the best known, and in many ways still the most interesting and important. 
In fact, in addition to their simplicity and their Hamiltonian nature, 
each has certain special properties that give them a ``universal'' character, 
so that they are almost sure to arise as approximate models in any physical 
situation that exhibits these properties. In this section I will try to say 
a little about these special features, and also explain how these equations 
have been used in both theoretical and applied mathematics.

  We have already discussed in some detail the historical background 
and many of the interesting features and applications of the KdV equation, so 
here I will only re-iterate the basic property responsible for its frequent
appearance in applied problems. In the KdV equation there is an extraordinary 
balance between the shock-forming tendency of its non-linear term $uu_x$ and 
the dispersive tendency of its linear term $u_{xxx}$, and this balance is responsible 
for the existence of remarkably stable configurations (solitons) that scatter elastically 
off one another under the KdV evolution. Moreover KdV is the simplest non-dissipative  
wave-equation equation with these properties.

 The Sine-Gordon equation is even older than the KdV equation; it arose first
in the mid-nineteenth century as the master equation for ``pseudo-spherical surfaces''
(i.e., surfaces of constant negative Gaussian curvature immersed in $\reals^3$).
Without going into the details (cf. [Da] and [PT], Part I, Section 3.2), the
Gauss-Codazzi  equations for such surfaces reduce to the Sine-Gordon equation, so that
by the  ``Fundamental Theorem of Surface Theory'', there is a bijective correspondence 
between isometry classes of isometric immersions of the hyperbolic plane into 
$\reals^3$ and solutions to the Sine-Gordon equation. Because of this (and the great 
interest in non-Euclidean geometry during the latter half of the last century) a 
prodigious amount of effort was devoted to the study of the Sine-Gordon equation by 
the great geometers of that period, resulting in a beautiful body of results, most 
of which can be found in G.~Darboux' superb treatise on surface theory {\it Le\c cons 
sur la Th\'eorie G\'en\'erale des Surfaces\/} [Da]. 

   One of the most notable features of this theory is the concept of a 
``B\"acklund Transformation''. Starting from any solution of the Sine-Gordon 
equation, this creates a two-parameter family of new solutions. One slight 
complication is that the construction of the new solutions requires solving a
certain ordinary differential equation. However the so-called ``Bianchi 
Permutability Formula'' allows us to easily compose B\"acklund Transformations.
That is, once we have found this first set of new solutions, we can apply 
another B\"acklund Transformations to any one of them to get still more solutions
of Sine-Gordon, and this second family of  new solutions can be written down explicitly
as algebraic functions of the first set, without solving any more ODEs. Moreover, we can 
continue inductively in this manner, getting an infinite sequence of families of more 
and more complex solutions to the Sine-Gordon equations (and related pseudospherical
surfaces). If we take as our starting solution the identically zero (or ``vacuum'') 
solution to the Sine-Gordon equation, this process can be carried out explicitly. At 
the first stage we get the so-called Kink (or one-soliton) solutions to the Sine-Gordon
equation, and the corresponding family of pseudospherical surfaces is the Dini 
family (including the well-known pseudosphere). Using the Bianchi Formula  
once gives rise to the two-soliton solutions of Sine-Gordon and the corresponding 
K\"uen Surface, and  repeated application leads in principle to all the higher 
soliton solutions of the Sine-Gordon equations (cf. [Da], [PT], {\it loc. cit.\/} for 
more details). In fact, the classical geometers knew so much about the ``soliton 
sector'' of solutions to Sine-Gordon that it might seem surprising at first that 
they did not go on to discover``soliton mathematics'' a century before it actually 
was. But of course they knew only half the story---they knew nothing of the dispersive,
non-soliton solutions to Sine-Gordon and had no imaginable way to discover the Inverse
Scattering Transform, which is the key to a full understanding of the space of all
solutions. (And finally, they probably never looked at Sine-Gordon as an evolution
equation for a one-dimensional wave, so they didn't notice the strange scattering
behavior of the solutions that they had calculated.)

    Nevertheless, their work did not go in vain. As soon as it was realized that 
Sine-Gordon was a soliton equation, it was natural to ask whether KdV also 
had an analogous theory of B\"acklund transformations that, starting from 
the vacuum solution marched up the soliton ladder. It was quickly discovered 
that this was in fact so, and while B\"acklund transformations have remained
until recently one of the more mysterious parts of soliton theory, each newly
discovered soliton equation was found to have an associated theory of B\"acklund 
transformations. Indeed this soon came to be considered a hallmark of the ``soliton 
syndrome'', and a test that one could apply to detect soliton behavior. A natural 
explanation of this relationship follows from the Terng-Uhlenbeck Loop Group 
approach to soliton theory, and we will remark on it briefly at the end of this 
article. For full details see [TU2].

  The Sine-Gordon equation has also been proposed as a simplified model 
for a unified field theory, and derived as the equation governing the propogation 
of dislocations in a crystal lattice, the propogation of magnetic flux in a Josephson 
junction transmission line, and many other physical problems.

  The Nonlinear Schr\"odinger Equation has an interesting pre-history. It was discovered 
``in disguise'' (and then re-discovered at least three times, cf. [Ri]) in the early part
of this century.  In 1906, Da Rios wrote a master's thesis [DaR] under the direction 
of Levi-Civita, in which he modeled the free evolution of a thin vortex-filament in a 
viscous liquid by a time-dependent curve $\gamma(x,t)$ in $\reals^3$ satisfying the 
equation $\gamma_t = \gamma_x \times \gamma_{xx}$. Now by the Frenet equations, 
$\gamma_x \times \gamma_{xx} = \kappa B$ where $\kappa = \kappa(x,t)$ is the curvature 
and $B$ the binormal, so the filament evolves by moving in the direction of its 
binormal with a speed equal  to its curvature. This is now often called the 
``vortex-filament equation'' or the  ``smoke-ring equation''. In 1971, Hasimoto 
noticed a remarkable gauge transformation that transforms the vortex-filament 
equation to the Nonlinear Schr\"odinger equation. In fact, if $\tau(\cdot,t)$ denotes 
the torsion of the curve $\gamma(\cdot,t)$,  then the  complex quantity 
$q(x,t) = \kappa(x,t) \exp(i\int\tau(\xi,t)\, d\xi)$  satisfies NLS if and only 
if $\gamma$ satisfies the vortex-filament equation. 

  But it is as an ``envelope equation'' that NLS has recently come into its own. 
If a one-dimensional, amplitude modulated, high-frequency wave is moving in a
highly dispersive and non-linear medium, then to a good approximation the
evolution of the wave envelope (i.e., the modulating signal) in a coordinate
system moving at the group velocity of the wave will satisfy NLS. Without going
into detail about what these hypotheses mean (cf. [HaK]) they do in fact apply 
to the light pulses travelling along optical fibers that are rapidly becoming the 
preferred means of communicating information at high bit-rates over long distances.
Solitons solutions of NLS seem destined play a very important r\^ole in keeping 
the Internet and the World Wide Web from being ruined by success. The story is only
half-told at present, but the conclusion is becoming clear and it is too good a story
to omit.

  For over a hundred years, analogue signals travelling over copper wires provided
the main medium for point-to-point communication between humans. Early implementations
of this medium (twisted pair) were limited in bandwidth (bits per second) to about
100 Kb/s per channel. By going over to digital signalling instead of analogue, one can 
get up to the 1 Mb/s range, and using coaxial cable one can squeeze out another several 
orders of magnitude. Until recently this seemed sufficient. A bandwidth of about 1 Gb/s 
is enough to satisfy the needs of the POTS (plain old telephone system) network
that handles voice communication for the entire United States, and that could be handled
with coaxial cable and primitive fiber optic technology for the trunk lines between central
exchanges, and twisted pairs for the low bandwidth ``last mile'' from the exchange to a 
user's home. And as we all know, a coaxial cable has enough bandwidth to provide us with 
several hundred channels of television coming into our homes.

  But suddenly all this has changed. As more and more users are demanding very high data-rate
services from the global Internet, the capacities of the communication providers have been
stretched to and beyond their limits, and they have been desperately trying to keep up. The
problem is particularly critical in the transoceanic links joining the US to Asia and Europe.
Fortunately, a lot of fiber optic cables have been laid down in the past decade, and even
more fortunately these cables are being operated at bandwidths that are very far below their
theoretical limits of about 100 GB/s. To understand the problems involved in using these
resources more efficiently, it is necessary to understand how a bit is transmitted along
an optical fiber. In principle it is very simple. In so-called RZ (return-to-zero) coding,
a pulse of high-frequency laser-light is sent to indicate a one, or not sent to indicate a 
zero. The inverse of the pulse-width in seconds determines the maximum bandwidth of the channel. 
A practical lower bound for the pulse-width is about a pico-second ($10^{-12}$ seconds) giving 
an upper bound of about 1000 GB/s for the bandwidth. But of course there are further practical
difficulties that limit data-rates to well below that figure (e.g., the pulses should be 
well-separated, and redundancy must be added for error correction) but actual data transmission
rates over optical fibers in the 100 GB/s range seems to be a reasonable goal 
(using wavelength-division-multiplexing).

  But there are serious technical problems. Over-simplifying somewhat, a major obstacle to
attaining such rates is the tendency of these very short pico-second pulses to disperse as
they travel down the optical fiber. For example, if an approximate square-wave pulse is sent,
then dispersion will cause very high error rates after only several hundreds of miles. 
However if the pulses are carefully shaped to that of an appropriate NLS soliton, then the
built-in stability of the soliton against dispersion will preserve the pulse shape over
very long distances, and theoretical studies show that error-free propogation at 10 GB/s 
across the Pacific is feasible with current technology, even without multi-plexing.
(For further details and references see [LA].)

\subsection {Nonlinear Schr\"odinger as a Hamiltonian Flow}

\noindent
  Let $\G$ denote the group $\SUn(2)$ of unitary $2 \times 2$ complex
matrices of determinant $1$, and and $\g$ its Lie algebra, $\sun(2)$,
of skew-adjoint complex matrices of trace $0$. The $3$-dimensional
real vector space $\g$ has a natural positive definite inner product
(the Killing form), defined by $\dip<<a;b>> = -{1\over 2} \tr(ab)$.
It is characterized (up to a constant factor) by the fact that it is
``Ad-invariant'', i.e., if $g \in \G$ then
$\dip<<\Ad(g)a;\Ad(g)b>> = \dip<<a;b>>$,  where
$\Ad(g) : \g \to \g$ is defined by $\Ad(g) a = gag^{-1}$.
Equivalently, for each element $c$ of $\g$,
$\ad(c) : \g \to \g$ defined by $\ad(c) a = [c,a]$
is skew-adjoint with respect to the Killing form:
$\dip<<[c,a];b>> + \dip<<a;[c,b]>> = 0$.

  We denote by $\T$ the standard maximal torus of $\G$, i.e.,
the group $\diag(e^{-i\theta},e^{i\theta})$ of diagonal, unitary
matrices of determinant $1$, 
and \t will denote its Lie algebra
$\diag(-i\theta,i\theta)$ of skew-adjoint, diagonal matrices
of trace zero. We define the specific element $\a$ of \t by
$\a = \diag(-i,i)$.

  The orthogonal complement, $\tperp$, of \t in $\g$ will play an
important r\^ole in what follows. It is clear that $\tperp$ is just
the space of ``off-diagonal''  skew-adjoint matrices, i.e.,
those with all zeros on the diagonal. (This follows easily from the
fact that the product of a diagonal matrix and a ``off-diagonal''
matrix is again off-diagonal, and so of trace zero.) Thus $\tperp$ is
the space of matrices of the form $\mat(0;q;-\bar q;0)$ where
$q \in \Cx$, and this gives a natural complex structure to the
$2$-dimensional real vector space $\tperp$.

  Note that \t is  just the kernel (or zero eigenspace) of $\ad(\a)$.
Since $\ad(\a)$  is skew-adjoint with respect to the Killing form, it
follows that  $\ad(\a)$ leaves $\tperp$ invariant, and we will denote
$\ad(\a)$ restricted to $\tperp$ by $J : \tperp \to \tperp$. A trivial
calculation shows that
$J \mat(0;q;-\bar q;0) = \mat(0;2iq;-\overline {2iq};0)$.

\noindent
{\bf Remark.} In the generalization to $\SUn(n)$, we choose $\a$
to be a diagonal element of $\sun(n)$ that is ``regular'', i.e., has
distinct eigenvalues. Then the Lie algebra \t of the maximal
torus $\T$ (the diagonal subgroup of $\SUn(n)$) is still all diagonal
skew-adjoint operators of trace zero, and is again the null-space of
$\ad(\a)$. Its orthogonal complement, $\tperp$, in $\sun(n)$ is thus still
invariant under $\ad(\a)$, but now it is no longer a single complex
eigenspace, but rather the direct sum of complex $\ad(\a)$ eigenspaces
(the so-called ``root spaces'').
\smallskip

  We define the phase space $P$ for the NLS Hierarchy by
$P = \S(\reals,\tperp)$, i.e., $P$ consists of all ``potentials''
$u$ that are Schwartz class maps of $\reals$ into $\tperp$:
$x \mapsto u(x) = \mat(0;q(x);-\bar q(x);0)$. Clearly $u \mapsto q$
establishes a canonical identification of $P$ with the space
$\S(\reals,\Cx)$ of all complex-valued Schwartz class functions
on the line. We define an $L^2$ inner product on $P$, making it
into a real pre-hilbert space, by
$\IP<u_1,u_2> = \int_{-\infty}^\infty\dip<<u_1(x);u_2(x)>>\,\,dx
          = -{1\over 2}\int_{-\infty}^\infty\tr(u_1(x) u_2(x))\,dx$.
When this is written in terms of $q$ we find,
$\IP<u_1,u_2>= \Re(\int_{-\infty}^\infty
q_1(x)\overline {q_2(x)}\,dx$).  And finally, if we decompose
$q_1$ and $q_2$ into their real and imaginary parts:
$q_j = v_j + i w_j$, then
$\IP<u_1,u_2>= \int_{-\infty}^\infty
(v_1 v_2 +w_1 w_2)\,dx$.

  We ``extend'' $J: \tperp \to \tperp$ to act pointwise on $P$,
i.e., $(Ju)(x)= J(u(x))$, and since $J$ is skew-adjoint, we can
define a skew bilinear form $\Omega$ on $P$ by
$$\eqalign{\Omega(u_1,u_2)
   &= \IP<J^{-1}u_1,u_2>\cr
   &= \Re\left(\int_{-\infty}^\infty {1\over 2i}q_1
          \overline {q_2}\,dx\right)\cr
   &= -{1\over 2} \Re \left(i\int_{-\infty}^\infty q_1
          \overline {q_2}\,dx\right)\cr
   &=  {1\over 2}\Im \left(\int_{-\infty}^\infty q_1
          \overline {q_2}\,dx\right).\cr
}$$

  Considered as a differential $2$-form on the real topological
vector space $P$, $\Omega$ is constant and hence closed. On the
other hand, since $J:P \to P$ is injective, it follows that $\Omega$
is weakly  non-degenerate, and hence a symplectic structure on $P$.

   From the definition of $\Omega$ we have
$\Omega(Ju_1,u_2) = \IP<u_1,u_2>$, thus if $F : P \to P$ has
a Riemannian gradient $\grad F$ then
$\Omega(J(\grad F)_{u_1},u_2) = \IP<(\grad F)_{u_1} ,u_2> =dF_{u_1}(u_2)$,
and so $\sgrad F = J\grad F$. In particular, if $F_1$ and $F_2$
are any two Hamiltonian functions on $P$ then their Poisson bracket
is given by the formula
$\PB(F_1,F_2) = \Omega(\sgrad F_2, \sgrad F_1)
              = \Omega(J\grad F_2, \grad F_1)
              = \IP<\grad F_2, J\grad F_1>
              = \IP<J\grad F_1, \grad F_2>$.

  A Calculus of Variations functional on $P$, $F:P\to\reals$, will be of the form
$F(u) = \int_{-\infty}^\infty\tilde F(v,w,v_x,w_x,\ldots)\,dx$, where
$q = v + i w$, and the differential of $F$ is given by
$dF_u(\delta u) = \int_{-\infty}^\infty \left(
{\delta F\over \delta v} \delta v +
{\delta F\over \delta w} \delta w\right) \,dx$, or equivalently
$$ dF_u(\delta u) ={1\over 2}\Re\left(\int_{-\infty}^\infty \left(
{\delta F\over \delta v}  +
i{\delta F\over \delta w}\right) (\delta v - i \delta w)\,dx\right) ,$$
where as usual
${\delta F\over \delta v} =
{\partial \tilde F\over \partial v} - {\partial\over \partial x}
\left({\partial \tilde F\over \partial v_x}\right) +
{\partial^2\over \partial x^2}
\left({\partial \tilde F\over \partial v_{xx}}\right)
-\ldots$, and a similar expression for ${\delta F\over \delta w}$.
However, it will be more convenient to give the polynomial
differential operator $\tilde F$ as a function of
$q = u+iv,\bar q = u - iv, q_x = u_x + i v_x, \bar q_x = u_x - i v_x,\ldots$
instead of as a function of $u,v$ and their derivatives.
Since $v = {1\over 2} (q + \bar q)$ and $w = {1\over 2i} (q - \bar q)$,
by the chain-rule,
${\partial\tilde F\over\partial \bar q}
= {1\over 2} \left({\partial\tilde F\over\partial v} +
i {\partial\tilde F\over\partial w}\right)$, with similar formulas for
${\partial\tilde F\over\partial \bar q_{x}}$,
${\partial\tilde F\over\partial \bar q_{xx}}$, etc.
Thus if we define
${\delta F\over \delta \bar q} =
{\partial \tilde F\over \partial \bar q} - {\partial\over \partial x}
\left({\partial \tilde F\over \partial \bar q_x}\right) +
{\partial^2\over \partial x^2}
\left({\partial \tilde F\over \partial \bar q_{xx}}\right) -\ldots$, then
${\delta F\over \delta \bar q} =
{\delta F\over \delta v} + i {\delta F\over \delta w}$,
and it follows that
$ dF_u(\delta u) ={1\over 2}\Re\left(\int_{-\infty}^\infty
{\delta F\over \delta \bar q} \overline{\delta q}\,dx\right) $,
where $\delta q = \delta v + i \delta w$, so
$\delta u = \mat(0;\delta q;-\overline {\delta q};0)$.
Recalling the formulae for $\IP<u_1,u_2>$, it follows that
$\grad F_u = \mat(0;{\delta F\over \delta \bar q};
-\overline {\delta F\over \delta \bar q};0)$,
and so
$\sgrad F_u = \mat(0;2i{\delta F\over \delta \bar q};
-\overline {2i{\delta F\over \delta \bar q}};0)$.
Thus, expressed in terms of $q$, the Hamiltonian flow in $P$
defined by $F$ is $q_t = 2i{\delta F\over \delta \bar q}$.

  If we take
$F(u) = -{1\over 2}\tr(u^4 + u_x^2) = {1\over 2}(|q|^4 + |q_x|^2)$,
then $\tilde F(q,\bar q, q_x, \bar q_x)
  = {1\over 2}(q^2\bar q^2 + q_x \bar q_x)$ and
${\delta F\over \delta \bar q} =
     q^2\bar q + {1\over 2}{\D\over \D x}(q_x)=
     ({1\over 2} q_{xx} + |q^2|q)$, and the
Hamiltonian equation is
$q_t = i(q_{xx} + 2|q^2|q)$, which is NLS.

\subsection {The Nonlinear Schr\"odinger Hierarchy}

\noindent
  For each potential $u$ in $P$ and complex number $\lambda$
we define an element $A(u,\lambda)$ of $C^\infty(\reals,\sln(2,\Cx))$ by
$A(u,\lambda) = \a\lambda + u = \mat(-i\lambda;q;-\bar q;i\lambda)$.
$A(u,\lambda)$ will play an important r\^ole in what follows, and you
should think of it as a  as a zero-order differential operator on
$C^\infty(\reals,\gln(n,\Cx))$, acting by pointwise multiplication on the
left. We are now going to imitate the construction of the
KdV Hierarchy. That is, we will look for a sequence of maps
$u \mapsto B_j(u,\lambda)$ of $P$ into $C^\infty(\reals,\sln(2,\Cx))$
(polynomials of degree $j$ in $\lambda$) such that the sequence
of ZCC Lax Equations $u_t = [B_j, {\D\over\D x} - A]$
is a sequence of commuting Hamiltonian flows on $P$, which
for $j=2$ is the NLS flow.

\proclaim {NLS Hierarchy Theorem}  For each $u$ in $P$ there exists
a sequence of smooth maps $Q_k(u) : \reals \to \sun(2)$ with the following
properties:
\roster
\item "{a)}" The $Q_k(u)$ can be determined recursively by:

             i\rom{)} $Q_0(u)$ is the constant matrix $\a$.

             ii\rom{)}
               $ [\a,Q_{k+1}(u)] =  (Q_k(u))_x + [Q_k(u),u]$,

             iii\rom{)}
               $(Q_k(u))_x + [Q_k(u),u]$ is off-diagonal.
\item "{b)}" If we define $B_j(u,\lambda) = \sum_{k=0}^j Q_k(u)\lambda^{k-j}$,
           and consider $B_k(u,\lambda)$ as a zero-order linear differential
           operator acting by pointwise matrix multiplication
           on elements $\psi$ of $C^\infty(\reals,\gln(2,\Cx))$,
           then the conditions ii)
           and iii) of a) are equivalent to demanding that the commutators
           $[B_j(u,\lambda), {\D\over\D x} - A(u,\lambda)]$ are independent
           of $\lambda$ and have only off-diagonal entries. In fact
           these commutators have the values:\hfil\break
           $[B_j(u,\lambda), {\D\over\D x} - A(u,\lambda)]
           =[\a,Q_{j+1}(u)] = (Q_j(u))_x -[u,Q_j(u)]$.
\item "{c)}" The matrix elements of $Q_k(u)$ can be determined so that they
           are polynomials in the derivatives (up to order $k-1$) of the
           matrix  entries of $u$, and this added requirement makes them
           uniquely determined. We can then
           regard $Q_k$ as a map of $P$ into $C^\infty(\reals,\sun(2))$.
           Similarly, for each real $\lambda$, $u \mapsto B_j(u,\lambda)$
           is a map of $P$ into $C^\infty(\reals,\sun(2))$.
\item "{d)}" If follows that the sequence of ZCC Lax Equations,
           $\left({\D\over\D x}- A\right)_t = [B_j, {\D\over\D x} - A]$,
           (or equivalently  $u_t =[\a,Q_{j+1}(u)]$) determine flows on $P$,
           the so-called higher order NLS flows. (The $j$-th
           of these is called the $j$-th NLS flow and the second is
           the usual NLS flow).
\item "{e)}" If we define Hamiltonians on $P$ by
          $H_k(u) = -{1\over k+1}\int_{-\infty}^\infty\tr(Q_{k+2}(u)\a)\,dx$,
           then $(\grad H_k)_u$ is the off-diagonal part of $Q_{k+1}(u)$.
\item "{f)}" It follows that the $j$-th NLS flow is Hamiltonian, and in fact
           is given by $u_t = (\sgrad H_k)_u$.
\item "{g)}" The Hamiltonian functions $H_k$ are in involution, i.e.,
           the Poisson brackets $\PB(H_k,H_l)$ all vanish, so that
           all the NLS flows on $P$ commute.
\endroster
\endproclaim

\noindent
{\bf Remark.} We will give part of the proof of this important theorem
here, and finish the proof later when we have developed more machinery.
However first we comment on the changes that are necessary when we go
from $2$ to $n$ dimensions, (i.e., replace $\gln(2,\Cx)$ by $\gln(n,\Cx)$,
and $\sun(2)$ by $\sun(n)$). In fact, surprising few changes are necessary.
The maximal torus $\T$  still consists of diagonal unitary matrices of trace
$1$ but now has dimension $(n-1)$ rather than $1$. We replace $\a$ by any
regular element of $\T$ (i.e., one with distinct elements on the diagonal).
This is equivalent to the key condition that $\T$ is the commutator of $\a$.
The biggest change is that to get the family of commuting Hamiltonian
flows we must now choose a second element $\b$ of $\T$, and replace
$Q_j(u) = Q_{\a,j}(u)$ by the more general $Q_{\b,j}(u)$, and the
$B_j(u,\lambda) = B_{\a,j}(u,\lambda)$ by the more general
$B_{\b,j}(u,\lambda) = \sum_{j=0}^j Q_{\b,k}(u)\lambda^{k-j}$.
The only further change is that i) of a) now reads
``$Q_{b,0}(u)$ is the constant matrix $\b$.'' {\it Mutatis mutandis\/},
everything else remains the same. For full details, see [Sa]\null.	

\smallskip\noindent
\proof
  Some easier parts of the proof will be indicated here, while other
more difficult steps will be deferred until after we discuss the
ZS-AKNS direct scattering theory, at which point they will be
much easier to demonstrate.

  The coefficient of $\lambda^{j-k}$ in the commutator
$[B_j(u,\lambda), {\D\over\D x} - A(u,\lambda)]$
is easily computed, and for $k = 0$ to $j-1$ we find
$-(Q_k(u))_k - [Q_k(u),u] - [Q_{k+1}(u),\a]$, while for $k=j$
(i.e., the term independent of $\lambda$) we get
$-(Q_j(u))_x - [Q_j(u),u]$, and c) is now immediate.

   If we write $Q_k(u)$ as the sum of its diagonal part, $T_k(u)$,
and its off-diagonal part, $P_k(u)$, then since $\ad(\a)$
annihilates diagonal matrices and is an isomorphism on the
off-diagonal matrices,  \hfil\break
\centerline{
$[\a,Q_{k+1}(u)] = \ad(\a)(T_{k+1}(u)) + \ad(\a)(P_{k+1}(u))
   = \ad(\a)(P_{k+1}(u))$, }
so by ii) of a):
$$P_{k+1}(u) = \ad(\a)^{-1}((P_k(u))_x + [T_k(u),u]).$$
(We have used the fact that, since $u$ is off-diagonal,
$[u,T_k(u)]$ is off-diagonal while $[u,P_k(u)]$ is diagonal.)

 Next note that condition iii) of statement a) can now be written as
$(T_j(u))_x = [u,P_j(u)]$
(because $[u,P_j(u)]$ is diagonal while  $[u,T_j(u)]$ is
off-diagonal). So we can write
$$T_{k+1}(u) = \int_{-\infty}^x [u,P_{k+1}(u)]\,dx,$$
where of course the indefinite integral is to be taken matrix element
by matrix element. Together, the latter two displayed equations
give an explicit recursive definition of
$Q_{k+1} = P_{k+1} + T_{k+1}$ in terms of $Q_k = P_k + T_k$.

  For example, since $Q_0(u) = \a$ we conclude that $P_0(u) = 0$ and $T_0(u) = \a$.
Then the formula for $P_{k+1}$ gives
$P_1(u) = \ad(\a)^{-1}( 0 + [\a,u]) = u$, and since $[u,u] = 0$,
the formula for $T_{k+1}$ gives
$T_1(u) = 0$, and therefore $Q_1(u) = P_1(u) = u$.

  Continuing, we find next that
$P_2(u) = \ad(\a)^{-1}(u_x) =\mat(0;-{i\over2}q_x;{i\over 2}\bar q_x;0)$,
and $(T_2(u))_x = [u,P_2(u)] =
\mat({i\over2}(q_x\bar q + q\bar q_x);0;0;{-{i\over2}(q_x\bar q + q\bar q_x)})$,
which gives by integration
$T_2(u) = \mat({i\over2}|q|^2;0;0;-{i\over2}|q|^2)$ and
$Q_2(u) = P_2(u) + T_2(u) =
\mat({i\over2}|q|^2;-{i\over2}q_x;{i\over 2}\bar q_x;-{i\over2}|q|^2)$.
(By what we have seen earlier, this shows that the second flow is
indeed the NLS flow).

  We could continue for another several steps, and at each stage,
after computing $P_j(u)$ and then $[P_j(u),u]$, the anti-derivative
of the latter turns out to be in $\S(\reals,\Cal{T})$,
so $Q_j(u) = P_j(u) + T_j(u)$ is in $\S(\reals,\sun(2))$.
(Note that this is clearly equivalent to the statement that
$\int_{-\infty}^\infty [u,P_{k+1}(u)]\,dx = 0$.)

  Unfortunately, no one has come up with a
simple inductive proof of that fact, so at this stage we are faced
with the unpleasant possibility that our recursive process might
lead to some $T_j(u)$ (and hence $Q_j(u)$) that does not vanish at
infinity. Later on, after we have discussed
the scattering theory for the ZS-AKNS Scheme, we will find a simple
argument to show that this cannot happen, and at that point we will
have a proof of statements a) through d). Similarly, I do not know
a proof of statement e) that avoids scattering theory, so I will again
defer the proof.

  Recalling that $\ad(\a)$, (i.e., bracketing with
$\a)$ annihilates diagonal matrices, it follows from e) that
$\sgrad H_k = J(\grad H_k) = [\a,Q_{k+1}]$, and so by d) the
$j$-th NLS flow is given by $u_t = (\sgrad H_k)_u$, which is f).

  For g), recall
$\PB(H_k,H_l) = \IP<J\grad H_k, \grad H_l> =
    \IP<[\a, Q_{k+1}(u)], Q_{l+1}(u)>$,
and using this formula, the $\ad$-invariance of the Killing form,
and the recursion relation
$ [\a,Q_{j+1}(u)] = (Q_j(u))_x - [u,Q_j(u)]$,
we will give an inductive argument that the $H_k$ are in involution.

\proclaim {Lemma 1}
\roster
\item "{a\rom{)}}" $\IP<[u,Q_j(u)],Q_k(u)> + \IP<Q_j(u),[u,Q_k(u)]> = 0$.
\item "{b\rom{)}}" $\IP<[u,Q_j(u)],Q_j(u)> = 0$.
\item "{c\rom{)}}" $\IP<(Q_j(u))_x,Q_k(u)> + \IP<Q_j(u),(Q_k(u))_x> = 0$.
\item "{d\rom{)}}"  $\IP<(Q_j(u))_x,Q_j(u)> = 0$.
\item "{e\rom{)}}" $\PB(H_j,H_{j-1}) = 0$.
\endroster
\endproclaim

\proof
  Statement a) is just a special case of the $\ad$ invariance of
the Killing form, and b) is a special case of a).

 Recalling that
$<u_1,u_2> = -\int_{-\infty}^\infty \tr(u_1,u_2)\, dx$,
it follows that
$$\IP<(Q_j(u))_x,Q_k(u)> + \IP<Q_j(u),(Q_k(u))_x> =
-\int_{-\infty}^\infty {d \over dx}\tr(Q_j(u),Q_k(u))\, dx,$$
which is clearly zero since $\tr(Q_j(u),Q_k(u))$ vanishes at infinity.
This proves c), and d) is just a special case of c).

  Since $\PB(H_j,H_{j-1}) = \IP<[\a,Q_{j+1}(u)],Q_j(u)>$, the
recursion formula for $[\a,Q_{j+1}(u)]$ gives
$\PB(H_j,H_{j-1}) = \IP<(Q_j(u))_x,Q_j(u)> - \IP<[u,Q_j(u)],Q_j(u)>$,
and e) now follows from b) and d).
\qed

\proclaim {Lemma 2}
   $\PB(H_k,H_l) = -\PB(H_{k-1}, H_{l+1})$.\endproclaim
\proof
 $\PB(H_k,H_l) = \IP<[\a, Q_{k+1}(u)], Q_{l+1}(u)>$, so that using the
recursion formula for $[\a, Q_{k+1}(u)]$ we find: \hfil\break
\centerline{$\PB(H_k,H_l) = \IP<(Q_k(u))_x, Q_{l+1}(u)> -
   \IP<[u,Q_k(u)], Q_{l+1}(u)>$,}
and  using a) of Lemma 1,\hfil\break
\centerline{$\PB(H_k,H_l) = \IP<(Q_k(u))_x, Q_{l+1}(u)> +
   \IP<(Q_k(u), [u,Q_{l+1}(u)]>$.}
Next, using the recursion formula for $[a,Q_{l+2}(u)]$,
we find that $\PB(H_k,H_l) = \hfil\break
\IP<(Q_k(u))_x, Q_{l+1}(u)> +
    \IP<(Q_k(u), (Q_{l+1}(u))_x> - \IP<(Q_k(u), [a,Q_{l+2}(u)]>$,
and we recognize the third term as $-\PB(H_{k+1},H_{l+1})$, while
the sum of the first two terms vanishes by c) of Lemma 1. \qed

  The proof that $\PB(H_k,H_l) = 0$ for any $k$ and $l$ is now easy.
We can suppose that $k\ge l$, and we apply Lemma 2 repeatedly,
decreasing the larger index by one and increasing the smaller by
one, until we ``meet in the middle''. At this point we have an identity
$\PB(H_k,H_l) = \pm \PB(H_m,H_n)$ where $m = n$ if $k$ and $l$ have
the same parity, while $m = n+1$ if the have opposite parity. In the
first case we get $\PB(H_k,H_l) = 0$ by the anti-symmetry of
Poisson Brackets, and in the second case $\PB(H_k,H_l) = 0$ by
e) of Lemma 1.

  This finishes our partial proof of the NLS Hierarchy Theorem; 
we will complete the proof later.
\goodbreak

\section {ZS-AKNS Direct Scattering Theory}

\subsection {Statements of Results}

\noindent
  For each potential $u$ in our phase space $\S(\reals,\tperp)$
we would like to define {\it scattering data\/}, by which we will
mean a measure of the asymptotic behavior of solutions of the
parallel transport equation,
$\psi_x = A(u,\lambda)\psi = (\a \lambda + u)\psi$,
for $x$ near $\pm \infty$.
Of course, to have a useful Inverse Scattering Method, the
scattering data for $u$ must be such that it allows us to
recover $u$. On the other hand, it is preferable to make the
scattering data as simple as possible, so it should be
``just enough'' to recover $u$. Direct Scattering Theory
refers to this search for such good minimal scattering data,
and for the explicit determination of the image of the Direct
Scattering Transform, (the map from $u \in \S(\reals,\tperp)$
to the scattering data of $u$). Identifying this image precisely is
of course essential for a rigorous definition of the Inverse
Scattering Transform that recovers $u$ from its scattering data.

   It turns out that, in discussing the asymptotic behavior of
solutions $\psi$ of the parallel transport equation near infinity,
it is more convenient to deal not with $\psi$ itself, but rather with
the related function $\phi = \psi(x)e^{-\a\lambda x}$, which
satisfies a slightly modified equation.

\proclaim {Proposition 1} If $\psi$ and $\phi$ are maps of $\reals$
into $\SLn(2,\Cx)$ that are related by
$\phi(x) = \psi(x)e^{-\a\lambda x}$,
then $\psi$ satisfies the parallel transport equation,
$\psi_x = (\a \lambda + u)\psi$, if and only if
$\phi$ satisfies what we shall call the
``modified parallel transport equation'',
$\phi_x = [\a\lambda,\phi] + u\phi$. \endproclaim

\proof
Clearly
$\phi_x = \psi_x e^{-\a\lambda x} -\psi  e^{-\a\lambda x}\a\lambda
   = (\a \lambda + u)\psi e^{-\a\lambda x} - \phi \a\lambda $,
and the result follows.  \qed

\smallskip\noindent
{\bf Definition.} For $u$ in $\S(\reals,\tperp)$, we will call
$m^u(x,\lambda)$ a {\it normalized eigenfunction of $u$
with eigenvalue $\lambda$\/}  if it satisfies the modified
parallel transport equation,
$m^u_x  = [\a\lambda,m^u] + u m^u$,
and if in addition:
\roster
\item "{1)}" $\lim_{x \to -\infty} m^u(x,\lambda) = I$.
\item "{2)}" $\sup_{x\in \reals}\norm|m^u(x,\lambda)| < \infty$
\endroster

\smallskip\goodbreak
It is these normalized eigenfunctions $m^u$ that will play the
r\^ole of scattering data in this theory; they are analogous to
the Jost solutions of the Schr\"odinger equation in the KdV theory.
Note that condition 2) just means that each matrix element of
$m^u(x,\lambda)$ is a bounded function of $x$.

   A complete theory of normalized eigenfunctions will be found in
[BC1]\null.	  We will next state the basic results proved there as three
theorems, Theorem A, Theorem B, and Theorem C,  reformulating
things somewhat so as to make the statements better adapted to the
Terng-Uhlenbeck version of inverse scattering theory that we will
explain later. Then we will sketch the proofs of these results,
leaving it to the interested reader to fill in many of the details
from the original paper of Beals and Coifman.

  We will denote $\S(\reals,\tperp)$ by $P$ in what follows.

\proclaim {Theorem A} For each $u$ in $P$ there is a unique normalized
eigenfunction $m^u(x,\lambda)$ for $u$ with eigenvalue $\lambda$,
except for $\lambda$ in $\reals \cup D^u$, where $D^u$ is
a bounded, discrete subset of $\Cx\setminus\reals$.
Moreover, as a function of $\lambda$, for each fixed $x$ in $\reals$,
$m^u(x,\lambda)$ is meromorphic in $\Cx\setminus\reals$
with poles at the points of $D^u$.\endproclaim
\noindent
Note that a matrix-valued function of a complex variable is said to be
holomorphic (resp., meromorphic) in a region $O$ if each of its matrix
elements is holomorphic (resp., meromorphic) in $O$, and a pole of
such a function is  a pole of any of its matrix elements.

\noindent
{\bf Definition.} An element $u$ of $P$ will be called a {\it regular
potential\/} if $D^u$ is a finite set and if, for all real $x$, the
function $m^u(x,\lambda)$ with $\lambda$ in the upper half-plane
$\Cx_+$ has smooth boundary values $m^u_+(x,r)$ on the real axis,
and similarly $m^u(x,\lambda)$ with $\lambda$ in the lower half-plane
$\Cx_-$ has smooth boundary values $m^u_-(x,r)$.
We will denote the set of regular potentials by $\Preg$.
\smallskip\noindent

\proclaim {Theorem B} The space $\Preg$ of regular potentials is open
and dense in the space $P = \S(\reals,\tperp)$ of all potentials.\endproclaim

  It is an essential fact that the normalized eigenfunctions
$m^u(x,\lambda)$ have asymptotic expansions as $|\lambda|$ tends
to infinity. Since the precise nature of these expansions will be
important, we will give the relevant definitions in some detail

  A matrix-valued function $f(\lambda)$ defined for complex
$\lambda$ with $|\lambda|$ sufficiently large is said to have
an asymptotic expansion at infinity if there exists a sequence of
matrices $f_n$ so that
$f(\lambda) - \sum_{j=0}^k f_j\lambda^{-j}
   = o(|\lambda|^{-k})$. It is easy to see inductively that
the $f_n$ are uniquely determined,
and we write $f \sim \sum_j f_j\lambda^{-j}$.

  Now suppose that we have matrix-valued functions $f(x,\lambda)$,
defined for all $x$ in $\reals$ and all $\lambda$ in $\Cx$ with
$|\lambda|$ sufficiently large. Suppose that we have matrix-valued
functions $f_n(x)$ such that for each $x$, $f(x,\lambda) \sim \sum_j
f_j(x)\lambda^{-j}$. We will write $f \simR \sum_j f_j\lambda^{-j}$
if this asymptotic expansion holds {\it uniformly in $x$\/}, i.e., if
$$\sup_x\norm|f(x,\lambda) - \sum_{j=0}^k f_j(x)\lambda^{-j}|
   = o(|\lambda|^{-k}).$$

  It is easy to explain the importance of the uniformity. Suppose
$f$ and the $f_n$ are differentiable functions of $x$. Then the
uniformity gives
$${f(x+\Delta x,\lambda) - f(x,\lambda)\over \Delta x}
   - \sum_{j=0}^k {f_j(x+\Delta x)-f_j(x)\over \Delta x}\lambda^{-j}
  = o(|\lambda|^{-k})$$
and letting $\Delta x$ approach zero gives
${\D f\over \D x} \simR \sum_j f^\prime_j\lambda^{-j}$,
i.e., we can differentiate such an asymptotic relation ``term by term''.

\proclaim {Theorem C} For $u$ in $\Preg$, the normalized eigenfunctions
$m^u(x,\lambda)$ have an asymptotic expansion as $\lambda$ tends
to infinity, $m^u  \simR \sum_j m^u_j\lambda^{-j}$.
In fact the $m^u_j$ are uniquely determined
inductively by the condition
$[\a,m^u_{j+1}(x)] = {d\over dx} m^u_j(x) - u(x) m^u_j(x)$.
\endproclaim

   The normalized eigenfunctions, $m^u(x,\lambda)$, satisfy a
simple relation, referred to as the ``reality condition'' that
follows as an easy consequence of the fact that $u(x)$ takes
its values in $\sun(2)$.

\proclaim {Proposition 2} If $u\in \Preg$ then the normalized
eigenfunctions $m^u$ satisfy the relation
$m^u(x,\bar\lambda)^* m^u(x,\lambda) = I$.\endproclaim

   So, passing to the limit as $\lambda \in \Cx_+$
approaches $r \in \reals$,

\proclaim {Corollary} $m^u_-(x,r)^* m^u_+(x,r) = I$.\endproclaim

  We will need one more property of the $m^u$ (or rather of their
boundary values, $m^u_\pm$).

\proclaim {Proposition 3} Let $u\in \Preg$ and $x \in \reals$,
and let $m^u_+(x,r) = g(x,r)h(x,r)$ be the canonical decomposition
of $m^u_+(x,r)$ into the product of a unitary matrix $g(x,r)$
and an upper-triangular matrix $h(x,r)$. Then $h(x,r) - I$ is
of Schwartz class in $r$.\endproclaim

\subsection {Outline of Proofs}

\noindent
   As was the case for the scattering theory for the
Schr\"odinger operator, it is a lot easier to see what
is happening for the special case of potentials with compact
support. It turns out for example that all such potentials
are regular. Below we will give most of the details of the
proofs of Theorems A, B, and C for the $2 \times 2$ case when
$u$ has compact support.

[In [BC1], the case of compactly supported potentials is
considered first, followed by the case of ``small potentials'',
i.e., those with $L^1$ norm less than $1$. For the latter,
it turns out that existence and uniqueness of the $m^u$
can be proved easily using the Banach Contraction Principle,
and moreover it follows that $D^u$ is empty. The case of
regular potentials (called ``generic'' in [BC1]) is then
handled by a limiting argument. [BC1] also consider the general
$n \times n$ case and does not assume that $u$ is necessarily
skew-adjoint. This latter generality adds substantial extra
complexity to the argument.]

  In any interval $[a,b]$ in which $u$ vanishes identically, the
modified parallel transport equation reduces to the Lax Equation
$\phi_x = [\a\lambda,\phi]$, so choosing an arbitrary $x_0$ in $[a,b]$,
the solution is
$\phi_(x) = e^{\a\lambda (x-x_0)}\phi(x_0) e^{-\a\lambda (x-x_0)}$, or
$\phi_(x) = e^{\a\lambda x}s e^{-\a\lambda x}$, where we define
$s = e^{-\a\lambda x_0}\phi(x_0) e^{\a\lambda x_0}$.
 This proves:

\proclaim {Proposition 4} Suppose $u$ in $P$ has compact support,
say $u(x) = 0$ for $|x| \ge M$. Then for each complex number $\lambda$
there is a unique solution $\phi^u(x,\lambda)$ of the modified
parallel transport equation with $\phi^u(x,\lambda) = I$ for $x \le -M$.
Moreover, for $x \ge M$, $\phi^u$ has the form
$\phi^u_(x.\lambda) = e^{\a\lambda x}s^u(\lambda) e^{-\a\lambda x}$
(where $s^u(\lambda) = e^{-\a\lambda M}\phi^u(M,\lambda) e^{\a\lambda M}$),
and for each real $x$, $\lambda \mapsto \phi^u(x,\lambda)$ is an
entire function (i.e., holomorphic in all of $\Cx$).
\endproclaim

  The fact that $\phi^u$ is holomorphic in $\lambda$ is a consequence
of the more general principle that if an ODE depends analytically
on a parameter $\lambda$, then the solution of the equation with
some fixed initial condition is analytic in $\lambda$. (In this case
the initial value condition is $\phi^u(-M,\lambda) = I$.)

\smallskip\noindent
{\bf Definition.} We will denote the matrix elements of $s^u(\lambda)$
by $s^u_{ij}(\lambda)$, and we define $D^u$ to be the set of all
$\lambda$ in the upper half-plane that are zeroes of $s^u_{11}$
union the set of all $\lambda$ in the lower half-plane that are
zeroes of $s^u_{22}$.

\smallskip\noindent
{\bf Remark.} It can be shown that the holomorphic functions $s^u_{12}$
and $s^u_{21}$ are not identically zero, so that $D^u$ is a discrete
set. In fact (cf.\ [BC1], section 4), $D^u$ is finite, and neither
$s^u_{11}$  nor $s^u_{22}$ has any zeroes on the real axis.
\smallskip

\proclaim {Proposition 5} Suppose $u$ in $P$ has compact support.
For each
$\lambda \in  \Cx \setminus (\reals \cup D^u)$
there is a unique
normalized eigenfunctions $m^u(x,\lambda)$.  For every $x$ in
$\reals$, $m^u(x,\lambda)$ is a meromorphic function of $\lambda$
for $\lambda$ in $\Cx\setminus \reals$, with poles at the points
of $D^u$. Finally, the restriction of $m^u(x,\lambda)$ to each
half-plane has a smooth extension to the real axis.
\endproclaim\goodbreak
\proof
 Since $\phi^u(x,\lambda)$ is invertible, there is no loss of
generality in assuming that a normalized eigenfunction has the
form $m^u(x,\lambda) = \phi^u(x,\lambda) \chi^u(x,\lambda)$. Then
$[\a\lambda,m^u] + u m^u = m^u_x = \phi^u_x \chi^u + \phi \chi^u_x$,
which simplifies to the same Lax Equation as before, namely
$\chi^u_x = [\a\lambda,\chi^u]$, but now valid on the whole of $\reals$,
and it follows that
$\chi^u(x,\lambda) = e^{\a\lambda x}\chi^u(\lambda) e^{-\a\lambda x}$,
and hence
$m^u(x,\lambda) =
   \phi^u(x,\lambda) e^{\a\lambda x}\chi^u(\lambda) e^{-\a\lambda x}$.

  Then, by Proposition 4, for $x \le -M$,
$m^u(x,\lambda) = e^{\a\lambda x}\chi^u(\lambda) e^{-\a\lambda x}$
while for $x \ge M$,
$m^u(x,\lambda) =
   e^{\a\lambda x}s^u(\lambda)\chi^u(\lambda) e^{-\a\lambda x}$.

  Let us write $\chi^u_{ij}(\lambda)$ for the matrix elements of
$\chi^u(\lambda)$, and try to determine them individually so that
Conditions 1) and 2) of the definition of generalized eigenfunctions
will be satisfied for the resulting $m^u(x,\lambda)$.

  Note that, since conjugating $\chi^u(\lambda)$ by a diagonal
matrix does not change its diagonal entries, the diagonal elements
of $m^u(x,\lambda)$ are just $\chi^u_{11}(\lambda)$ and
$\chi^u_{22}(\lambda)$ for $x \le -M$. Since Condition 1)
requires that $m^u(x,\lambda)$ converge to the identity matrix
as $x$ approaches $ -\infty$, it follows that we must take
$\chi^u_{11}(\lambda) = \chi^u_{22}(\lambda) = 1$, and conversely
with this choice Condition 1) is clearly satisfied.

  On the other hand, an easy calculation shows that the off-diagonal
elements, $m_{12}^u(x,\lambda)$ and $m_{21}^u(x,\lambda)$, are
given respectively by $e^{-2i\lambda x}\chi^u_{12}(\lambda)$
and $e^{2i\lambda x}\chi^u_{21}(\lambda)$, when $x \le -M$. If
$\lambda = \sigma + i\tau$,
$m_{12}^u(x,\lambda) =
   e^{-2i\sigma x}e^{2\tau x}\chi^u_{12}(\lambda)$,
and
$m_{21}^u(x,\lambda) =
   e^{2i\sigma x}e^{-2\tau x}\chi^u_{21}(\lambda)$.
Since Condition 2) requires that these remain bounded when
$x$ approaches $-\infty$, it follows that when $\lambda$ is
in the lower half-plane (i.e., $\tau < 0$) then
$\chi^u_{12}(\lambda) = 0$, and similarly,
$\chi^u_{21}(\lambda) = 0$ for $\lambda$ in the upper half-plane.

  Next, take $x > M$, so that
$m^u(x,\lambda) =
    e^{\a\lambda x}s^u(\lambda)\chi^u(\lambda) e^{-\a\lambda x}$.
Then another easy computation shows that if $\lambda$ is in the upper
half-plane, then
$m_{12}^u(x,\lambda) =
 e^{-2i\lambda x}(
       s_{11}^u(\lambda)\chi^u_{12}(\lambda) +s_{11}^u(\lambda))$,
while $m_{12}^u(x,\lambda) = 0$.
Since $m_{11}^u(x,\lambda) = s_{11}^u(\lambda)$ and
$m_{22}^u(x,\lambda) = s_{22}^u(\lambda)$ are independent of $x$,
the condition for $m^u(\lambda)$ to remain bounded when $x$ approaches
$+\infty$  is just
$s_{11}^u(\lambda)\chi^u_{12}(\lambda) + s_{12}^u(\lambda) = 0$, and
this uniquely determines $\chi_{12}^u(\lambda)$, namely
$\chi_{12}^u(\lambda) = - s_{12}^u(\lambda)/s_{11}^u(\lambda)$.
So for $\lambda$ in the upper half-plane
$\chi^u(\lambda) =
\mat(1;- s_{12}^u(\lambda)/s_{11}^u(\lambda);0;1)$
is the unique choice of $\chi^u$ satisfying Conditions 1) and 2).
A similar computation shows that for $\lambda$ in the lower half-plane
$\chi^u(\lambda) =
\mat(1;0;s_{21}^u(\lambda)/s_{22}^u(\lambda);1)$. All 
conclusions of the proposition follow from these explicit formulas
and the fact that $s_{11}$ and $s_{22}$ have no zeroes on the
real axis.
\qed

\proclaim {Lemma} If $\psi_x = A\psi$ and $\phi_x = -\phi A$
then $\phi \psi$ is constant.\endproclaim
\proof
$(\phi \psi)_x =\phi_x \psi + \phi \psi_x = 0.$ \qed

  We can now prove Proposition 2.

\proof  It will suffice to prove that
$m^u(x,\bar\lambda)^* m^u(x,\lambda)$ is constant, since
we know that as $x$ approaches $-\infty$ the product
converges to $I$. If we define
$\psi(x,\lambda) = m^u(x,\lambda)e^{\a\lambda x}$,
then $\psi(x,\bar\lambda)^* = e^{-\a\lambda x}m^u(x,\bar\lambda)$,
and therefore
$m^u(x,\bar\lambda)^* m^u(x,\lambda)
= \psi(x,\bar\lambda)^* \psi^u(x,\lambda)$
and it will suffice to prove that
$\psi(x,\bar\lambda)^* \psi^u(x,\lambda)$
is constant. By Proposition 1,
$\psi_x(x,\lambda) = (\a\lambda + u)\psi(x,\lambda)$.
Since $u^* = -u$ and $(\a\bar\lambda)^* = -\a\lambda$,
$\psi_x(x,\bar\lambda)^* = \psi(x,\bar\lambda)^*(\a\bar\lambda + u)^*
= -\psi(x,\bar\lambda)^*(\a\lambda + u)$, and the preceding lemma
completes the proof.  \qed

  Our Theorem C is just Theorem 6.1, page 58, of [BC1]\null. While the
proof is not difficult, neither is it particularly illuminating,
and we will not repeat it here. Similarly, our Proposition 3
follows from Theorem E${}^\prime$, page 44 of [BC1]\null.

   This completes our discussion of the proofs of Theorems A, B, C,
and Propositions 2 and 3. In the remainder of this section we will
see how these results can be used to complete the proof of the
NLS Hierarchy Theorem.

 Since $m^u(x,\lambda)^{-1} = m^u(x,\bar\lambda)^*$, it follows that
 $m^u(x,\lambda)\a (m^u(x,\lambda))^{-1}$
has an asymptotic expansion.

\noindent
{\bf Definition.} We denote the function
$m^u(x,\lambda)\a (m^u(x,\lambda))^{-1}$ by
$Q^u(x,\lambda)$.

\smallskip\noindent
So by the preceding remark,

\proclaim {Corollary} $Q^u(x,\lambda)$ has an asymptotic expansion
$Q^u \simR \sum_{j=0}^\infty Q^u_j\lambda^{-j}$, with $Q^u_0 = \a$,
hence also
$Q^u_x \simR \sum_{j=0}^\infty (Q^u_j)_x\lambda^j$
\endproclaim

\proclaim {Lemma} If we define
$\psi(x,\lambda) = m^u(x,\lambda)e^{\a\lambda x}$ then
$Q^u(x,\lambda) = \psi\a \psi^{-1}$. \endproclaim
\proof Immediate from the fact that all diagonal matrices commute.  \qed

  Now
$(\psi\a \psi^{-1})_x = \psi_x \a\psi^{-1} + \psi\a(\psi^{-1})_x$,
and by Proposition 1, $\psi_x = (\a\lambda + u) \psi$. Also, from
$\psi \psi^{-1} = I$ we get
$\psi_x \psi^{-1} + \psi (\psi^{-1})_x = 0$. Combining all
these facts gives
$(\psi\a \psi^{-1})_x = [\a \lambda + u, \psi\a \psi^{-1}]$,
and hence, by the lemma,
$Q^u_x(x,\lambda) = [\a \lambda + u, Q^u(x,\lambda)]$.
If we insert in this identity the asymptotic expansion
$Q^u \simR \sum_{j=0}^\infty Q^u_j\lambda^j$
we find a second asymptotic expansion for
$Q^u_x(x,\lambda)$, in addition to the one from the above
Corollary, namely
$Q^u_x \simR
   \sum_j ([a,Q^u_{j+1}] + [u,Q^u_j])\lambda^j$.
Therefore, by uniqueness of asymptotic expansions we have proved:

\proclaim {Proposition 6} The recursion relation
$(Q^u_j)_x = [a,Q^u_{j+1}] + [u,Q^u_j]$, is satisfied by the
coefficients $Q^u_j$ of the asymptotic expansion of
$Q^u(x,\lambda)$,
and hence they are identical with the functions
$Q_j(u) :\reals \to \sun(2)$ defined in the NLS Hierarchy Theorem.
\endproclaim

  We are now finally in a position to complete the proof of the
NLS Hierarchy Theorem.

  Since $\a^2 = I$, it follows that also
$Q^u(x,\lambda)^2 = (m^u\a (m^u)^{-1})^2 = I$, and hence
$I \sim (\sum_{j=0}^\infty Q^u_j\lambda^j)^2$.
Expanding and comparing coefficients of $\lambda^{-k}$,
uniqueness of asymptotic expansions gives
$\a Q_k(u) + Q_k(u)\a = -\sum_{j = 1}^{k-1} Q_j(u) Q_{k-j}(u)$.
Recall that we needed one fact
to complete the proof of statements a) through d)
of the NLS Hierarchy Theorem, namely that
if $Q_k(u) = P_k(u) + T_k(u)$ is the decomposition
of $Q_k(u)$ into its off-diagonal part and its diagonal part,
then the matrix elements of $T_k(u)$ are polynomials in the
matrix elements of $u$ and their derivatives. Moreover, we saw
that we could assume inductively that this was true for the
matrix elements of $Q_j(u)$ for $j < k$.
But if $T_k =\mat(t_k;0;0;-t_j)$, then
$\a Q_k(u) + Q_k(u)\a = -2iT_k =\mat(-2it_k;0;0;-2it_j)$
and the desired result is now immediate from the inductive
assumption.

  The other statement of the NLS Hierarchy Theorem that
remains to be proved is e).

  Define a function $\tilde F^u(x,\lambda) = \tr(Q^u(x,\lambda) \a)$.
Clearly $\tilde F^u(x,\lambda)$ has an asymptotic expansion,
$\tilde F^u \simR \sum_j \tilde F^u_j\lambda^{-j}$, where
$\tilde F^u_j = \tr(Q^u_j\a)$.

  From what we have just seen, $\tilde F^u(x,\lambda)$
is Schwartz class in $x$, so we can define a map
$F(u,\lambda) = \int_{-\infty}^\infty \tilde F^u(x,\lambda)\, dx
=\int_{-\infty}^\infty\tr(Q^u \a)\, dx$, and
$F(u,\lambda) \sim \sum_j F_j(u) \lambda^{-j}$ where
$F_j(u) =
\int_{-\infty}^\infty  \tilde F^u_j(x)\,dx$.

  If we consider $u \mapsto F(u,\lambda)$ as a function on
$P$, then $\grad F$ is a vector field on $P$, and
$(\grad F)_u \sim \sum_j (\grad F_j)_u \lambda^{-j}$.
We claim that statement e) of The NLS Hierarchy Theorem
follows from the following proposition. (For a proof of which,
see Proposition 2.4 of [Te2]\null.)

\proclaim {Proposition 7} If $v$ in $\Preg$ then
$${d\over d\epsilon}\Bigl|_{\epsilon = 0}
F(u+\epsilon v,\lambda) =
\int_{-\infty}^\infty
\tr\left({dQ^u(x,\lambda)\over d\lambda} v(x)\a\right)\,dx.$$
\endproclaim

\noindent
Indeed, expand both sides of the latter equality
in asymptotic series in $\lambda$,
and compare coefficients of $\lambda^{-j}$.  Since
${dQ^u(x,\lambda)\over d\lambda} = \sum_j -jQ^u_j\lambda^{-j-1}$,
we find
$(dF_j)_u(v) = \int_{-\infty}^\infty
\tr((-(j-1)Q_{j-1}(u)(x) v(x)\a)\,dx.$
Recalling the definition of the inner product in $P$, we see
$-{1\over j-1}(\grad F_j)_u$ is the projection of
$Q_{j-1}(u)$ on $\tperp$, i.e., the
off-diagonal part of $Q_{j-1}(u)$.
So if we define
$H_j(u) =  -{1\over j+1}\int_{-\infty}^\infty\tr((Q_{j+2}(u) \a)\, dx
= -{1\over j+1}F_{j+2}(u)$, then
$(\grad H_j)_u =  -{1\over j+1}(\grad F_{j+2})_u$ is the
off-diagonal part of $Q_{j+1}(u)$, which is statement e)
of The NLS Hierarchy Theorem.

\section {Loop Groups, Dressing Actions, and Inverse Scattering}

\subsection {Secret Sources of Soliton Symmetries}

\noindent
  This article is titled  ``The Symmetries of Solitons'', and
we have been hinting that many of the remarkable properties of
soliton equations are closely  related to the existence of large
and non-obvious groups of symplectic automorphisms that act
on the phase spaces of these Hamiltonian systems and leave the
Hamiltonian function invariant. We are now finally in a position
where we can describe these groups and their symplectic actions.

  The groups themselves are so-called loop groups. While they
have been around in various supporting r\^oles for much longer,
in the past three decades they have have been increasingly studied
for their own sake, and have attained a certain prominence. See for
example [PrS]\null.

  Given any Lie group $\G$, we can define its associated loop
group, $L(\G)$ as the group of all maps (of some appropriate
smoothness class) of $\sph^1$ into $\G$, with pointwise composition.
For our purposes we will always assume that $\G$ is a matrix group
and the maps are smooth (i.e., infinitely differentiable).

  The theory gets more interesting when we regard the loops in $\G$
as boundary values of functions that are holomorphic (or meromorphic)
in the interior (or exterior) of the unit disk and take subgroups
by restricting the analytic properties of these analytic extensions.
That is, we concentrate on the analytic extensions rather than the
boundary value.

   Once we take this point of view, it is just as natural to
pre-compose with a fixed linear fractional transformation mapping
the real line to the unit circle, (say $z \mapsto (1 + iz)/(1 - i z)$),
so that elements of the loop groups become maps of $\reals$ into $\G$
that are boundary values of certain analytic functions in the upper or
lower half-plane, and this is the point of  view we will adopt. Note
that the above linear fractional transformation take $-1$ in $\sph^1$
to infinity, and for certain purposes it is important to know how the
nature of the original map of $\sph^1$ into $\G$ at $-1$ translates to
properties of the transformed map of $\reals$ into $\G$ at $\pm\infty$.
A straightforward calculation gives the following answer:

\proclaim {Proposition 1} \rom{([TU], Proposition 7.7)}
   Given $g: \sph^1  \to \GLn(n,\Cx)$, define
$\Phi(g): \reals \to \GLn(n,\Cx)$ by
$\Phi(g) (r) = g({1+ i r \over 1- i r} )$. Then:
\roster
\item "{\rom{(}i\rom{)}}" g is smooth if and only if $\Phi(g)$ is smooth
            and has asymptotic expansions at $+\infty$ and at
            $-\infty$ and these expansions agree.
\item "{\rom{(}ii\rom{)}}" $g - I$ is infinitely flat at $z = -1$ if and only
             if $\Phi(g) - I$ is of Schwartz class.
\item "{\rom{(}iii\rom{)}}" $g :\Cx \to \GLn(n,\Cx)$ satisfies the reality
               condition $g({1\over \bar z})^*g(z) = I$ if and only if
               $\Phi(g)(\lambda) = g({1+ i\lambda \over 1- i\lambda})$
               satisfies
               $\Phi(g)(\bar\lambda)^*\Phi(g)(\lambda) = I$.
\endroster
\endproclaim

   The first, and most important, loop group we will need is called
$\Dm$. The analytic properties of its elements are patterned
after those proved to hold in the preceding section for the
normalized eigenfunctions $m^u(x,\lambda)$ as functions of $\lambda$.

\noindent
{\bf Definition.} We will denote by $\Dm$ the group of all meromorphic
maps $f:\Cx\setminus \reals \to \GLn(n,\Cx)$
having the following properties:
\roster
\item "{1)}" $f(\bar\lambda)^*f(\lambda) = I$.
\item "{2)}" $f$ has an asymptotic expansion
           $f(\lambda) \sim I + f_1\lambda^{-1}+ f_2\lambda^{-2} + \cdots$.
\item "{3)}" The set $D^f$ of poles of $f$ is finite.
\item "{4)}" $f$ restricted to the upper half-plane, $\Cx_+$, extends
           to a smooth function on the closure of the upper-half plane,
           and similarly for the lower half-plane. The boundary values
           are then maps $f_{\pm}:\reals \to \GLn(n,\Cx)$, and by 1)
           they satisfy $f_+(r)^*f_-(r) = I$.
\item "{5)}" If $f_+(r) = g(r)h(r)$ is the factorization of $f_+(r)$
           as the product of a unitary matrix $g(r)$ and an upper
           triangular $h(r)$, then $h-I$ is of Schwartz class.
\endroster

\smallskip\noindent
{\bf Definition.}
  We define a map, $\Fs:P_0 \to \Dm$, the {\it Scattering
Transform\/},  by $\Fs(u)(\lambda) = f^u(\lambda) = m^u(0,\lambda)$.
\smallskip\noindent
  That $m^u(0,\lambda)$ is in fact an element of $\Dm$
is a consequence of the definition of the set $P_0$ of regular
potentials, and Theorems A and C and Propositions 2 and 3 of the
preceding section. There is nothing special about $0$ in the
above definition. We could have equally well chosen any other
fixed real number $x_0$ and used $m^u(x_0,\lambda)$ instead of
$m^u(0,\lambda)$.

\subsection {Terng-Uhlenbeck Factoring and the Dressing Action}

\noindent
  There are three other loop groups that play an essential r\^ole in
the definition of the Inverse Scattering Transform, $\IFs$, and we
define these next.

\noindent
{\bf Definition.}
We will denote by $\Gp$ the loop group of all entire functions
$h:\Cx \to \GLn(n,\Cx)$, and by $\Hp$ the abelian subgroup of $\Gp$
consisting of all elements of the form $e^{\a P(\lambda)}$
where $P:\Cx \to \Cx$ is a polynomial in $\lambda$. Finally, we
define $\Hm$ to be the subgroup of $\Dm$ consisting of those
elements $f$ taking values $f(\lambda)$ in the
diagonal subgroup of $\GLn(n,\Cx)$.
For each $x$ in $\reals$ we define $e_\a(x)$ in $\Hp$ by
$e_\a(x)(\lambda) = e^{\a\lambda x}$, and for each positive
integer $j$ we define a one-parameter subgroup $e_{\a,j}$
of $\Hp$ by $e_{\a,j}(t) = e^{\a\lambda^j t}$.
(Note that $e_\a(x) = e_{\a,1}(x)$.)

\smallskip\noindent
  The following theorem is one of the basic results of [TU]\null.
As we shall see, it provides an alternative, group theoretic
approach to ZS-AKNS Inverse Scattering Theory. (In fact,
conversely, it can be proved using earlier approaches to
ZS-AKNS Inverse Scattering Theory).

\proclaim {Terng-Uhlenbeck Factoring Theorem} \rom{([TU], 7.11 
and 7.16)} If $f\in\Dm$ then:
\roster
\item "{1)}" for any $h\in \Hp$,
           $hf^{-1}: \Cx\setminus (\reals \cup D^f) \to \GLn(n,\Cx)$
           can be factored uniquely in the form
           $hf^{-1} = M^{-1}E$, with $M$ in $\Dm$ and $E$ in $\Gp$.
\item "{2)}" Taking $h = e_{\a,1}(x)$ in 1)
           (i.e., $h(\lambda) = e^{\a\lambda x}$),
           we get a one parameter family of such factorings,
           $e_{\a,1}(x) f^{-1} = M^{-1}(x)E(x)$ and, writing
           $E_x$ for the derivative of $E$, it follows that
           $E_x = (\a\lambda + u)E$ for a unique, regular
           potential $u$ in $P_0$.
\endroster
\endproclaim

\noindent
   We note that in 1) uniqueness is easy and only existence of the
decomposition needs proof. Indeed, uniqueness is equivalent to the
statement that $\Dm \cap \Gp =I$, and this is immediate from
from Liouville's Theorem that bounded holomorphic functions are
constant (recall that elements of $\Dm$ converge to $I$ as
$\lambda \to \infty$). The existence part of 1) follows from the 
two classical Birkhoff Decomposition Theorems, and statement 2) gives
the dependence of this factorization on the parameter $x$.

\noindent
{\bf Definition.} We define a left action of $\Hp$ on $\Dm$,
called the {\it dressing action\/}, and denoted by
$(h,f) \mapsto h*f$. It is defined by $h*f = M$, where $M$ is
given by the factoring of $hf^{-1}$ in 1) of the previous theorem.

\smallskip\noindent
  Of course we must check that $(h_1 h_2)*f = h_1*(h_2*f)$,
but this is easy. Suppose $h_2 f^{-1} = M_2^{-1}E_2$,
i.e., $h_2*f = M_2$, and use the factoring theorem again to
write $h_1 M_2^{-1}$ as a product,
$h_1 M_2^{-1} = M_1^{-1} E_1$, i.e., $h_1*M_2 = M_1$.
Then $(h_1 h_2)f^{-1} = h_1(h_2 f^{-1}) = h_1 M_2^{-1}E_2
= M_1^{-1} E_1 E_2$, so $(h_1 h_2) * f = M_1 =
h_1*M_2=h_1*(h_2*f)$.

  Now that we have an action of $\Hp$ on $\Dm$, it follows
that every one-parameter subgroup of $\Hp$ defines a flow
on $\Dm$. In particular the one-parameter subgroups $e_{\a,j}$
define an important sequence of flows on $\Dm$.

\noindent
{\bf Definition.}  For each positive integer $j$ we define a
flow on $\Dm$, called {\it the $j$-th flow\/}, by
$(t,f) \mapsto e_{\a,j}(t)*f$.

\smallskip\noindent
  Of course, since $\Hp$ is an abelian group and all the
$e_{\a,j}$ are one-parameter subgroups of $\Hp$, it follows that
this sequence of flows all mutually commute.

\subsection {The Inverse Scattering Transform}

\noindent
We are now in a position to define the Inverse Scattering
Transform.

\noindent
{\bf Definition.} We define a map $\IFs:\Dm \to P_0$,
called {\it Inverse Scattering Transform\/}, by associating
to $f$ in $\Dm$ the regular potential $u = \IFs(f)$ in
$P_0$ given by 2) of the Terng-Uhlenbeck Factoring Theorem.
That is, if we define
$\psi(x,\lambda) = (e_\a(x)*f)(\lambda) e^{\a\lambda x}$, then
$u$ is characterized by the fact that $\psi$ satisfies the
parallel transport equation with potential $u$,
$\psi_x = (\a\lambda + u)\psi$.

\smallskip\noindent
\proclaim {Theorem D} The maps $\Fs: P_0 \to \Dm$ and
$\IFs : \Dm \to P_0$ satisfy:
\roster
\item "{a)}" $\IFs \circ \Fs = {\roman {identity}}$.
\item "{b)}" $\Fs \circ \IFs(f) \in f\Hm$.
\endroster
Thus, the map $P_0 \to \Dm \to  \Dm/\Hm$ that
is the composition of $\Fs$ and the natural projection
of $\Dm$ on $\Dm/\Hm$ is a bijection. \endproclaim

   Recall that in the NLS-Hierarchy Theorem we defined a sequence
of flows on $P_0$, the $j$-th of which we also called the
``$j$-th flow''. As you probably suspect:

\proclaim {Theorem E} \rom{([TU] Theorem 8.1)} The transforms $\Fs: P_0 \to \Dm$
and $\IFs : \Dm \to P_0$ are equivariant with respect to
the $j$-th flow on $\Dm$ and the $j$-th flow on $P_0$.
In particular if $u(t)$ in $P_0$ is a solution of the
$j$-th flow, then $\Fs(u(t)) = e_{\a,j}(t)*\Fs(u(0))$. \endproclaim

\proclaim {Corollary}  The following algorithm finds the
solution $u(x,t)$ for the $j$-th flow in $P_0$ with initial
condition $u(0) =u(x,0)$:
\roster
  \item "{1)}" Compute the parallel translation operator
           $\psi(x,0,\lambda)$ having the correct asymptotic
           behavior. That is, solve the following {\bf linear}
           ODE problem:

          {a\rom{)} $\psi_x(x,0,\lambda) =
                      (\a\lambda + u(x,0))\psi(x,0,\lambda)$

          {b\rom{)} $\lim_{x\to -\infty}
                          \psi(x,0,\lambda)e^{-\a\lambda x} = I$.

          {c\rom{)} $\psi(x,0,\lambda) e^{-\a\lambda x}$ is bounded.

\item "{2)}"  Define $f$ in $\Dm$  by $f(\lambda) = \psi(0,0,\lambda)$.
\item "{3)}" Factor $e_{a,j}(t)e_{a,1}(x)f^{-1}$ as 
           $M(x,t)^{-1}E(x,t)$
           where $M(x,t) \in \Dm$ and $E(x,t) \in \Gp$.
\item "{4)}" Then, putting $\psi(x,t,\lambda) = 
             M(x,t)(\lambda) e^{\a\lambda x + \lambda^j t}$,\hfil\break
             $u(x,t) =
               \psi_x(x,t,\lambda) \psi^{-1}(x,t,\lambda) - \a \lambda$.
         (The RHS is independent of $\lambda$.)
\endroster
\endproclaim
\proof  This just says that $u(t) = \IFs(e_{\a,j}(t) * \Fs(u(0)))$. \qed

\subsection {ZS-AKNS Scattering Coordinates}

\noindent
An important ingredient of the KdV Inverse Scattering Method,
based on the Schr\"odinger operator, was the that the
``coordinates'' of the scattering data evolved by a linear
ODE with constant coefficients, and so this evolution could
be solved explicitly. Recall that this allowed us to derive
an explicit formula for the KdV multi-solitons. Such scattering
coordinates (or ``action-angle variables'') also exist for the
ZS-AKNS Hierarchy, and even for the more general $n \times n$
systems, but the story is somewhat more complicated in this
case and we will only outline the theory here and refer to
[ZS] and [BS] for more complete descriptions.

  Another advantage of the loop group approach is that it permits
us to factor the scattering data into discrete and continuous parts.
To a certain extent this allows us to discuss separately the
scattering coordinates and evolution of each part.

 \noindent
 {\bf Definition.}  We define two subgroups,
 $\Dmd$ and $\Dmc$, of $\Dm$, by
 $$
 \eqalign{
 \Dmc =&
     \{ f\in \Dm \mid f \ \hbox{is holomorphic in } \Cx\setminus \reals\},\  \hbox{and} \cr 
 \Dmd =&
     \{ f\in \Dm \mid f \ \hbox{is meromorphic in } \Cx\}}
 $$

\smallskip\noindent
{\bf Remark.} Since elements of $\Dm$ approach $I$ at infinity,
it follows that any $f$ in $\Dmd$ is actually meromorphic
on the whole Riemann sphere, and hence a rational function
of the form $f_{ij}(\lambda) = P_{ij}(\lambda)/Q_{ij}(\lambda)$, 
where the polynomial maps $P_{ij}$ and of $Q_{ij}$ have the 
same degrees for a given diagonal entry, and $Q_{ij}$ has larger degree
for an off-diagonal entry. For this reason, $\Dmd$ is also
referred to as the rational subgroup of $\Dm$. Also, since
$f$ satisfies the reality condition,
$f\bar\lambda)^*f(\lambda) = I$ and is holomorphic on the real
line, it follows that for $r$ in $\reals$, $f(r)^*f(r) = I$
(i.e., $f$ is unitary on $\reals$), and the boundary
values $f_+$ of $f$ from $\Cx_+$ and $f_-$ from $\Cx_-$
are equal, so that the ``jump'', $v^f(r) = f_-^{-1}(r) f_+(r)$
is the identity.

\smallskip\noindent
\proclaim {Theorem F} \rom{(TU Theorem 7.5)} Every $f$ in $\Dm$
can be factored uniquely as a product $f = hg$ where
$h\in\Dmc$ and $g\in \Dmc$. In fact the multiplication map
$\Dmc \times \Dmd \to \Dm$ is a diffeomorphism. \endproclaim
\proof This is an immediate consequence of Proposition 1
of the previous section and the following classical theorem
of G.~D.~Birkhoff.  \qed

\proclaim {Birkhoff Decomposition Theorem}
\rom{([PrS], Theorem 8.1.1)} \hfil\break
Let $L(\GLn(n,\Cx))$
denote the loop group of all smooth maps of $\sph^1$ into
$\GLn(n,\Cx)$, $\Omega\!\Un(n)$ the subgroup of all smooth
maps $g$ of $\sph^1$ into
$\Un(n)$ such that $g(-1) = I$, and $L^+(\GLn(n,\Cx))$ the subgroup
of $L(\GLn(n,\Cx))$ consisting of all $g$ that are the boundary
values of holomorphic maps of the open unit disk into $\GLn(n,\Cx)$.
Then any $f$ in $L(\GLn(n,\Cx))$ can be factored uniquely as a
product $f = gh$ where $g \in L^+(\GLn(n,\Cx))$ and
$h\in \Omega\!\Un(n)$. In fact the multiplication map
$L^+(\GLn(n,\Cx)) \times \Omega\!\Un(n) \to L(\GLn(n,\Cx))$
is a diffeomorphism.  \endproclaim

\noindent
{\bf Definition.} Given $z\in \Cx$ and an orthogonal projection
$\pi$ in $\GLn(n,\Cx)$ we define $g_{z,\pi}$ in $\Dmd$ by
$g_{z,\pi}(\lambda) = I + {z-\bar z\over \lambda - z}\pi$

\smallskip\noindent
\proclaim {Theorem G} \rom{(Uhlenbeck [U1])}
The elements $g_{z,\pi}$ for $z \in \Cx\setminus\reals$
generate the group $\Dmd$.  \endproclaim

  It follows easily from Theorem G and the Bianchi Permutability
Formula ([TU] Theorem 10.13) that at each simple pole $z$ of
an element of $\Dmd$ we can define a ``residue'', which is just
the image of a certain orthogonal projection, $\pi$. To be precise:

\proclaim {Theorem H} If $f \in \Dm$ and $z$ is a simple pole
of $f$, then there exists a unique orthogonal projection $\pi$
such that $fg_{z,\pi}^{-1}$ is holomorphic at $z$.  \endproclaim

  The set of $f$ in $\Dm$ for which all the poles are simple
is open and dense, and it is for these $f$ that we will define
``scattering coordinates'', $S^f$.
\smallskip

\noindent
{\bf Definition.} Given $f$ in $\Dm$ with only simple poles,
the {\it scattering coordinates of $f$\/}, $S^f$ consists of
the following data:
\roster
\item "{a)}" The set $D^f = \{z_1,\ldots,z_N\}$ of poles of $f$.
\item "{b)}" For each $z$ in $D^f$, the ``residue'' of $f$ at $z$,
          i.e., the image $V^f_z$ of the unique orthogonal projection, 
           $\pi =\pi^f_z$ such that $fg_{z,\pi}^{-1}$ is holomorphic at $z$.
\item "{c)}" The jump function of $f$, i.e., the map
           $v^f : \reals \to \GLn(n,\Cx)$
           defined by $v^f(r) = f_-^{-1}(r) f_+(r)$.
\endroster

\smallskip\noindent
The following theorem describes the evolution of the
scattering coordinates $S^f$.

\proclaim {Theorem I} \rom{([TU1])} If $f(t) \in \Dm$ evolves by
the $j$-th flow and $f(0)$ has only simple poles, then $S^{f(t)}$
evolves as follows:
\roster
\item "{a\rom{)}}" $D^{f(t)} = D^{f(0)}$,
\item "{b\rom{)}}" For $z$ in $D^{f(0)}$,
           $V^{f(t)}_z = e^{-\a z^jt}(V_z^{f(0)})$,
\item "{c\rom{)}}" $v^{f(t)}(r) = e^{\a r^jt}v^{f(0)}(r)e^{-\a r^jt}$.  
\endroster
\endproclaim

  We next explain how to recover $f\in \Dm$ from $S^f$.
To do this first write $f=gh$ with $g\in \Dmd$ and $h\in \Dmc$.
Then,
$v^f = f_-^{-1}f_+ = (g_- h_-)^{-1}(g_+h_+) = h_-^{-1}h^+$,
since as we saw above, $g_- = g_+$. It follows from uniqueness
of the Birkhoff decomposition that $v^f$ determines $h_-$ and
$h_+$ and hence $h$. (Recall that $h$ in $\Cx_+$ (respectively
$\Cx_-$) is the unique meromorphic extension of $h_+$ (respectively
$h_-$).) On the other hand, from the poles $z$ of $g$ and the residues
$\pi^f_z$ of $g$ at these poles we can recover $g$ and hence
$f = gh$.

   There is again an explicit formula for ``pure solitons'',
or ``reflectionless potentials'' (i.e., $u \in P_0$ such that
$f^u$ is in $\Dmd$). We will content ourselves here with writing
the formula for the $1$-solitons of NLS, i.e., a single simple pole,
say at $z = r + is$, with residue the projection of $\Cx^2$ onto the
vector $(\sqrt{1 - |b|^2}, b)$, where $b\in \Cx$ with $|b| < 1$.
Then the solution $q(x,t)$ of NLS is:
$${ 4 s b\sqrt{1-|b|^2}e^{(-2irx+(r^2-s^2)t)}
\over
e^{-2(sx + 2rst)} (1-|b|^2) + e^{2(sx + 2rst)} |b|^2 }.$$
(For $n$-soliton formulas, see [FT] for the $\sun(2)$ case  
and [TU2] for the $\sun(n)$ case.)

 Recall that we have
a natural bijection: $P_0 \to \Dm \to \Dm/\Hm$, where the first
arrow is the Scattering Transform, $\Fs$, and the second is the
natural coset projection. Since we have a natural action
of $\Dm$ on its coset space $\Dm/\Hm$, this induces an action
of $\Dm$ on $P_0$, and so the subgroups $\Dmc$ and $\Dmd$
also act on $P_0$. The orbit of $0$ under $\Dmd$ give the
reflectionless potentials or pure solitons, while the orbit of
$0$ under $\Dmc$ gives the potentials without poles.

  We can now at last explain how the notion of B\"acklund
transformation fits into this picture; namely  
the action of the generators $g_{z,\pi}$ of $\Dmd$ on
$P_0$ are just the  classical B\"acklund transformations. Typically they 
add one to the number of solitons in a solution.

\def\Bibliography{%
\def \bookitem //##1//##2//##3//##4//##5//##6//##7//##8//
{\ref%
\key \ignorespaces##1%
\by \ignorespaces##2%
\book \ignorespaces##3%
\publ \ignorespaces##4%
\yr \ignorespaces##6%
\endref}
\def \b{\bookitem}
\def \articleitem //##1//##2//##3//##4//##5//##6//##7//##8//
{\ref%
\key \ignorespaces##1%
\by \ignorespaces##2%
\paper \ignorespaces##3%
\jour \ignorespaces##4%
\vol \ignorespaces##5%
\yr \ignorespaces##6%
\pages \ignorespaces##7%
\endref}
\def \a{\articleitem}
\def \preprintitem //##1//##2//##3//##4//##5//##6//##7//##8//
{\ref%
\key \ignorespaces##1%
\by \ignorespaces##2%
\paper \ignorespaces##3%
\toappear%
\endref}
\def \p{\preprintitem}
}

\Refs

\widestnumber\key{AKNS2}

\Bibliography

\b //AC//Ablowitz, M.J., Clarkson, P.A.//Solitons, non-linear evolution
equations
and inverse scattering//Cambridge Univ. Press////1991//////

\a //AKNS1//Ablowitz, M.J., Kaup, D.J., Newell, A.C. and Segur, H.//Method for
solving the Sine-Gordon equation//Phys. Rev. Lett.//30//1973//1262--1264////

\a //AKNS2//Ablowitz, M.J., Kaup, D.J., Newell, A.C. and Segur, H.//The
inverse scattering
trans\-form---Fourier analysis for nonlinear problems//Stud. Appl.
Math.//53//1974//249--315////

\b //AbM//Abraham, R., Marsden, J.E.//Foundations of Mechanics//
Benjamin/Cummings////1978//////

\a //Ad//Adler, M.//On a trace functional for formal pseudo-differential
operators and the symplectic structure of the Korteweg-de Vries equation//Invent.
Math//50//1979//219--249////

\a //AdM//Adler, M., van Moerbeke, P.//Completely integrable systems,
Euclidean Lie algebras and curves//Adv. Math.//38//1980//267-317////

\b //Ar//Arnold, V.I.//Mathematical Methods of Classical
Mechanics//Springer-Verlag////1978//////

\b //AA//Arnold, V.I., Avez, A.//Ergodic Problems of Classical Mechanics//
W. A. Benjamin, Inc., New York////1968//////

\b //Au//Audin, M.//Spinning Tops//Cambridge Univ. Press////1996//////

\a //BC1//Beals, R., Coifman, R.R.//Scattering and inverse scattering for
first order systems//
Commun. Pure Appl. Math.//37//1984//39--90////

\a //BC2//Beals, R., Coifman, R.R.//Inverse scattering and evolution
equations//Commun. Pure Appl. Math.//
38//1985//29--42////

\a //BC3//Beals, R., Coifman, R.R.//Linear spectral problems, non-linear
equations and the $\bar\partial$-method//
Inverse Problems//5//1989//87--130////

\a //BS//Beals, R., Sattinger. D.H.//On the complete integrability of
complete integrable systems//
Commun. Math. Phys.//138//1991//409--436////

\a //BDZ//Beals, R., Deift, P., Zhou, X.//The inverse scattering transform
on the line////////7--32////

\a //Bi//Birkhoff, G.D.//Proof of the Ergodic Theorem//Proc. Nat. Acad. Sci.
USA//17//1931//650--660////

\a //BS//Bona, J.L. and Smith, R.//The Initial-Value Problem for the
Korteveg-de Vries Equation//
Philos. Trans. Royal Soc. London, Series A//278//1975//555--604////

\a //Bu//Budagov, A.S.//A completely integrable model of classical field theory
with nontrivial particle
interaction in two-dimensional space-time. Questions in quantum field theory
and statistical physocs//Zap. Nauchn. Sem. Leningrad. Otdel. Mat. Inst.
Steklov (LOMI)//77//1978//24--56, 229//(in Russian)//

\b //BuC//Bullough, R.K., Caudrey, P.J.//Solitons//Topics in Current
Physics, vol. 117, Springer-Verlag////
1980//////

\b //Da//G. Darboux//Le\c cons sur la th\'eorie g\'en\'erale des surfaces//
Chelsea////1972////3rd edition//

\a //DaR//Da Rios////Rend. Circ. Mat. Palermo//22//1906//117-135////

\b //DJ//Drazin, P.G., Johnson, R.S.//Solitons: an introduction//Cambridge
Univ. Press////1989//////

\a //Dr//Drinfel'd, V.G.//Hamiltonian structures on Lie groups, Lie
bialgebras and the geometrical
meaning of Yang-Baxter equations//Sov. Math. Doklady//27//1983//68--71////

\a //DS//Drinfel'd, V.G., and Sokolov, V.V.//Equations of Korteweg-de Vries
type and simple Lie algebras//
Dokl. Akad. Nauk SSSR//258//1981//11--16//(Trans. as Soviet Math. Dokl. 23,
457--462//

\a//Fe//Fermi, E.//Beweis dass ein mechanisches Normalsysteme im Allgemeinen
quasi-ergodisch ist//Phys, Zeit.//
24//1923//261--265////

\b //FT//Faddeev, L.D., Takhtajan, L.A.//Hamiltonian Methods in the theory
of Solitons//\break
Springer-Verlag////1987//////

\b //FPU//Fermi, E., Pasta, J., Ulam, S.//Studies of Nonlinear Problems.
I//in Nonlinear Wave Motion, Lectures in Applied Math., vol. 15,
Amer. Math. Soc.////1974//143--155////

\a //FNR1//Flaschka, H., Newell, A.C., Ratiu, T.//Kac-Moody Lie algebras
and soliton equations,
II. Lax equations associated with $A^{(1)}_1$//Physica//9D//1983//303--323////

\a //FNR2//Flaschka, H., Newell, A.C., Ratiu, T.//Kac-Moody Lie algebras
and soliton equations,
IV. Lax equations associated with $A^{(1)}_1$//Physica//9D//1983//333--345////

\a //FRS//Frenkel, I.E., Reiman, A.G., Semenov-Tian-Shansky, M.A.//Graded
Lie algebras and completely
integrable dynamical systems//Soviet Math. Dokl.//20//1979//811--814////

\a //G//Gardner, C.S.//The Korteweg-de Vries Equation as a Hamiltonian
system//J. Math. Physics//12//
1971//1548--1551////

\a //GGKM//Gardner, C.S., Greene, J.M., Kruskal, M.D., Miura, R.M.//Method
for solving the
Korteweg-de Vries equation//Physics Rev. Lett.//19//1967//1095--1097////

\a //GDi//Gel'fand, I.M., Dikii, L.A.//Fractional Powers of Operators and
Hamiltonian Systems//
Funkcional'nyi Analiz i ego Prilozhenija//10//1976//////

\a //GDo//Gel'fand, I.M., Dorfman, I. Ya//Hamiltonian operators and
algebraic structures related to them//
Functional Anal. Appl.//13//1979//248--261////

\a //GL//Gel'fand, I.M., Levitan, B. M.//On the determination of a differential
equation from its spectral function//
Izv. Akad. Nauk SSSR Ser. Mat.//15//1951//309--366////

\a//Ha//Hasimoto, H.//Motion of a vortex filament and its relation to elastic
//J. Phys. Soc. Japan//31//1971//293--295////

\b //HaK//Hasimoto, H., Kodama, Y.//Solitons in Optical Communications//
Clarendon Press, Oxford////1995//////

\a //Ka1//Kato, T.//On the Cauchy Problem for the (Generalized)
Korteweg-de Vries Equation//
Studies in Applied Math., Adv. in Math. Supp. Stud.//8//1983//93--128////

\b //Ka2//Kato, T.//Quasi-linear equations of evolution, with applications to
 partial differential equations//Lecture Notes in Math., vol. 448, 
 Springer-Verlag, Berlin and New York////
1988//27--50////

\a //KdV//Korteweg, D.J., de Vries, G.//On the change of form of long waves
advancing in a rectangular canal,
and on a new type of long stationary waves//Philos. Mag. Ser.
5//39//1895//422--443////

\a //Kos//Kostant, B.//The solution to a generalized Toda lattice and
representation theory//
Adv. Math.//34//1979//195--338////

\a //KM//Kay, B., Moses, H.E.//The determination of the scattering potential from
the spectral measure function, III//Nuovo Cim.//3//1956//276--304////

\b //KS//Klein, F., Sommerfeld A.//Theorie des Kreisels//Teubner,
Liepzig////1897//////

\b //L//Lamb, G.L., Jr.//Elements of Soliton Theory//
John Wiley \& Sons, New York////1980//////

\a //La1//Lax, P.D.//Integrals of nonlinear equations of evolution and
solitary waves//Comm. Pure. Appl. Math.//
31//1968//467--490////

\b //La2//Lax, P.D.//Periodic Solutions of the KdV Equations//
in Nonlinear Wave Motion, Lectures in
Applied Math., vol. 15, Amer. Math. Soc.////1974//85--96////

\b //La3//Lax, P.D.//Outline of a theory of the KdV equation//in Recent
Mathematical Methods in Nonlinear Wave Propogation, Lecture Notes in Math.,
 vol. 1640,  Springer-Verlag, Berlin and New York////1996//70--102////

\b //LA//Luther, G.G., Alber, M.S.//Nonlinear Waves, Nonlinear Optics, and
Your Communications Future//in Nonlinear Science Today, Springer-Verlag 
New York,Inc.////1997//////

\a //M//Marchenko,V.A.//On the reconstruction of the potential energy from
phases of the scattered waves//
Dokl. Akad. Nauk SSSR //104//1955//695--698////

\b //N//Newell, A.C.//Solitons in Mathematics and Physics//SIAM, CBMS-NSF
vol. 48////1985//////

\b //NMPZ//Novikov, S., Manakov, S., Pitaevskii, L.B., Zakharov,
V.E.//Theory of Solitons//\break
Plenum, New York////1984//////

\a //OU//Oxtoby, J.C., Ulam, S.M.//Measure Preserving Homeomorphisms and
Metric Transitivity//Annals of Math.//
42//1941//874--920////

\b //PT//Palais R.S., and Terng C.L.
//Critical Point Theory and Submanifold Geometry//Lecture Notes in Math., vol. 1353, 
Springer-Verlag, Berlin and New York////1988//////

\b //Pe//Perelomov, A.M.//Integrable Systems of Classical Mechanics and Lie
Algebras//
Birk\-h\"auser Verlag, Basel////1990//////

\b //PrS//Pressley, A. and Segal, G. B.//Loop Groups//Oxford Science Publ.,
Clarendon Press, Oxford////1986//////

\a //RS//Reyman, A.G., Semonov-Tian-Shansky//Current algebras and
non-linear partial differential equations//
Sov. Math., Dokl.//21//1980//630--634////

\a //Ri//Rica, R.L.//Rediscovery of the Da Rios Equation
//Nature//352//1991//561--562////

\b //Ru//Russell, J.S.//Report on Waves//14th Mtg. of the British Assoc.
for the Advance. of Science,
John Murray, London, pp. 311--390 + 57 plates////1844//////

\a //Sa//Sattinger, D.H.//Hamiltonian hierarchies on semi-simple Lie
algebras//Stud. Appl. Math.//
72//1984//65--86////

\a //SW//Segal, G., Wilson, G.//Loop groups and equations of KdV
type//Publ. Math. IHES//61//1985//5--65////

\a //Se1//Semenov-Tian-Shansky, M.A.//Dressing transformations and Poisson
group actions//Publ. RIMS Kyoto Univ.//
21//1985//1237--1260////

\a //Se2//Semenov-Tian-Shansky, M.A.//Classical r-matrices, Lax equations,
Poisson Lie groups, and dressing
transformations// Lecture Notes in Physics,
Springer-Verlag//280//1986//174--214////

\a //Sh//Shabat, A.B.//An inverse scattering problem//Diff.
Equ.//15//1980//1299--1307////

\a //St//Strang, G.//On the Construction and Comparison of Difference Schemes//SIAM
J. Numerical Analysis//5//1968// 506--517////

\a //Sy//Symes, W.W.//Systems of Toda type, Inverse spectral problems, and
representation theory//
Inventiones Math.//59//1980//13--51////

\a //Ta//Tappert, F.//Numerical Solutions of the Korteweg-de Vries
Equations and its Generalizations by the
Split-Step Fourier Method//in Nonlinear Wave Motion, Lectures in
Applied Math., vol. 15, Amer. Math. Soc.////1974//215--216////

\a //Te1//Terng, C.L.//A higher dimensional generalization of the
Sine-Gordon equation and its soliton theory//
Ann. Math.//111//1980//491--510////

\a //Te2//Terng, C.L.//Soliton equations and differential geometry//J.
Differential Geometry//45//1996//407--445////

\p //TU1//Terng, C.L., Uhlenbeck, K.//Poisson Actions and Scattering Theory
for Integrable Systems, dg-ga/9707004//////////preprint//

\p //TU2//Terng, C.L., Uhlenbeck, K.//B\"acklund transformations and loop group actions
//////////preprint//

\a //U1//Uhlenbeck, K.//Harmonic maps into Lie group (classical solutions
of the chiral model)//
J. Differential Geometry//30//1989//1-50////

\a //U2//Uhlenbeck, K.//On the connection between harmonic maps and the
self-dual Yang-Mills and the
Sine-Gordon equations//Geometry \& Physics//2//1993//////

\b //Ul//Ulam, S. M.//Adventures of a Mathematician//Univ. of Calif.
Press////1991//////

\a //Wa//Wadati, M.//The modified Korteweg-de Vries equation//J. Phys. Soc.
Japan//34//1973//380--384////

\a //Wi//Wilson, G.//The modified Lax equations and two dimensional
Toda lattice equations associated with simple
Lie algebras//Ergodic Theory and Dynamical Systems I//30//1981//361--380////

\a //ZK//Zabusky, N.J., Kruskal, M.D.//Interaction of solitons in a
collisionless plasma and the recurrence
of initial states//Physics Rev. Lett.//15//1965//240--243////

\a //ZF//Zakharov, V.E., Faddeev, L.D.//Korteweg-de Vries equation, a
completely integrable Hamiltonian system//
Func. Anal. Appl.//5//1971//280--287////

\a //ZMa1//Zakharov, V.E., Manakov, S.V.//On resonant interaction of wave
packets in non-linear media//JETP Letters//
18//1973//243--247////

\a //ZMa2//Zakharov, V.E., Manakov, S.V.//The theory of resonant
interaction of wave packets in non-linear media//
Sov. Phys. JETP//42//1975//842--850////

\a //ZMi1//Zakharov, V.E., Mikhailov, A.V.//Example of
nontrivial interaction of solitons in two-dimensional classical
field theory// JETP Letters//27//1978//42--46////

\a //ZMi2//Zakharov, V.E., Mikhailov, A.V.//Relativistically invariant
two-dimensional models of field theory
which are integrable by means of the inverse scattering problem
method//Soviet Physics JETP//47//1978//
1017--1027////

\a //ZS//Zakharov, V.E., Shabat, A.B.//Exact theory of two-dimensional
self-focusing and one-dimensional self-modulation of waves
in nonlinear media//Sov. Phys. JETP//34//1972//62--69////

\endRefs

\enddocument

\bye